\documentclass[numberwithinsect,cleveref,thm-restate,USenglish,nameinlink]{lipics-v2021}

\hypersetup{colorlinks=true,citecolor=purple,linkcolor=purple,hypertexnames=false}
\usepackage{todonotes}
\usepackage{tensor}
\usepackage{mathtools}

\usepackage{tikz}
\usetikzlibrary{positioning}

\nolinenumbers

\makeatletter
\let\c@author\relax
\makeatother
\usepackage[sortcites=true,firstinits=true,maxnames=10]{biblatex}
\usepackage{xparse}
\makeatletter
\NewDocumentCommand{\xsideset}{mmme{_^}}{
  \mathop{
    \settowidth{\dimen0}{$\m@th\displaystyle#3$}
    \dimen0=.5\dimen0
    \settowidth{\dimen2}{$
      \m@th\displaystyle#3
      \IfValueT{#4}{_{#4}}
      \IfValueT{#5}{^{#5}}
    $}
    \dimen2=.5\dimen2
    \advance\dimen2 -\dimen0
    \sbox6{\scriptspace\z@$\displaystyle{\vphantom{#3}}#1$}
    \sbox8{\scriptspace\z@$\displaystyle{\vphantom{#3}}#2$}
    \ifdim\wd6>\dimen2 \kern\dimexpr\wd6-\dimen2\relax\fi
    {
     \mathop{\llap{\copy6}{\displaystyle#3}\rlap{\copy8}}\limits
     \IfValueT{#4}{_{#4}}
     \IfValueT{#5}{^{#5}}
    }
    \ifdim\wd8>\dimen2 \kern\dimexpr\wd8-\dimen2\relax\fi
  }
}
\makeatother

\newtheorem{defn}[theorem]{Definition}

\title{A characterization of wreath products where knapsack is decidable}

\author{Pascal Bergstr\"{a}{\ss}er}%
{Fachbereich Informatik, Technische Universit\"{a}t Kaiserslautern, Germany}%
{}%
{https://orcid.org/0000-0002-4681-2149}%
{}%

\author{Moses Ganardi}%
{Max Planck Institute for Software Systems (MPI-SWS), Kaiserslautern, Germany}%
{}%
{https://orcid.org/0000-0002-0775-7781}%
{}%

\author{Georg Zetzsche}%
{Max Planck Institute for Software Systems (MPI-SWS), Kaiserslautern, Germany}%
{}%
{https://orcid.org/0000-0002-6421-4388}%
{}%

\authorrunning{Pascal Bergstr\"{a}{\ss}er, Moses Ganardi, and Georg Zetzsche}
\Copyright{Pascal Bergstr\"{a}{\ss}er, Moses Ganardi, Georg Zetzsche}



\newcommand{\KP}{\mathsf{KP}}
\newcommand{\IKPpm}{\mathsf{KP}^{\pm}}
\newcommand{\IKPp}{\mathsf{KP}^{+}}

\newcommand{\BS}{\mathsf{BS}}

\newcommand{\ExpEq}{\mathsf{ExpEq}}
\newcommand{\N}{\mathbb{N}}
\newcommand{\Z}{\mathbb{Z}}

\newcommand{\NP}{\mathsf{NP}}

\newcommand{\HKP}{\mathsf{HKP}}
\newcommand{\supp}{\mathsf{supp}}
\newcommand{\pc}{\mathsf{pc}}
\newcommand{\lo}[2]{\tensor*[^{#1}]{{#2}}{}}
\newcommand{\loi}[4]{\tensor*[^{#1}]{{#2}}{^{#3}_{#4}}}

\newcommand{\Q}{\mathbb{Q}}
\newcommand{\MKP}{\mathsf{MKP}}
\newcommand{\SAT}{\mathsf{SAT}}
\newcommand{\GL}{\mathsf{GL}}
\newcommand{\Th}{\mathsf{Th}}
\newcommand{\ord}{\mathsf{ord}}
\newcommand{\sol}{\mathsf{sol}}

\newcommand\tostar[1][w]{\xrightarrow{#1}\mathrel{\vphantom{\to}^*}}

\newcommand{\xsubparagraph}[1]{\subparagraph{#1}}

\newcommand{\gz}[1]{}

\begin{document}

\maketitle

\begin{abstract}
  The knapsack problem for groups was introduced by Miasnikov,
  Nikolaev, and Ushakov. It is defined for each finitely generated
  group $G$ and takes as input group elements
  $g_1,\ldots,g_n,g\in G$ and asks whether there are
  $x_1,\ldots,x_n\ge 0$ with $g_1^{x_1}\cdots
  g_n^{x_n}=g$. We study the knapsack problem for wreath products
  $G\wr H$ of groups $G$ and $H$.
  
  Our main result is a characterization of those wreath products
  $G\wr H$ for which the knapsack problem is decidable. The
  characterization is in terms of decidability properties of the
  indiviual factors $G$ and $H$.  To this end, we introduce two
  decision problems, the \emph{intersection knapsack problem} and its
  restriction, the \emph{positive intersection knapsack problem}.

  Moreover, we apply our main result to $H_3(\Z)$, the discrete
  Heisenberg group, and to Baumslag-Solitar groups $\BS(1,q)$ for
  $q\ge 1$. First, we show that the knapsack problem is undecidable
  for $G\wr H_3(\Z)$ for any $G\ne 1$. This implies that for $G\ne 1$
  and for infinite and virtually nilpotent groups $H$, the knapsack
  problem for $G\wr H$ is decidable if and only if $H$ is virtually
  abelian and solvability of systems of exponent equations is
  decidable for $G$. Second, we show that the knapsack problem is
  decidable for $G\wr\BS(1,q)$ if and only if solvability of systems
  of exponent equations is decidable for $G$.
\end{abstract}

\newpage
\section{Introduction}

\xsubparagraph{The knapsack problem} The knapsack problem is a decision
problem for groups that was introduced by Miasnikov, Nikolaev, and
Ushakov~\cite{MiNiUs14}.  If $G$ is a finitely generated group, then
the knapsack problem for $G$, denoted $\KP(G)$, takes group elements
$g_1,\ldots,g_n,g\in G$ as input
(as words over the generators) and
it asks whether there are natural numbers $x_1,\ldots,x_n\ge 0$ such
that $g_1^{x_1}\cdots g_n^{x_n}=g$. Since its introduction, a
significant amount of attention has been devoted to understanding for
which groups the problem is decidable and what the resulting
complexity
is~\cite{Loh19hyp,LohreyZ18,GanardiKLZ18,MiTr17,KoenigLohreyZetzsche2015a,FrenkelNU15,LohreyZ20,FigeliusGLZ20}.
For matrix semigroups, the 
knapsack problem has been studied implicitly
by Bell, Halava, Harju, Karhumäki, and Potapov \cite{BellHHKP08},
Bell, Potapov, and Semukhin \cite{BellPS19}, and for commuting
matrices by Babai, Beals, Cai, Ivanyos, and Luks~\cite{BabaiBCIL96}.

There are many groups for which knapsack has been shown decidable. For
example, knapsack is decidable for virtually special
groups~\cite[Theorem 3.1]{LohreyZ18}, co-context-free
groups~\cite[Theorem 8.1]{KoenigLohreyZetzsche2015a}, hyperbolic
groups~\cite[Theorem 6.1]{MiNiUs14}, the discrete Heisenberg
group~\cite[Theorem 6.8]{KoenigLohreyZetzsche2015a}, and
Baumslag-Solitar groups $\BS(p,q)$ for co-prime $p,q>1$~\cite[Theorem
2]{dudkin2018knapsack} and for $p=1$~\cite[Theorem
4.1]{LohreyZ20}. Moreover, the class of groups where knapsack is
decidable is closed under free products with
amalgamation~\cite[Theorem 14]{LohreyZetzsche2016a} and HNN
extensions~\cite[Theorem 13]{LohreyZetzsche2016a} over finite
identified subgroups. On the other hand, there are nilpotent groups
for which knapsack is undecidable~\cite[Theorem
6.5]{KoenigLohreyZetzsche2015a}.

\xsubparagraph{Wreath products} A prominent construction in group
theory and semigroup theory is the wreath product $G\wr H$ of two
groups $G$ and $H$.  Wreath products are important algorithmically,
because the Magnus embedding theorem~\cite[Lemma]{Mag39} states that
for any free group $F$ of rank $r$ and a normal subgroup $N$ of $F$,
one can find $F/[N,N]$ as a subgroup of $\Z^r\wr (F/N)$, where $[N,N]$
is the commutator subgroup of $N$.  This has been used by several
authors to obtain algorithms for groups of the form $F/[N,N]$, and in
particular free solvable groups. Examples include the word problem
(folklore, see~\cite{gul2017magnus}), the conjugacy
problem~\cite{matthews1966conjugacy,remeslennikov1970some,gul2017magnus,MVW19},
the power problem~\cite{gul2017magnus}, and the knapsack
problem~\cite{FigeliusGLZ20,GanardiKLZ18}.

For groups $G$ and $H$,
their wreath product $G\wr H$ can be roughly described as follows. An
element of $G\wr H$ consists of (i)~a labeling, which maps each
element of $H$ to an element of $G$ and (ii)~an element of $H$, called
the \emph{cursor}. Here, the labeling has finite support, meaning all
but finitely many elements of $H$ are mapped to the identity of
$G$. Moreover, each element of $G\wr H$ can be written as a product of
elements from $G$ and from $H$. Multiplying an element $g\in G$ will
multiply $g$ to the label of the current cursor position. Multiplying
an element $h\in H$ will move the cursor by multiplying $h$.

Understanding the knapsack problem for wreath products is challenging
for two reasons. First, the path that the expression
$g_1^{x_1}\cdots g_n^{x_n}g^{-1}$ takes through the group $H$ can have
complicated interactions with itself: The product can place elements
of $G$ at (an \emph{a priori} unbounded number of) positions $h\in H$
that are later revisited. At the end of the path, each position of $H$
must carry the identity of $G$ so as to obtain
$g_1^{x_1}\cdots g_k^{x_k}g^{-1}=1$. The second reason is that the
groups $G$ and $H$ play rather different roles: \emph{A priori}, for each
group $G$ the class of all $H$ with decidable $\KP(G\wr H)$ could be
different, resulting in a plethora of cases.

Decidability of the knapsack problem for wreath products has been
studied by Ganardi, König, Lohrey, and Zetzsche~\cite{GanardiKLZ18}. They focus on the case
that $H$ is knapsack-semilinear, which means that the solution sets of
equations $g_1^{x_1}\cdots g_n^{x_n}=g$ are (effectively) semilinear.
A set $S \subseteq \N^n$ is {\em semilinear} if it is a finite union
of {\em linear} sets $\{ u_0 + \lambda_1 u_1 + \dots + \lambda_k u_k \mid \lambda_1, \dots, \lambda_k \in \N  \}$
for some vectors $u_0, \dots, u_k \in \N^n$.
Under this assumption, they show that $\KP(G\wr H)$ is decidable if
and only if solvability of systems of exponent equations is decidable
for $G$~\cite[Theorem 5.3]{GanardiKLZ18}.  Here, an exponent equation
is one of the form $g_1^{x_1}\cdots g_n^{x_n}=g$, where variables
$x_i$ are allowed to repeat. The problem of solvability of systems of
exponent equations is denoted $\ExpEq(G)$. Moreover, it is shown there
that for some number $\ell\in\N$, knapsack is undecidable for
$G\wr (H_3(\Z)\times\Z^\ell)$, where $H_3(\Z)$ denotes the discrete
Heisenberg group and $G$ is any non-trivial group~\cite[Theorem
5.2]{GanardiKLZ18}. Since $\KP(H_3(\Z)\times\Z^\ell)$ is decidable for
any $\ell\ge 0$~\cite[Theorem 6.8]{KoenigLohreyZetzsche2015a}, this
implies that wreath products do not preserve decidability of knapsack
in general.  However, apart from the latter undecidability result,
little is known about wreath products $G\wr H$ where $H$ is not
knapsack-semilinear.  As notable examples of this, knapsack is
decidable for solvable Baumslag-Solitar
groups~$\BS(1,q)$~\cite[Theorem 4.1]{LohreyZ20} and for the discrete
Heisenberg group $H_3(\Z)$~\cite[Theorem
6.8]{KoenigLohreyZetzsche2015a}, but it is not known for which $G$ the
knapsack problem is decidable for $G\wr H_3(\Z)$ or for $G\wr \BS(1,q)$.

The only other paper which studies the knapsack problem over wreath products is \cite{FigeliusGLZ20}.
It is concerned with complexity results (for knapsack-semilinear groups)
whereas in this paper we are concerned with decidability results.

\xsubparagraph{Contribution} Our main result is a characterization of
the groups $G$ and $H$ for which $\KP(G\wr H)$ is
decidable. Specifically, we introduce two problems, \emph{intersection
  knapsack} $\IKPpm(H)$ and the variant \emph{positive intersection
  knapsack} $\IKPp(H)$ and show the following. Let $G$ and $H$ be
finitely generated, with $G$ non-trivial and $H$ infinite. Then
knapsack for $G\wr H$ is decidable if and only if $\ExpEq(G)$ is
decidable and either (i)~$G$ is abelian and $\IKPp(H)$ is decidable or
(ii)~$G$ is not abelian and $\IKPpm(H)$ is decidable. Note that the
case of finite $H$ is not interesting: For $|H|=m$, $\KP(G\wr H)$ is
equivalent to $\KP(G^{m})$ (see \cref{sec:results}).

Thus, our
result relieves us from considering every pair $(G,H)$ of groups and
allows us to study the factors separately.
It is not hard to see that decidability of $\ExpEq(G)$ is necessary
for decidability of $\KP(G\wr H)$ if $H$ is infinite. It is
surprising that the only other property of $G$ that is relevant for
decidability of $\KP(G\wr H)$ is whether $G$ is abelian or not. This
is in contrast to the effect of other structural properties of $G$ on
the complexity of $\KP(G\wr\Z)$: If $G\ne 1$ is a finite nilpotent
group, then $\KP(G\wr\Z)$ is $\NP$-complete~\cite[Theorem
2]{FigeliusGLZ20}, whereas for finite and non-solvable $G$, the
problem $\KP(G\wr\Z)$ is $\Sigma_2^p$-complete~\cite[Corollary
25]{FigeliusGLZ20}.

\xsubparagraph{Applications} We also obtain two applications. First,
we deduce that $\KP(G\wr H_3(\Z))$ is undecidable for every
$G\ne 1$. This implies that if $G\ne 1$ and $H$ is virtually nilpotent
and infinite, then $\KP(G\wr H)$ is decidable if and only if $H$ is
virtually abelian and $\ExpEq(G)$ is decidable.  Moreover, we show
that $\KP(G\wr \BS(1,q))$ is decidable if and only if $\ExpEq(G)$ is.

\xsubparagraph{Ingredients} For the ``if'' direction of our main
result, we reduce $\KP(G\wr H)$ to $\ExpEq(G)$ and $\IKPpm(H)$
(respectively $\IKPp(H)$) using extensions of techniques used by
Figelius, Ganardi, Lohrey, and Zetzsche~\cite{FigeliusGLZ20}. Roughly speaking, the problem
$\IKPpm(H)$ takes as input an expression
$h_0g_1^{x_1}h_1\cdots g_n^{x_n}h_n$ and looks for numbers
$x_1,\ldots,x_n\ge 0$ such that the walk defined by the product
$h_0g_1^{x_1}h_1\cdots g_n^{x_n}h_n$ meets specified constraints about
self-intersections. Such a constraint can be either (i)~a \emph{loop
  constraint}, meaning the walk visits the same point after two
specified factors or (ii)~a \emph{disjointness constraint} saying that
the $(x_i+1)$-many points visited when multiplying $g_i^{x_i}$ do not
intersect the $(x_j+1)$-many points visited while multiplying $g_j^{x_j}$.

The ``only if'' reductions in our main result involve substantially
new ideas.  The challenge is to guarantee that the constructed
instances of $\KP(G\wr H)$ will leave an element $\ne 1$ somewhere, as
soon as any constraint is violated. In particular, the loop
constraints have to be checked independently of the disjointness
constraints. Moreover, if several constraints are violated, the
resulting elements $\ne 1$ should not cancel each other.
Furthermore, this has to be achieved despite almost no information on the structure of $G$ and $H$.
This requires
an intricate construction that uses various patterns in the Cayley
graph of $H$ for which we show that only very specific arrangements
permit cancellation. To this end, we introduce the notion of \emph{periodic complexity},
which measures how many periodic sequences are needed to cancel out a sequence of elements of a group.
Roughly speaking, for the loop constraints we use patterns of high periodic complexity,
whereas for the disjointness constraints we use patterns with low periodic complexity but many large gaps.
This ensures that the disjointness patterns cannot cancel the loop patterns or vice versa.

\section{Preliminaries}
\xsubparagraph{Knapsack problems}
For a group $G$ and a subset $S \subseteq G$ we write $S^*$
for the submonoid generated by $S$, i.e. the set of products of elements from $S$.
Let $G$ be a group with a finite {\em (monoid) generating set} $\Sigma \subseteq G$, i.e. $G = \Sigma^*$.
Such groups are called {\em finitely generated}.
An {\em exponent expression} over $G$ is an expression $E = e_1 \dots e_n$
consisting of {\em atoms} $e_i$
where each atom $e_i$ is either a {\em constant} $e_i = g_i \in G$
or a {\em power} $e_i = g_i^{x_i}$ for some $g_i \in G$ and variable $x_i$.
Here the group elements $g_i$ are given as words over $\Sigma$.
We write $\gamma(e_i) = g_i$ for the constant or the base of the power.
Furthermore let $P_E \subseteq [1,n]$ be the set of indices of the powers in $E$
and $Q_E = [1,n] \setminus P_E$ be the set of indices of the constants in $E$.
If $\nu \in \mathbb{N}^X$ is a valuation of the variables $X$ that occur in $E$, then for each $i\in [1,n]$, we define
$\nu(e_i)=\gamma(e_i)^{\nu(x_i)}$ if $i\in P_E$; and $\nu(e_i)=e_i$ if $i\in Q_E$.
Moreover, $\nu(E) := \nu(e_1) \cdots \nu(e_n)$ and the set of $G$-{\em solutions} of $E$ as
$\sol_G(E) := \{\nu \in \N^X \mid \nu(E) = 1\}$.

For a group $G$, the problem of \emph{solvability of exponent equations} $\ExpEq(G)$ is defined as:
\begin{description} 
\item[Given] a finite list of exponent expression $E_1,\dots,E_k$ over $G$.
\item[Question] Is $\bigcap_{i=1}^k \sol_G(E_i)$ non-empty?
\end{description}

An exponent expression is called a {\em knapsack expression}
if all variables occur at most once.
The {\em knapsack problem} $\KP(G)$ over $G$ is defined as follows:
\begin{description}
\item[Given] a knapsack expression $E$ over $G$.
\item[Question] Is there a valuation $\nu$ such that $\nu(E) = 1$?
\end{description}
The definition from~\cite{MiNiUs14} asks whether
$g_1^{x_1}\cdots g_n^{x_n}=g$ has a solution for given
$g_1,\ldots,g_n,g\in G$. The two versions are inter-reducible in
polynomial time~\cite[Proposition 7.1]{KoenigLohreyZetzsche2015a}.

\xsubparagraph{Wreath products}
Let $G$ and $H$ be groups. Consider the direct sum $K = \bigoplus_{h
  \in H} G_h$, where $G_h$ is a copy of $G$. We view $K$ as the set $G^{(H)}$ of
all mappings $f\colon H\to G$ such that $\supp(f) := \{h\in H \mid f(h)\ne
1\}$ is finite, together with pointwise multiplication as the group
operation.  The set $\supp(f)\subseteq H$ is called the
\emph{support} of $f$. The group $H$ has a natural left action on
$G^{(H)}$ given by $\lo{h}{f}(a) = f(h^{-1}a)$, where $f \in G^{(H)}$ and
$h, a \in H$.  The corresponding semidirect product $G^{(H)} \rtimes
H$ is the (restricted) \emph{wreath product} $G \wr H$.  In other words:
\begin{itemize}
\item
Elements of $G \wr H$ are pairs $(f,h)$, where $h \in H$ and
$f \in G^{(H)}$.
\item
The multiplication in $G \wr H$ is defined as follows:
Let $(f_1,h_1), (f_2,h_2) \in G \wr H$. Then
$(f_1,h_1)(f_2,h_2) = (f, h_1h_2)$, where
$f(a) = f_1(a)f_2(h_1^{-1}a)$.
\end{itemize}
There are canonical mappings
$\sigma \colon G \wr H \to H$ with $\sigma(f,h) = h$ and  
$\tau \colon G \wr H \to G^{(H)}$ with $\tau(f,h) = f$ for $f\in G^{(H)}$, $h\in H$.
In other words: $g = (\tau(g), \sigma(g))$ for $g \in G \wr H$.
Note that $\sigma$ is a homomorphism whereas $\tau$ is in general not a homomorphism.
Throughout this paper, the letters $\sigma$ and $\tau$ will have the above meaning
(the groups $G,H$ will be always clear from the context).
We also define $\supp(g) = \supp(\tau(g))$ for all $g \in G \wr H$.

The following intuition might be helpful:
An element $(f,h) \in G\wr H$ can be thought of
as a finite multiset of elements of $G \setminus\{1_G\}$ that are sitting at certain
elements of $H$ (the mapping $f$) together with the distinguished
element $h \in H$, which can be thought of as a \emph{cursor}
moving in $H$.
We can compute the product $(f_1,h_1) (f_2,h_2)$
as follows: First, we shift the finite collection of $G$-elements that
corresponds to the mapping $f_2$ by $h_1$: If the element $g \in G\setminus\{1_G\}$ is
sitting at $a \in H$ (i.e., $f_2(a)=g$), then we remove $g$ from $a$ and
put it to the new location $h_1a \in H$. This new collection
corresponds to the mapping $f'_2 \colon  a \mapsto f_2(h_1^{-1}a)$.
After this shift, we multiply the two collections of $G$-elements
pointwise: If $g_1\in G$ and $g_2\in G$ are sitting at $a \in H$
(i.e., $f_1(a)=g_1$ and $f'_2(a)=g_2$), then we put 
$g_1g_2$ into the location $a$. The new distinguished
$H$-element (the new cursor position) becomes $h_1 h_2$.

Clearly, $H$ is a subgroup of $G \wr H$. We also regard $G$ as a subgroup of $G \wr H$ by identifying 
$G$ with the set of all $f \in G^{(H)}$ with $\supp(f) \subseteq \{1\}$. This copy of $G$ together 
with $H$ generates $G \wr H$.   In particular, if $G = \langle \Sigma \rangle$ and $H = \langle \Gamma \rangle$ 
with $\Sigma \cap \Gamma = \emptyset$ then $G \wr H$ is generated by $\Sigma \cup \Gamma$.
With these embeddings, $GH$ is the set of $(f,h)\in G\wr H$ with $\supp(f)\subseteq\{1\}$ and $h\in H$.

\xsubparagraph{Groups}
Our applications will involve two well-known types of groups: the {\em discrete Heisenberg group} $H_3(\Z)$, which consists of the %
matrices $\left(\begin{smallmatrix} 1 & a & c \\ 0 & 1 & b \\ 0 & 0 & 1
  \end{smallmatrix}\right)$
with $a,b,c \in \Z$, and the {\em Baumslag-Solitar
  groups}~\cite{baumslag1962some} $\BS(p, q)$ for $p,q \in \N$, where
$\BS(p,q) = \langle a,t \mid ta^pt^{-1} = a^q \rangle$.

A subgroup $H$
of $G$ is called \emph{finite-index} if there are finitely many cosets
$gH$. If $ab=ba$ for every $a,b\in G$, then $G$ is \emph{abelian}.
A group has a property \emph{virtually} if it has a
finite-index subgroup $H$ with that property. For example, a group is
virtually abelian if it has a finite-index abelian subgroup.
For two elements $a,b\in G$, we write $[a,b]=aba^{-1}b^{-1}$ and call
this the \emph{commutator} of $a,b$.
If $A,B$ are subgroups of $G$, then $[A,B]$ is the subgroup
generated by all $[a,b]$ with $a\in A$ and $b\in B$. For $g,h\in G$,
we write $\loi{h}{g}{}{}=hgh^{-1}$. In particular, if $g\in G$ and $h\in H$,
then $\loi{h}{g}{}{}$ is the element $(f,1)\in G\wr H$ with $f(h)=g$
and $f(h')=1$ for $h'\ne h$.

\section{Main results}\label{sec:results}
We first introduce the new (positive) intersection knapsack problem.  A
solution to a knapsack expression $E$ describes a walk in the Cayley
graph that starts and ends in the group identity. Whereas the
ordinary knapsack problem only asks for the expression to yield the
identity, our extended version can impose constraints on how this walk
intersects itself.

A {\em walk} over $G$ is a nonempty sequence $\pi = (g_1, \dots, g_n)$ over $G$.
Its support is $\supp(\pi) = \{g_1, \dots, g_n\}$.
It is a {\em loop} if $g_1 = g_n$.
Two walks are {\em disjoint} if their supports are disjoint.
We define a partial concatenation on walks:
If $\pi = (g_1, \dots, g_n)$ and $\rho = (h_1, \dots, h_m)$ with $g_n = h_1$
then $\pi \rho = (g_1, \dots, g_n, h_2, \dots, h_m)$.
A {\em progression} with period $h \in G$ over $G$ is a walk of the form $\pi = (g, gh, gh^2, \dots, gh^\ell)$
for some $g \in G$ and $\ell \ge 0$.
We also call the set $\supp(\pi)$ a progression, whose period may not be unique.
If $h \neq 1$ we also call $\pi$ a {\em ray}.

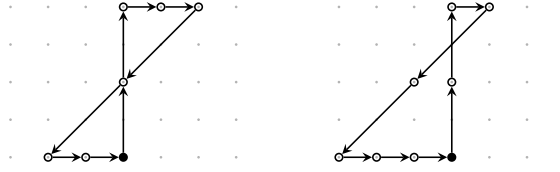
\begin{figure}
\centering
\raisebox{-0.5\height}{
\begin{tikzpicture}[semithick,scale=.5]
   \foreach \x in {-3,...,3} {
        \foreach \y in {0,...,4} {
            \fill[color=black!30] (\x,\y) circle (1pt);
        }
    }
\tikzset{every node/.style={circle, draw, inner sep = 1pt}}
\node (A) [fill, circle, inner sep = 1pt] at (0,0) {};
\node (B) at (0,2) {};
\node (C) at (0,4) {};
\node (C2) at (1,4) {};
\node (D) at (2,4) {};
\node (E2) at (-1,0) {};
\node (E) at (-2,0) {};
\draw[->,>=stealth] (A) edge (B) (B) edge (C) (C) edge (C2) (C2) edge (D) (D) edge (B) (B) edge (E) (E) edge (E2) (E2) edge (A);
\end{tikzpicture}
}\hspace{3em}\raisebox{-0.5\height}{
\begin{tikzpicture}[semithick,scale=.5]
   \foreach \x in {-3,...,3} {
        \foreach \y in {0,...,4} {
            \fill[color=black!30] (\x,\y) circle (1pt);
        }
    }
\tikzset{every node/.style={circle, draw, inner sep = 1pt}}
\node (A) [fill, circle, inner sep = 1pt] at (0,0) {};
\node (B) at (0,2) {};
\node (C) at (0,4) {};
\node (C2) at (1,4) {};
\node (D) at (-1,2) {};
\node (D2) at (-3,0) {};
\node (E2) at (-1,0) {};
\node (E) at (-2,0) {};
\draw[->,>=stealth] (A) edge (B) (B) edge (C) (C) edge (C2) (C2) edge (D) (D) edge (D2) (D2) edge (E) (E) edge (E2) (E2) edge (A);
\end{tikzpicture}
}
\caption{Consider the knapsack equation $g_1^{x_1}g_2^{x_2}g_3^{x_3}g_4^{x_4} = 1$ over $\Z^2$
written multiplicatively, where $g_1 = (0,2)$, $g_2 = (1,0)$, $g_3 = (-2,-2)$ and $g_4 = (1,0)$
and the disjointness condition $D = \{(1,3)\}$.
The solid dot represents the origin $(0,0)$.
The knapsack equation is satisfied by $(x_1,x_2,x_3,x_4) = (2,2,2,2)$ but it violates $D$,
as illustrated on the left.
On the right the solution $(x_1,x_2,x_3,x_4) = (2,1,2,3)$ is depicted, which satisfies $D$.}
\label{fig}
\end{figure}

A {\em factorized walk} is a walk $\pi$ equipped with a {\em factorization} $(\pi_1, \dots, \pi_n)$,
i.e. $\pi = \pi_1 \dots \pi_n$.
One also defines the concatenation of factorized walks in the straightforward fashion.
If $E = e_1 \dots e_n$ is an exponent expression and $\nu$ is a valuation over $E$
we define the factorized walk $\pi_{\nu,E} = \pi_1 \dots \pi_n$ induced by $\nu$ on $E$
where
\[
	\pi_i = \begin{cases}
	(\nu(e_1 \dots e_{i-1}) \, g_i^k)_{0 \le k \le \nu(x_i)}, & \text{if } e_i = g_i^{x_i} \\
	(\nu(e_1 \dots e_{i-1}), \nu(e_1 \dots e_{i-1}) \, g_i), & \text{if } e_i = g_i.
	\end{cases}
\]
The {\em intersection knapsack problem} $\KP^\pm(G)$ over $G$ is defined as follows:
\begin{description}
\item[Given] a knapsack expression $E$ over $G$, a set $L \subseteq [0,n]^2$ of loop constraints,
and a set $D \subseteq [1,n]^2$ of disjointness constraints.
\item[Question] Is there a valuation $\nu$ such that $\nu(E) = 1$
and the factorized walk $\pi_{\nu,E} = \pi_1 \dots \pi_n$ induced by $\nu$ on $E$ satisfies
the following conditions:
\begin{itemize}
\item $\pi_{i+1} \dots \pi_j$ is a loop for every $(i,j) \in L$
\item $\pi_i$ and $\pi_j$ are disjoint for every $(i,j) \in D$.
\end{itemize}
\end{description}
The {\em positive intersection knapsack problem} $\KP^+(G)$ over $G$
is the restriction of $\KP^\pm(G)$ to instances where $D = \emptyset$.
We denote the set of solutions of a $\KP^\pm(G)$-instance (resp. $\KP^+(G)$-instance) $(E,I,D)$ (resp. $(E,I)$) as $\sol_G(E,I,D)$ (resp. $\sol_G(E,I)$).
\Cref{fig} shows an example for the intersection knapsack problem over $\Z^2$.

The following is our main result.
\begin{theorem}\label{main-result}
  Let $G$ and $H$ be f.g.\ groups such that $G$ is non-trivial and $H$
  is infinite. Then $\KP(G\wr H)$ is decidable if and only if
  $\ExpEq(G)$ is decidable and either
  \begin{enumerate}
  \item $G$ is abelian and $\IKPp(H)$ is decidable or
  \item $G$ is not abelian and $\IKPpm(H)$ is decidable.
  \end{enumerate}
\end{theorem}
Here, we assume $H$ to be infinite, because the case of finite $H$ is
not interesting: If $|H|=m$, then $G\wr H$ has $G^m$ as a finite-index
subgroup~\cite[Proposition 1]{LohreySteinbergZetzsche2015a}, meaning
$\KP(G\wr H)$ is decidable if and only if $\KP(G^m)$ is~\cite[Theorem
7.3]{KoenigLohreyZetzsche2015a}.

If $H$ is knapsack-semilinear, it is easy to see that both $\IKPp(H)$
and $\IKPpm(H)$ are decidable via an encoding in Presburger
arithmetic.  Hence, the main decidability result of
\cite{GanardiKLZ18}, saying that for knapsack-semilinear $H$,
$\KP(G\wr H)$ is decidable if and only if $\ExpEq(G)$ is decidable, is
generalized by \cref{main-result}.

\xsubparagraph{Logical version of $\KP^+$ and $\KP^{\pm}$} For our
applications of \cref{main-result}, it is often convenient to use a
formulation of $\KP^+(G)$ and $\KP^\pm(G)$ in terms of logics over an
extended Cayley graph of $G$. The \emph{Cayley graph of $G$} is the
logical structure
$\mathcal{C}(G) = (G,(\xrightarrow{g})_{g \in G})$, with domain $G$
and with the relation $\xrightarrow{g}$ for each\footnote{Customarily, one
  only includes the edge relations $(\xrightarrow{s})_{s \in S}$ for some
  finite generating set $S$ of $G$. We choose $S = G$ to make
  the presentation in the following cleaner.}  $g\in G$, where
$g_1 \xrightarrow{g} g_2$ if and only if $g_1 g = g_2$.  We define the
extension
$\mathcal{C}^+(G) = (G,(\xrightarrow{g})_{g \in G},(\tostar[g])_{g \in
  G})$ where $\tostar[g]$ is the reflexive transitive closure of
$\xrightarrow{g}$.  Finally, we define a further extension
$\mathcal{C}^\pm(G) = (G,(\xrightarrow{g})_{g \in G},(\tostar[g])_{g
  \in G},(\bot_{g,h})_{g,h \in G})$ with \emph{disjointness relations}
$\bot_{g,h}$, which are binary relations on pairs $G^2$: For any
$g,h \in G$ and $(g_1,g_2),(h_1,h_2) \in G^2$ we have that
$(g_1,g_2) \bot_{g,h} (h_1,h_2)$ if and only if for some
$k,\ell \in \N$, we have $g_1 g^k = g_2$, $h_1 h^\ell = h_2$, and the
walks $(g_1,g_1 g,\dots,g_1 g^k)$ and $(h_1,h_1 h,\dots,h_1 h^\ell)$
are disjoint.
We denote by $\mathcal{F}^\pm$ the set of positive existential
first-order formulas over $\mathcal{C}^\pm(G)$, i.e. formulas
$\exists y_1 \dots \exists y_m \varphi(y_1, \dots, y_m)$ where
$\varphi(y_1, \dots, y_m)$ is a positive Boolean combination of atomic
formulas.  Then $\SAT^\pm(G)$ is the decision problem that asks if a
closed formula in $\mathcal{F}^\pm$ holds in $\mathcal{C}^\pm(G)$.
The fragment $\mathcal{F}^+$ and the
 problem $\SAT^+(G)$ are defined similarly.
Clearly, $\KP^\pm(G)$ (resp. $\KP^+(G)$) reduces to $\SAT^{\pm}(G)$ (resp. $\SAT^+(G)$).
In \cref{appendix-SAT-KP}, we show:
 \begin{theorem}\label{thm:SAT-KP}
For any finitely generated group $G$, the problem $\SAT^\pm(G)$ (resp. $\SAT^+(G)$) is decidable if and only if  $\KP^\pm(G)$ (resp. $\KP^+(G)$) is decidable.
\end{theorem}

\xsubparagraph{Virtually nilpotent groups}  
It was shown by Ganardi, König, Lohrey, and Zetzsche\ that for some number $\ell\in\N$ and all groups $G \neq 1$,
$\KP(G\wr (H_3(\Z)\times\Z^\ell))$ is undecidable~\cite[Theorem
5.2]{GanardiKLZ18}, but essentially nothing is known so far about the
groups $G$ for which the problem $\KP(G\wr H_3(\Z))$ is decidable.
Using \cref{main-result}, this can be settled.
\begin{theorem}\label{application-heisenberg}
  For every non-trivial $G$, the problem $\KP(G\wr H_3(\Z))$ is
  undecidable.
\end{theorem}
This is in contrast to decidability of $\KP(H_3(\Z))$~\cite[Theorem
6.8]{KoenigLohreyZetzsche2015a}.  We show
\cref{application-heisenberg} by proving in \cref{sec:application-heisenberg} that
$\SAT^+(H_3(\Z))$ (and thus $\KP^+(H_3(\Z))$) is undecidable.

The interest in the Heisenberg group stems from its special role
inside the class of virtually nilpotent groups. This class, in turn,
consists exactly of the finite extensions of groups of unitriangular
integer matrices (see, for example, \cite[Theorem
17.2.5]{KaMe1979}). Furthermore, a celebrated result of
Gromov~\cite{gromov1981groups} states that the f.g.\ virtually
nilpotent groups are precisely the f.g.\ groups with polynomial
growth. In some sense, the discrete Heisenberg group is the smallest
f.g.\ virtually nilpotent group that is not virtually abelian.
Therefore, \cref{application-heisenberg} implies the following
characterization of all wreath products $G\wr H$ with decidable
$\KP(G\wr H)$ where $H$ is infinite and virtually nilpotent. See
\cref{appendix-application-virtually-nilpotent} for details.
\begin{corollary}\label{application-virtually-nilpotent}
  Let $G,H$ be f.g.\ non-trivial groups. If $H$ is virtually nilpotent
  and infinite, then $\KP(G\wr H)$ is decidable if and only if $H$ is
  virtually abelian and $\ExpEq(G)$ is decidable.
\end{corollary}
By undecidability of $\ExpEq(H_3(\Z))$, this implies: If
$G\ne 1$ and $H$ are f.g.\ virtually nilpotent and $H$ is infinite,
then $\KP(G\wr H)$ is decidable if and only if $G$ and $H$ are
virtually abelian.

\xsubparagraph{Solvable Baumslag-Solitar groups}
Our second application of \cref{main-result} concerns wreath products
$G\wr \BS(1,q)$.  It is known that knapsack is decidable for
$\BS(1,q)$~\cite[Theorem 4.1]{LohreyZ20}, but again, essentially
nothing is known about $\KP(G\wr\BS(1,q))$ for any
$G$.
\begin{theorem}\label{application-bs}
  For any f.g.\ group $G$ and $q\ge 1$, the problem
  $\KP(G\wr\BS(1,q))$ is decidable if and only if $\ExpEq(G)$ is
  decidable.
\end{theorem}
Extending methods from Lohrey and Zetzsche~\cite{LohreyZ20}, we show that
$\IKPpm(\BS(1,q))$ is decidable for any $q\ge 1$ and thus obtain
\cref{application-bs} in \cref{sec:applications}.

\xsubparagraph{Magnus embedding} Another corollary concerns groups of
the form $F/[N,N]$, where $F$ is a f.g.\ free group and $N$ is a
normal subgroup.  Recall that any f.g.\ group can be written as $F/N$,
where $F$ is an f.g.\ free group and $N$ is a normal subgroup of
$F$. Dividing by $[N,N]$ instead of $N$ yields $F/[N,N]$, which is
subject to the Magnus embedding~\cite[Lemma]{Mag39} of $F/[N,N]$ into
$\Z^r\wr (F/N)$, where $r$ is the rank of $F$. We show
in \cref{appendix-application-magnus}:
\begin{corollary}\label{application-magnus}
  Let $F$ be a finitely generated free group and $N$ be a normal
  subgroup of~$F$.  If $\IKPp(F/N)$ is decidable, then so is
  $\KP(F/[N,N])$. 
\end{corollary}

\xsubparagraph{Knapsack vs. intersection knapsack} Introducing the
problems $\KP^+$ and $\KP^{\pm}$ raises the question of whether they
are substantially different from the similar problems $\KP$ and
$\ExpEq$: Is $\KP^+(G)$ or $\KP^{\pm}(G)$ perhaps inter-reducible with
$\KP(G)$ or $\ExpEq(G)$?  Our applications show that this is not the
case.  Since $\KP(H_3(\Z))$ is decidable~\cite[Theorem
6.8]{KoenigLohreyZetzsche2015a}, but $\KP^+(H_3(\Z))$ is not, neither
$\KP^+(G)$ nor $\KP^{\pm}(G)$ can be inter-reducible with $\KP(G)$ in
general. Moreover, one can show\footnote{Since there is no published
  proof available, we include a proof in \cref{appendix-expeq-bs},
  with kind permission of Moses Ganardi and Markus Lohrey.} that
$\ExpEq(\BS(1,2))$ is undecidable~\cite{GanardiLohrey2020}, whereas
$\KP^{\pm}(\BS(1,q))$ is decidable for any $q\ge 1$. Hence, neither
$\KP^+(G)$ nor $\KP^{\pm}(G)$ can be inter-reducible with $\ExpEq(G)$
in general. However, we leave open whether there is a f.g.\ group $G$
for which $\KP^+(G)$ is decidable, but $\KP^{\pm}(G)$ is undecidable
(see \cref{sec:conclusion}).

\section{From wreath products to intersection knapsack}\label{sec:wr-kp}

In this section, we prove the ``if'' direction of \cref{main-result}
by deciding $\KP(G\wr H)$ using $\ExpEq(G)$ and either $\KP^{\pm}(H)$
or $\KP^+(H)$ (depending on whether $G$ is abelian).

\iffalse
To simplify the construction, instead of $\KP^{\pm}(H)$ and $\KP^+(H)$,
the proof uses slight variations: The \emph{modified intersection
  knapsack problem} $\MKP^\pm(G)$ has the same input instances, but a
disjointness constraint is considered satisfied if the corresponding
walks are disjoint, disregarding the last point. The \emph{positive}
variant $\MKP^+(G)$ of $\MKP^{\pm}$ again differs only by having no
disjointness constraints. See \Cref{app:MKP} for the formal
definitions and a proof that $\MKP^{\pm}(H)$ (resp. $\MKP^+(H)$) is
inter-reducible with $\KP^{\pm}(H)$ (resp. $\KP^+(H)$).
\fi

\xsubparagraph{Normalization}\label{sec:normalization}
We fix a wreath product $G \wr H$ with $G$ and $H$ finitely generated groups. Note that we may assume that $\KP(H)$ is decidable.
In our reduction, we will augment the $\KP(G\wr H)$-instance with positive intersection constraints regarding the cursor in $H$. This results in instances of the \emph{hybrid intersection knapsack problem} $\HKP^\pm(G \wr H)$ over $G \wr H$: It is defined as $\KP^\pm(G \wr H)$ but the loop and disjointness constraints consider the $\sigma$-image of elements.
Let us make this more precise. If $E = \alpha_1 \cdots \alpha_n$ is a knapsack expression over $G \wr H$, then we define for all $i \in [1,n]$ and $\nu \in \N^X$ the set
\[\supp_E^\nu(i) := \{\sigma(\nu(\alpha_1 \cdots \alpha_{i-1}) \gamma(\alpha_i)^k) \mid 0 \leq k \leq \nu(x_i)-1\}\]
if $i \in P_E$ and
\[\supp_E^\nu(i) := \{\sigma(\nu(\alpha_1 \cdots \alpha_{i-1}))\}\]
if $i \in Q_E$. For a walk $w = (w_1,\dots,w_k)$ over $G \wr H$ we write $\sigma(w) := (\sigma(w_1),\dots,\sigma(w_k))$.
Then the \emph{hybrid intersection knapsack problem} $\HKP^\pm(G \wr H)$ over $G \wr H$ is defined as follows:
\begin{description}
\item[Given] a knapsack expression $E$ over $G$, a set $L \subseteq [0,n]^2$ of loop constraints,
and a set $D \subseteq [1,n]^2$ of disjointness constraints.
\item[Question] Is there a valuation $\nu \in \N^X$ with factorized walk $\pi_{\nu,E} = \pi_1 \dots \pi_n$ induced by $\nu$ on $E$ such that the following conditions are fulfilled:
\begin{itemize}
\item $\nu(E) = 1$
\item $\sigma(\pi_{i+1} \dots \pi_j)$ is a loop for all $(i,j) \in L$
\item $\supp_E^\nu(i) \cap \supp_E^\nu(j) = \emptyset$ for all $(i,j) \in D$.
\end{itemize}
\end{description}
Its \emph{positive} version $\HKP^+(G\wr H)$ is again defined by having no disjointness constraints. The set $\sol_{G\wr H}$ is defined accordingly. Note that to simplify the constructions in the proofs, the disjointness constraints in an $\HKP^\pm(G \wr H)$-instance disregard the last point of walks.

In the following, when we write a knapsack expression as $E = \alpha_1 \cdots \alpha_n \alpha_{n+1}$, we assume w.l.o.g. that $\alpha_{n+1}$ is a constant.
Two elements $g,h \in H$ are called {\em commensurable} if $g^x = h^y$ for some $x,y \in \Z \setminus \{0\}$.
It is known that if $g_1,g_2$ have infinite order and are not commensurable, then there is at most one solution $(x_1,x_2) \in \Z^2$
for the equations $g_1^{x_1} g_2^{x_2} = g$~\cite[Lemma~9]{FigeliusGLZ20}.

\iffalse
We say that a knapsack expression $E = \alpha_1 \cdots \alpha_n$ is \textit{torsion-free} if for all $i \in P_E$ it holds that $\sigma(\gamma(\alpha_i)) = 1$ or $\sigma(\gamma(\alpha_i))$ has infinite order.

A knapsack expression $E = \alpha_1 \cdots \alpha_n \alpha_{n+1}$ is in $GH$\textit{-form} if for all $i \in P_E$ it holds that $\sigma(\gamma(\alpha_i)) = 1$ or $\gamma(\alpha_i) \in GH$ and for all $i \in Q_E \setminus \{n+1\}$ it holds that $\alpha_i \in H$.
\fi
 
Let $E = \alpha_1 \cdots \alpha_n\alpha_{n+1}$ be a knapsack
expression and write $g_i=\gamma(\alpha_i)$ for $i\in[1,n+1]$. The
expression (resp. the corresponding $\HKP^{\pm}(G\wr H)$-instance) is
\textit{c-simplified} if for any
$i, j \in P_E$ with $g_i\notin H$ and $g_j\notin H$, we
have that commensurability of $\sigma(g_i)$ and $\sigma(g_j)$ implies
$\sigma(g_i) = \sigma(g_j)$.  We call the expression (resp. the corresponding
$\HKP^{\pm}(G\wr H)$-instance) \emph{normalized} if it is c-simplified and
each atom $\alpha_i$ with $i\in[1,n]$ is of one of the following
types: We either have (a)~$i \in Q_E$ and $g_i \in H$ or
(b)~$i \in P_E$ and $\sigma(g_i) = 1$ or (c)~$i \in P_E$,
$g_i \in GH$ and $\sigma(g_i)$ has infinite
order. Using generalizations of ideas from \cite{GKLZ17} and
\cite{FGLZ20}, we show:
\begin{theorem}\label{thm:normal}
Given an instance of $\KP(G\wr H)$, one can effectively construct an equivalent finite set of normalized $\HKP^+(G \wr H)$-instances.  
\end{theorem}
Here, a problem instance $I$ is \emph{equivalent} to a set
$\mathcal{I}$ of problem instances if $I$ has a solution if and only
if at least one of the instances in $\mathcal{I}$ has a solution. 

\iffalse
Thus, we may assume that for an $\HKP^+(G \wr H)$-instance $(E = \alpha_1 \cdots \alpha_n \alpha_{n+1}, L)$ we have that for all $i \in [1,n]$ one of the following properties holds:
\begin{enumerate}[(a)]
\item 
\item 
\item 
\end{enumerate}
A knapsack expression (resp. $\HKP^\pm(G \wr H)$-instance) of this form is said to be \textit{normalized}.
\fi

%
%

\xsubparagraph{Non-abelian case} Note that in a normalized knapsack
expression, atoms of type (b) and (c) and the last atom $\alpha_{n+1}$
may place non-trivial elements of $G$. Our next step is to transform the
input instance further so that only the atoms of type (c) can place
non-trivial elements of $G$, which leads to the notion of stacking-freeness.

Let $E = \alpha_1 \cdots \alpha_n \alpha_{n+1}$ be a knapsack expression over $G \wr H$ and let $g_i := \gamma(\alpha_i)$ for all $i \in [1,n+1]$.
We call an index $i \in [1,n+1]$ {\em stacking} if either $i \in P_E$ and $\sigma(g_i) = 1$, or $i = n+1$ and $g_{n+1} \notin H$.
We say that $E$ is \textit{stacking-free} if it has no stacking indices.
Thus, a normalized expression $E$ is stacking-free if each atom is either of type (c) or a constant in $H$.

\iffalse
value in $G$ that walks through the Cayley graph of $H$ and this walk forms a ray. If $(i,h)$ is an address with $i \in P_E$ and $\sigma(\gamma(\alpha_i)) \neq 1$ and $\nu \in \N^X$ is a valuation, then
$(\sigma(\nu(\alpha_1 \cdots \alpha_{i-1})) h (h^{-1} \sigma(\gamma(\alpha_i)) h)^j)_{0 \leq j \leq \nu(x_i)-1}$
is the ray associated to $(i,h)$.
\fi
\begin{lemma}\label{lem:stacking-free2}
Given a normalized $\HKP^\pm(G \wr H)$-instance, one can effectively construct an equivalent finite set of stacking-free, normalized $\HKP^\pm(G \wr H)$-instances.
\end{lemma}
Let us sketch the proof of \cref{lem:stacking-free2}.
We use the notion of an address from \cite{GKLZ17}. An \textit{address} of $E$ is a pair $(i,h)$ with $i \in [1,n+1]$ and $h \in H$ such that $h \in \supp(\gamma(\alpha_i))$. The set of addresses $A_E$ of $E$ is finite and can be computed. Intuitively, an address represents a position in a knapsack expression where a point in $H$ can be visited.

Intuitively, instead of placing elements of $G$ by atoms of type~(b)
and by $\alpha_{n+1}$, we introduce loop and disjointness constraints
guaranteeing that in points visited by these atoms, a solution would
have placed elements that multiply to $1\in G$. To this end, we pick
an address $(i,h)\in A$ of a stacking index $i$ and then guess a set
$C\subseteq A$ of addresses such that the point $h'\in H$ visited at
$(i,h)$ is visited by exactly the addresses in $C$. The latter
condition is formulated using loop and disjointness constraints in an
$\HKP^\pm(G \wr H)$-instance $I_C$. In $I_C$, we do not place elements
at $C$ anymore; instead, we construct a set $S_C$ of exponent
equations over $G$ that express that indeed the point $h'$ carries
$1\in G$ in the end. Note that this eliminates one address with stacking index. We repeat this until we are left with a set of stacking-free instances of
$\HKP^\pm(G\wr H)$, each together with an accumulated set of exponent
equations over $G$. We then take the subset $\mathcal{I}$ of
$\HKP^{\pm}(G \wr H)$-instances whose associated $\ExpEq(G)$-instance has a
solution.  This will be our set for \cref{lem:stacking-free2}.

The last step of the non-abelian case is to construct $\KP^{\pm}(H)$-instances.
\begin{lemma}\label{lem:interval-KP2}
Given a stacking-free, normalized $\HKP^\pm(G \wr H)$-instance, one can effectively construct an equivalent finite set of $\KP^\pm(H)$-instances.
\end{lemma}

We are given an instance $(E,L,D)$ with $E=\alpha_1\cdots\alpha_n$ and
write $g_i=\gamma(\alpha_i)$ for $i\in[1,n]$.  As $(E,L,D)$ is
normalized and stacking-free, only atoms of type (c) with
$g_i\notin H$ can place non-trivial elements of $G$. Moreover,
if $\alpha_i$ and $\alpha_j$ are such atoms, then the elements
$\sigma(g_i)$ and $\sigma(g_j)$ are either non-commensurable or equal.
In the first case, the two rays produced by $\alpha_i$ and $\alpha_j$
can intersect in at most one point; in the second case, they intersect
along subrays corresponding to intervals $I_i\subseteq [0,\nu(x_i)]$
and $I_j\subseteq[0,\nu(x_j)]$.

Thus, the idea is to split up each ray wherever the intersection with another ray starts or ends: We guess for each ray as above the number $m \leq 2 \cdot |A_E| - 1$ of subrays it will be split into and replace $g_i^{x_i}$ with $g_i^{y_1} \cdots g_i^{y_m}$. After the splitting, subrays are either equal or disjoint. We guess an equivalence relation on the subrays; using loop constraints, we ensure that subrays in the same class are equal; using disjointness constraints, we ensure disjointness of subrays in distinct classes. Finally, we have to check that for each equivalence class $C$,
the element of $G$ produced by the rays in $C$ does indeed multiply to $1\in G$. This can be
checked because $\ExpEq(G)$ (and thus the word problem for $G$) is decidable.

\xsubparagraph{Abelian case} We now come to the case of abelian $G$:
We show that $\KP(G\wr H)$ is decidable, but only using instances of
$\KP^+(H)$ instead of $\KP^\pm(H)$. Here, the key insight is that we
can use the same reduction, except that we just do not impose the
disjointness constraints.  In the above reduction, we use disjointness
constraints to control exactly which positions in our walk visit the
same point in $H$. Then we can check that in the end, each point in
$H$ carries $1\in G$. However, if $G$ is abelian, it suffices to make
sure that the set of positions in our walk decomposes into subsets,
each of which produces $1\in G$: If several of these subsets
do visit the same point in $H$, the end result will still be $1\in G$.

We illustrate this in a slightly simpler setting. Suppose we have a
product $g=\loi{h_1}{a}{}{1}\cdots \loi{h_n}{a}{}{n}$ with
$h_1,\ldots,h_n\in H$ and $a_1,\ldots,a_n\in G$. Then $g$ is obtained
by placing $a_1$ at $h_1\in H$, then $a_2$ at $h_2\in H$, etc.  For a
subset $S=\{s_1,\ldots,s_k\}\subseteq[1,n]$ with
$s_1<\cdots<s_k$, we define
$g_S=\loi{h_{s_1}}{a}{}{s_1}\cdots \loi{h_{s_k}}{a}{}{{s_k}}$. Hence,
we only multiply those factors from $S$. An equivalence relation
$\equiv$ on $[1,n]$ is called \emph{cancelling} if $g_C=1$ for every
class $C$ of $\equiv$. Moreover, $\equiv$ is called \emph{equilocal}
if $i\equiv j$ if and only if $h_i=h_j$.  It is called \emph{weakly
  equilocal} if $i\equiv j$ implies $h_i=h_j$.  Now observe that for
any $G$, we have $g=1$ if and only if there is an equilocal cancelling
equivalence on $[1,n]$. However, if $G$ is abelian, then $g=1$ if and
only if there is a \emph{weakly} equilocal equivalence on $[1,n]$.
Since weak equilocality can be expressed using only equalities (and no
disequalities), it suffices to impose loop conditions in our instances.

\xsubparagraph{Comparison to previous approach in \cite{FGLZ20}} The
reduction from $\KP(G\wr H)$ to $\ExpEq(G)$ and $\KP^{\pm}(H)$
($\KP^{+}(H)$ respectively) uses similar ideas as the proof of
\cite[Theorem 4]{FGLZ20}, where it is shown $\ExpEq(K)$ is in $\NP$ if
$K$ is an iterated wreath product of $\Z^r$ for some $r\in\N$.

Let us compare our reduction with the proof of \cite[Theorem
4]{FGLZ20}.  In \cite{FGLZ20}, one solves $\ExpEq(K)$ by writing
$K=G\wr H$ where $G$ is abelian and $H$ is orderable and
knapsack-semilinear.  In both proofs, solvability of an instance (of
$\ExpEq(G\wr H)$ in \cite{FGLZ20} and $\KP(G\wr H)$ here) is
translated into a set of conditions by using similar decomposition
arguments. Then, the two proofs differ in how satisfiability of these
conditions is checked.

In
\cite{FGLZ20}, this set of conditions is expressed in Presburger
arithmetic, which is possible due to knapsack-semilinearity of $H$.
In our reduction, we have to translate the conditions in $\ExpEq(G)$ and
$\KP^+(H)$ ($\KP^{\pm}(H)$) instances. Here, we use loop constraints
where in Presburger arithmetic, once can compare variables directly.
Moreover, our reduction uses disjointness constraints to express solvability
in the case that $G$ is non-abelian. This case does not occur in \cite[Theorem 4]{FGLZ20}.
Finally, we have to check whether the elements from $G$ written at the
same point of $H$ multiply to 1. The reduction of \cite{FGLZ20} can
express this directly in Presburger arithmetic since $G$ is
abelian. Here, we use instances of $\ExpEq(G)$.

\section{From intersection knapsack to wreath products}\label{sec:to-wreath}

In this section, we prove the ``only if'' direction of
\cref{main-result}. Since it is known that for infinite $H$,
decidability of $\KP(G\wr H)$ implies decidability of
$\ExpEq(G)$~\cite[Proposition.~3.1, Proposition~5.1]{GanardiKLZ18}, it
remains to reduce (i)~$\KP^+(H)$ to $\KP(G \wr H)$ for any group
$G \neq 1$, and (ii)~ $\KP^\pm(H)$ to $\KP(G \wr H)$ for any
non-abelian group $G$.  In the following, let $G$ be a non-trivial
group and $H$ be any group and suppose $\KP(G\wr H)$ is decidable.

First let us illustrate how to reduce $\KP^+(H)$ to $\KP(G \wr H)$.
Suppose we want to verify whether a product $h_1 \dots h_m = 1$ over $H$
satisfies a set of loop constraints
$L \subseteq [0,m]^2$,
i.e. $h_{i+1} \dots h_j = 1$ for all $(i,j) \in L$.
To do so we insert into the product for each $(i,j) \in L$ a function $f \in G^{(H)}$ after the element $h_i$
and its inverse $f^{-1}$ after the element $h_j$.
We call these functions {\em loop words} since their supports are contained
in a cyclic subgroup $\langle t \rangle$ of $H$.
We can choose the loop words such that
this modified product evaluates to 1 if and only if the loop constraints are satisfied.
For the reduction from $\KP^\pm(H)$ we need to make the construction more robust
since we simultaneously need to simulate disjointness constraints.

If $H$ is a torsion group then $\KP^+(H)$ and $\KP^\pm(H)$ are
decidable if the word problem of $H$ is decidable: For each exponent,
we only have to check finitely many candidates.  Since $\KP(G\wr H)$
is decidable, we know that $\KP(H)$ is decidable and hence also the
word problem.  Thus, we assume $H$ not to be a torsion group and may
fix an element $t \in H$ of infinite order.

\xsubparagraph{Periodic complexity}
Let $K$ be a group. The following definitions will be employed with $K=\Z$ or $K=H$.
For any subset $D \subseteq K$, let $G^{(D)}$ be the group of all functions $u \colon K \to G$
whose support $\supp(u) = \{ h \in K \mid u(h) \neq 1 \}$ is finite and contained in $D$.
A function $f \in G^{(K)}$ is {\em basic periodic}
if there exists a progression $D$ in $K$ and $c \in G$ such that $f(h) = c$ for all $h \in D$
and $f(h) = 1$ otherwise. 
The {\em value} of such a function $f$ is the element $c$;
a {\em period} of $f$ is a period of its support.
We will identify a word $u = c_1 \dots c_n \in G^*$
with the function $u \in G^{(\Z)}$ where $u(i) = c_i$ for $i \in [1,n]$ and $u(i) = 1$ otherwise.
Recall that for $u \in G^{(\Z)}$ and $s\in\Z$, we have $\lo{s}{u}(n)=u(n-s)$. We extend this to $s\in \Z_\infty := \Z\cup\{\infty\}$ by setting $\lo{\infty}{u}(n)=1$ for all $n\in\Z$.
The {\em periodic complexity} of $u \in G^{(\Z)}$ is the minimal number $\pc(u) = k$
of basic periodic functions $u_1, \dots, u_k$
such that $u = \prod_{i=1}^k u_i$.
Given a progression $D = \{ p+qn \mid n \in [0,\ell] \}$ in $\Z$ and a function $u \in G^{(\Z)}$ 
we define $\pi_D(u)(n) = u(p+qn)$ for all $n \in \Z$ and say that $\pi_D(u)$ is a {\em periodic subsequence} of $u$.
Note that periodic subsequences of basic periodic functions are again basic periodic.
Furthermore, since $\pi_D \colon G^{(\Z)} \to G^{(\Z)}$ is a homomorphism,
taking periodic subsequences does not increase the periodic complexity.

\begin{lemma}
	\label{lem:interval-words}
	Given $n,k \in \N$ and $a \in G \setminus \{1\}$,
	one can compute $u_1, \dots, u_n \in \langle a \rangle^{(\N)}$
	such that $\prod_{i=1}^n \loi{p_i}{u}{}{i} \loi{q_i}{u}{-1}{i}$ has periodic complexity $\ge k$
	for all $(p_1, \dots, p_n) \neq (q_1, \dots, q_n) \in \Z_\infty^n$.
\end{lemma}

Here is a proof sketch for \cref{lem:interval-words}.
The case $n = 1$ can be shown by taking any function $v = a_1 \dots a_m \in \langle a \rangle^{(\N)}$
with large periodic complexity
and defining $u_1 = a_1 (1)^{m-1} a_2 (1)^{m-1} \dots a_m (1)^{m-1} a_1 \dots a_m$
where $(1)^{m-1}$ is the sequence consisting of $m-1$ many $1$'s.
If $p,q \in \Z_\infty$ are distinct then $\loi{p}{u}{}{1} \loi{q}{u}{-1}{1}$ always
contains $v$ or $v^{-1}$ as a periodic subsequence and thus has large periodic complexity.
For $n > 1$ we define $u_i$ ($i > 1$) to be stretched versions of $u_1$
such that the supports of any two functions $\loi{p}{u}{}{i}$, $\loi{q}{u}{}{j}$ where $i \neq j$
intersect in at most one point.
This allows to argue that $\prod_{i=1}^n \loi{p_i}{u}{}{i} \loi{q_i}{u}{-1}{i}$ still has
large periodic complexity as soon as $p_i \neq q_i$ for some $i$.

\xsubparagraph{Expressing loop constraints} We now show how to use
\cref{lem:interval-words} to encode loop constraints over a
product $h_1 \dots h_m$ over $H$ in an instance of $\KP(G \wr H)$.

Recall that a loop constraint $(i,j)$ stipulates that $\sigma(g_{i+1} \dots g_j) = 1$.
If we only want to reduce $\KP^+(H)$, it is not hard to see that it
would suffice to guarantee
$\prod_{i=1}^n \loi{p_i}{u}{}{i}\loi{q_i}{u}{-1}{i}\ne 1$ in
\cref{lem:interval-words}. In that case, we could essentially use the
functions $u_i$ as loop words. However, in order to express
disjointness constraints in $\KP^\pm(H)$, we will construct expressions
over $G\wr H$ that place additional ``disjointness patterns'' in the
Cayley graph of $H$.  %
We shall make sure that the disjointness patterns
are tame: Roughly speaking, this means they are basic periodic and
either (i)~place elements from a fixed subgroup $\langle a\rangle$ or (ii)~can
intersect a loop word at most once. Here, the high periodic complexity
of $\prod_{i=1}^n \loi{p_i}{u}{}{i}\loi{q_i}{u}{-1}{i}$ will allow
us to conclude that tame patterns cannot make up for a violated loop
constraint.

Let us make this precise.
Recall that two elements $g,h \in H$ are called {\em commensurable}
if $g^x = h^y$ for some $x,y \in \Z \setminus \{0\}$.
Let $a \in G \setminus \{1\}$.
Let $\mathsf{P}_{a,t}(G \wr H)$
be the set of elements $g \in G \wr H$
such that $\tau(g)$ is basic periodic and either,
(i)  its value belongs to $\langle a \rangle$, or
(ii) its period is not commensurable to $t$.
In particular, a power $(ch)^k$ (where $c \in G$, $h \in H$, $k \in \N$)
belongs to $\mathsf{P}_{a,t}(G \wr H)$ if $c \in \langle a \rangle$
or $h$ is not commensurable to $t$. Note that since loop words are
always placed along the direction $t$, this guarantees tameness: In
case (ii), the period of $\tau(g)$ being non-commensurable to $t$
implies that the support of any $h'g$, $h'\in H$, can intersect the
support of a loop word in $\langle a\rangle^{(\langle t\rangle)}$ at
most once.
Using \cref{lem:interval-words}, we show the following.
\begin{lemma}
  \label{lem:loop}
  Given $a \in G \setminus \{1\}$, $m \in \N$ and $L \subseteq [0,m]^2$
  we can compute $f_0, \dots, f_m \in \langle a \rangle^{(t^*)}$
  such that:
  \begin{enumerate}
  \item Let $h_1, \dots, h_m \in H$.
    Then $h_1 \dots h_m = 1$ and $h_{i+1} \dots h_j = 1$ for all $(i,j) \in L$
    if and only if $f_0 h_1 f_1 \dots h_m f_m = 1$. 
  \item Let $g_1, \dots, g_m \in \mathsf{P}_{a,t}(G \wr H)$ 
    such that $\sigma(g_{i+1} \dots g_j) \neq 1$ for some $(i,j) \in L$.
    Then $f_0 g_1 f_1 \dots g_m f_m \neq 1$. 
  \end{enumerate}
\end{lemma}
Observe that the first constraint says that if we only use the loop
words $f_i$, then they allow us to express loop constraints. The
second constraint tells us that a violated loop constraint cannot be
compensated even with perturbations $g_1,\ldots,g_m$, provided that
they are tame.

\xsubparagraph{The abelian case}
\cref{lem:loop} provides a simple reduction from $\KP^+(H)$ to $\KP(G \wr H)$.
Given an instance $(E = e_1 \dots e_n, L)$ of $\KP^+(H)$
we compute $f_0, \dots, f_m \in \langle a \rangle^{(t^*)}$
using \cref{lem:loop}.
Then  $\nu\colon X\to\N$ satisfies $\nu(E) = 1$ and $\nu(e_{i+1} \dots e_j)$ for all $(i,j) \in L$
if and only if $\nu(f_0 e_1 f_1 \dots e_n f_n) = 1$.
Hence $(E,L)$ has a solution if and only if $\nu(f_0 e_1 f_1 \dots e_n f_n) = 1$ does.

\xsubparagraph{The non-abelian case}
Now let $G$ be a non-abelian group.
In the following we will reduce $\KP^\pm(H)$ to $\KP(G \wr H)$.
The first step is to construct from an $\KP^\pm(H)$-instance $I$
an equivalent $\HKP^+(G \wr H)$-instance $\hat I$
using a nontrivial commutator $[a,b] \neq 1$ in $G$.
In a second step we apply the ``loop words''-construction from \cref{lem:loop} (point 2) to $\hat I$,
going to a (pure) knapsack instance.
It guarantees that, if a loop constraint is violated, then the knapsack instance does not evaluate to 1.
Furthermore, if a disjointness constraint is violated
then there exists a large number of pairwise distant points
in the Cayley graph of $H$ which are labeled by a nontrivial element.
These points cannot be canceled by the functions $f_i$ from \cref{lem:loop}.
Finally, if all loop and disjointness constraints are satisfied
then the induced walk in the Cayley graph provides enough ``empty space'' 
such that the loop words can be shifted to be disjoint
from the original walk induced by $\hat I$ (encoding the disjointness constraints).

\xsubparagraph{Normalization} Let $I = (E = e_1 \dots e_n,L,D)$ be a
$\KP^\pm(H)$-instance where $e_i$ is either a
constant $e_i = h_i$ or a power $e_i = h_i^{x_i}$.  We will start by
establishing the following useful properties.  We call $I$ {\em
  torsion-free} if $h_i$ has infinite order for all $i \in P_E$.  Call
$I$ {\em orthogonalized} for all $(i,j) \in D \cap P_E^2$ such that we
have $\langle h_i \rangle \cap \langle h_j \rangle = \{1\}$.  If $I$
is torsion-free and orthogonalized then it is called {\em normalized}.
The orthogonality will be crucial for the tameness of the disjointness patterns
since at most one of the elements $h_i,h_j$ for $(i,j) \in D \cap P_E^2$ is commensurable to $t$.
Furthermore, it guarantees that there is at most one intersection point
for any pair $(i,j) \in D$.

\begin{lemma}
  \label{lem:kppm-normal}
  One can compute a finite set $\mathcal{I}$ of normalized
  instances of $\KP^\pm(H)$ such that
  $I$ has a solution if and only if there exists $I' \in \mathcal{I}$ which has a solution.
\end{lemma}
Here, torsion-freeness is easily achieved: If $h_i$ has finite order,
then $h_i^{x_i}$ can only assume finitely many values, so we replace
$h_i^{x_i}$ by one of finitely many constants. Orthogonality requires
an observation: If $\langle h_i\rangle\cap\langle h_j\rangle\ne\{1\}$,
then any two intersecting progressions $\pi_i, \pi_j$ with periods $h_i$ and $h_j$, respectively,
must intersect periodically, meaning
there exists an intersection point that is close to an endpoint of $\pi_i$ or $\pi_j$.
This means, in lieu of $(i,j)\in D$, we can require
disjointness of one power with a constant.

\xsubparagraph{Expressing disjointness constraints} Hence we can assume
that $I$ is normalized.  To express disjointness constraints, we must
assume that $G$ is non-abelian. Let $a,b \in G$ with
$aba^{-1}b^{-1}=[a,b] \neq 1$.  Our starting point is the following
idea. To express that two progressions $\pi_i$ and $\pi_j$, induced by a valuation of $E$,
are disjoint, we construct an
expression over $G\wr H$ that first places $a$ at each point in $\pi_i$,
then $b$ at each point in $\pi_j$, then again $a^{-1}$ at each point in
$\pi_i$, and finally $b^{-1}$ at each point in $\pi_j$,
see~\eqref{eq:commutators}.
Here we need loop constraints that express that the start and endpoints
of the two traversals of $\pi_i$ (and $\pi_j$) coincide.
Then, if $\pi_i$ and
$\pi_j$ are disjoint, the effect will be neutral; otherwise
any intersection point will carry $aba^{-1}b^{-1}\ne 1$.

However, this leads to two problems.  First, there might be more than
one disjointness constraint: If $k$ disjointness constraints are
violated by the same point $h''\in H$, then $h''$ would carry
$[a,b]^k$, which can be the identity (for example, $G$ may be
finite). Second, when we also place loop words (which multiply elements from $\langle a\rangle$), those could also interfere with the commutator (for
example, instead of $aba^{-1}b^{-1}$, we might get
$aba^{-1}(a)b^{-1}(a^{-1})=1$).

Instead, we do the following. Let $t \in H$ be the element of
infinite order used for the loop words. Moreover, let
$D = \{ (i_1,j_1), \dots, (i_d,j_d) \}$.  For each $(i_k,j_k)\in D$,
instead of performing the above ``commutator construction'' once, we
perform it $n+d$ times, each time shifted by $t^{N_k}\in H$ for some
large $N_k$.  The numbers $N_0<N_1<\cdots$ are chosen so large that
for at least one commutator, there will be no interference from other
commutators or from loop words.

Let us make this precise. Since $I$ is orthogonalized, we may assume
that for each $(i,j)\in D\cap P_E^2$, the elements $h_j$ and $t$ are
not commensurable; otherwise we swap $i$ and $j$.  The resulting
$\HKP^+(G \wr H)$-instance $\hat I$ will have length
$m = n + 4d(n+d)(n+2)$.  In preparation, we can compute a number $N$
such that the functions $f_0, \dots, f_m$ from \cref{lem:loop} for any
$L \subseteq [0,m]^2$ satisfy
$\supp(f_i) \subseteq \{t^j \mid j \in [0,N-1] \}$.
For each $i \in [1,n]$, $c \in G$, $s \in \N$, we define
the knapsack expression $E_{i,c,s}$ over $G \wr H$ as
\begin{equation}
	\label{eq:eics}
	E_{i,c,s} = \begin{cases}
		e_1 \dots e_{i-1} \, (t^s) \, (c \, t^{-s} h_i t^s)^{x_i} (c t^{-s}) \, e_{i+1} \dots e_n, & \text{if } e_i = h_i^{x_i},\\
		e_1 \dots e_{i-1} \, (t^s) \; (c \, t^{-s} h_i t^s) \;\; (c t^{-s}) \, e_{i+1} \dots e_n, & \text{if } e_i = h_i.
	\end{cases}
\end{equation}
The parentheses indicate the atoms.
We define
\begin{equation}
	\label{eq:commutators}
	\hat E = E \cdot \prod_{k = 1}^d \prod_{s \in S_k}
	\Big(E_{i_k,a,s} \cdot E_{j_k,b,s} \cdot E_{i_k,a^{-1},s} \cdot E_{j_k,b^{-1},s} \Big)
\end{equation}
where $S_k = \{ j (n+d)^{2k} N \mid j \in [1,n+d] \}$ for all $k \in [1,d]$,
and all occurrences of expressions of the form $E_{i,c,s}$ use fresh variables.
Note that $E_{i_k,a,s} \cdot E_{j_k,b,s} \cdot E_{i_k,a^{-1},s} \cdot E_{j_k,b^{-1},s}$
performs the commutator construction for $(i_k,j_k)$, shifted by $t^s$.
Let $\hat E = \hat e_1 \dots \hat e_m$ be the resulting expression.
Notice that its length is indeed $m = n + 4d(n+d)(n+2)$ as claimed above.

Finally, in our $\HKP^+(G\wr H)$ instance, we also add a set
$J\subseteq[0,m]^2$ of loop constraints stating that for each
$k\in[1,d]$ and $s\in S_k$, the $i_k$-th atom in $E_{i_k,a,s}$ arrives
at the same place in $H$ as the $i_k$-th atom in $E$ (and analogously
for $E_{j_k,b,s}$, $E_{i_k,a^{-1},s}$, $E_{j_k,b^{-1},s}$). See
\cref{to-wreath:loop-constraints} for details.

Let $f_0, \dots, f_m \in \langle a \rangle^{(t^*)}$ be the loop words
from \cref{lem:loop} for the set $J \subseteq [0,m]^2$. It is now
straightforward to verify that the elements $\hat{e}_i$ are all tame
as explained above. In other words, for every valuation $\nu$ and
$i\in[1,m]$, we have $\nu(\hat{e}_i)\in \mathsf{P}_{a,t}$ (see
\cref{lem:t-com}).

\xsubparagraph{Shifting loop words} By construction, we now know that
if the instance $f_0\hat{e}_1f_1\cdots \hat{e}_mf_m$ of $\KP(G\wr H)$
has a solution, then so does our normalized instance $I$ of
$\KP^{\pm}(H)$.  However, there is one last obstacle: Even if all loop and
disjointness constraints can be met for $I$, we cannot guarantee that
$f_0\hat{e}_1f_1\cdots \hat{e}_mf_m$ has a solution: It is possible
that some loop words interfere with some commutator constructions so
as to yield an element $\ne 1$.

The idea is to \emph{shift} all the loop words $f_0,\ldots,f_m$ in
direction $t$ by replacing $f_i$ by
$t^{r}f_it^{-r}=\loi{t^r\!\!}{f}{}{i}$ for some $r\in\N$. We shall
argue that for some $r$ in some bounded interval, this must result in
an interference free expression; even though the elements $\hat{e}_i$
may modify an unbounded number of points in $H$. To this end, we use
again that the $\hat{e}_i$ are tame: Each of them either (i)~places
elements from $\langle a\rangle$, or (ii)~has a period
non-commensurable to $t$. In the case (i), there can be no
interference because the $f_i$ also place elements in
$\langle a\rangle$, which is an abelian subgroup. In the case (ii),
$\hat{e}_i$ can intersect the support of each $f_j$ at most once.
Hence, there are at most $m$ points each $f_j$ has to avoid after
shifting.
The following simple lemma states that one can always shift finite sets $F_i$ in parallel
to avoid finite sets $A_i$, by a bounded shift.
Notice that the bound does not depend on the size of the elements in the sets $F_i$ and $A_i$.

\begin{lemma}
	\label{lem:shifting}
	Let $F_1,\ldots,F_m\subseteq\Z$ with $|F_i|\le N$ and
	$A_1,\ldots,A_m\subseteq\Z$ with $|A_i|\le \ell$.
	There exists a shift $r\in[0,Nm\ell]$ such that
	$(r+F_i)\cap A_i=\emptyset$ for each $i\in[1,m]$.
\end{lemma}

\begin{proof}
For every $a \in \Z$ there exist at most $|F_i| \le N$
many shifts $r \in \N$ where $a \in r+F_i$.
Therefore there must be a shift $r\in[0,Nm\ell]$ such that
$(r+F_i)\cap A_i=\emptyset$ for each $i\in[1,m]$.
\end{proof}

We can thus prove the following lemma, which clearly
completes the reduction from $\KP^\pm(H)$ to $\KP(G \wr H)$.

\begin{lemma}
	\label{lem:main-lem}
	$I = (E,L,D)$ has a solution if and only if
	$\loi{t^r\!\!}{f}{}{0} \hat e_1 \loi{t^r\!\!}{f}{}{1} \dots \hat e_m \loi{t^r\!\!}{f}{}{m}$
	has a solution for some $r \in [0,Nm^2]$.
\end{lemma}

\section{Applications}\label{sec:applications}

\label{sec:application-heisenberg}
\xsubparagraph{The discrete Heisenberg group}
Here, we prove that $\SAT^+(H_3(\Z))$ is undecidable. Together with \cref{main-result} and \cref{thm:SAT-KP}, this directly implies \cref{application-heisenberg}. Define the matrices
$A = \begin{psmallmatrix} 1 & 1 & 0 \\ 0 & 1 & 0 \\ 0 & 0 &
  1 \end{psmallmatrix}$,
$B = \begin{psmallmatrix} 1 & 0 & 0 \\ 0 & 1 & 1 \\ 0 & 0 &
  1 \end{psmallmatrix}$, and
$C = \begin{psmallmatrix} 1 & 0 & 1 \\ 0 & 1 & 0 \\ 0 & 0 &
  1 \end{psmallmatrix}$.  The group $H_3(\Z)$ is generated by $A$ and
$B$ and we have $AC=CA$ and $BC=CB$. It is well-known that
(I)~$A^iC^j=A^{i'}C^{j'}$ iff
$i=i'$ and $j=j'$; and (II)~$B^iC^j=B^{i'}C^{j'}$ iff $i=i'$ and
$j=j'$; and (III)~$A^iB^jA^{-i'}B^{-j'}=C^k$ if and only if $i=i'$, $j=j'$, and
$k=ij$. For proofs, see \cref{sec:appendix-h3}.

We show undecidability of
$\SAT^+(H_3(\Z))$ by reducing from solvability of Diophantine
equations over natural numbers. Hence, we are given a finite system
$\bigwedge_{j=1}^m E_j$ of equations of the form $x=a$, $z=x+y$, and $z=xy$.
It is well-known that solvability of such equation systems is undecidable~\cite{Mat93}.
Given such an equation system over a set of variables $X$
we define a $\mathcal{C}^+(H_3(\Z))$-formula containing the variables $\{g_x \mid x \in X \} \cup \{g_0\}$
with the interpretation that $g_x=g_0 C^x$.
First we state that $g_0 \tostar[C] g_x$ for all $x \in X$.
Expressing $x=a$ is done simply with $g_0 \xrightarrow{C^a} g_x$. For $z=x+y$, we use
\[ C^xA^* \cap A^{x'}C^*\cap (AC)^*\ne \emptyset ~~\wedge~~
  A^{x'}C^*\cap C^zA^*\cap C^y (AC)^*\ne \emptyset. \] This can
be expressed in $\mathcal{C}^+(H_3(\Z))$ with a fresh variable $f_{x'}$ for
$g_0 A^{x'}$:
For example, the first conjunct holds iff there exists $h \in H_3(\Z)$
such that $g_0 \tostar[A] f_{x'}$, $g_x\tostar[A] h$,
$f_{x'} \tostar[C] h$, $g_0 \tostar[AC] h$.
By (I) and $AC=CA$, the first conjunct holds iff $x=x'$. Similarly, the second conjunct holds iff
$z=x'+y$, hence $z=x+y$. For $z=xy$, we use:
\[ C^xA^* \cap A^{x'}C^*\cap (AC)^*\ne \emptyset ~~\wedge~~B^{y'}C^* \cap C^yB^* \cap (BC)^*\ne \emptyset \\
  ~~\wedge~~ A^{x'}B^*(A^{-1})^*\cap B^{y'}C^*\cap
  C^zB^*\ne\emptyset.\] Like above, the first and second conjunct
express $x'=x$ and $y'=y$. The third says that
$A^{x'}B^r(A^{-1})^{s}=B^{y'}C^z$ for some $r,s\ge 0$, so by (III),
it states $z=x'y'$, hence $z=xy$.

\xsubparagraph{Solvable Baumslag-Solitar groups} We show that
$\SAT^{\pm}(\BS(1,q))$ is decidable for every $q\ge 1$. By
\cref{main-result} and \cref{thm:SAT-KP}, this proves
\cref{application-bs}. Our proof is based on the following observation, which is shown in \cref{sec:appendix-bs}.
\begin{proposition}\label{fo-cplus-bs}
  The first-order theory of $\mathcal{C}^+(\BS(1,q))$ is decidable.
\end{proposition}
For \cref{fo-cplus-bs}, we show that given any finite subset
$F\subseteq\BS(1,q)$, the structure
$(\BS(1,q),(\xrightarrow{g})_{g \in F},(\tostar[g])_{g \in F})$ is
effectively an automatic structure, which implies that its first-order
theory is decidable~\cite[Corollary 4.2]{khoussainov1994automatic}.
This uses a straightforward extension of the methods
in~\cite{LohreyZ20}.  In~\cite[proof of Theorem 4.1]{LohreyZ20}, it is shown that
$\KP(\BS(1,q))$ can be reduced to the existential fragment of the
structure $(\Z, +, V_q)$, where $V_q(n)$ is the largest power of $q$
that divides $n$. The structure $(\Z,+,V_q)$ is called \emph{B\"{u}chi
  arithmetic} and is well-known to be automatic. Here, we show that
$(\BS(1,q),(\xrightarrow{g})_{g \in F},(\tostar[g])_{g \in F})$ can be
interpreted in a slight extension of B\"{u}chi arithmetic that is
still automatic.  From \cref{fo-cplus-bs}, we can derive a stronger
statement, which clearly implies decidability of $\SAT^\pm(\BS(1,q))$:
\begin{theorem}\label{fo-cplusminus-bs}
  The first-order theory of $\mathcal{C}^\pm(\BS(1,q))$ is decidable.
\end{theorem}
Indeed, since $\BS(1,q)$ is torsion-free, we can express the predicate $\bot_{g,h}$ using
universal quantification: We have $(g_1,g_2)\bot_{g,h} (h_1,h_2)$ if
and only if $g_1\tostar[g] g_2$ and $h_1\tostar[h] h_2$
and
\[ \forall f,f'\in\BS(1,q)\colon \left(g_1\tostar[g] f \wedge f\tostar[g] g_2 \wedge h_1\tostar[h] f' \wedge f'\tostar[h] h_2 \right) \to f\ne f'. \]

\section{Conclusion}\label{sec:conclusion}
We have shown that for infinite groups $H$, the problem $\KP(G\wr H)$
is decidable if and only if $\ExpEq(G)$ is decidable and either
(i)~$G$ is abelian and $\KP^+(H)$ is decidable or (ii)~$G$ is
non-abelian and $\KP^{\pm}(H)$ is decidable. This reduces the study of
decidablity of $\KP(G\wr H)$ to decidability questions about the
factors $G$ and $H$.

However, we leave open whether there is a group $H$ where $\KP^+(H)$
is decidable, but $\KP^{\pm}(H)$ is undecidable.
It is clear that both
are decidable for all groups in the class of knapsack-semilinear
groups. This class contains a large part of the groups for which
knapsack has been studied. For example, it contains graph
groups~\cite[Theorem 3.11]{LohreyZ18} and hyperbolic
groups~\cite[Theorem 8.1]{Loh19hyp}.  Moreover, knapsack-semilinearity
is preserved by a variety of constructions: This includes wreath
products~\cite[Theorem 5.4]{GanardiKLZ18}, graph products~\cite{FL19},
free products with amalgamation and HNN-extensions over finite
identified subgroups~\cite{FL19}, and taking finite-index
overgroups~\cite{FL19}. Moreover, the groups $H_3(\Z)$ and $\BS(1,q)$
for $q\ge 2$ are also unable to distinguish $\KP^+$ and $\KP^{\pm}$:
We have shown here that $\KP^+$ is undecidable in $H_3(\Z)$ and
$\KP^{\pm}$ is decidable in $\BS(1,q)$. To the best of the authors'
knowledge, among the groups for which knapsack is known to be
decidable, this only leaves $\BS(p,q)$ for $p,q$ coprime, and
$G\wr \BS(1,q)$ (with decidable $\ExpEq(G)$) as candidates to
distinguish $\KP^+$ and $\KP^{\pm}$.

\defbibheading{bibliography}[\refname]{\section*{#1}}
\printbibliography

\clearpage
\appendix
\section{Proofs from Section~\ref{sec:results}}\label{appendix-results}

\subsection{Proof of Theorem \ref{thm:SAT-KP}}\label{appendix-SAT-KP}
The goal of this section is to show that $\SAT^\pm(G)$ is effectively equivalent to $\KP^\pm(G)$ and $\SAT^+(G)$ is effectively equivalent to $\KP^+(G)$ for any finitely generated group $G$. We begin with the equivalence of the more general problems. The first direction is shown in the following lemma:
\begin{lemma}\label{lem:KP-SAT}
For any finitely generated group $G$ it holds that if $\SAT^\pm(G)$ is decidable, then $\KP^\pm(G)$ is decidable as well.
\end{lemma}
\begin{proof}
Let $(E = \alpha_1 \cdots \alpha_n \alpha_{n+1},L,D)$ be a $\KP^\pm(G)$-instance with $\alpha_{n+1}$ a constant and variables in $X = \{x_1,\dots,x_n\}$. We write $g_i := \gamma(\alpha_i)$ for all $i \in [1,n+1]$ and define the following formula in $\mathcal{F}^\pm$:
\begin{equation*}
\begin{split}
\varphi := \exists y_0,\dots,y_n \colon & \bigwedge_{i \in P_E} y_{i-1} \tostar[g_i] y_i \wedge \bigwedge_{i \in Q_E \setminus \{n+1\}} y_{i-1} \xrightarrow{g_i} y_i \wedge y_n \xrightarrow{g_{n+1}} y_0 \wedge \\
& \bigwedge_{(i,j) \in L} y_i \xrightarrow{1} y_j \wedge \bigwedge_{(i,j) \in D} (y_{i-1},y_i) \bot_{g_i,g_j} (y_{j-1},y_j).
\end{split}
\end{equation*}
Let $\varphi(y_0,\dots,y_n)$ be the part of $\varphi$ without the existential quantifiers which means that $y_0,\dots,y_n$ are free variables in $\varphi(y_0,\dots,y_n)$. For an assignment $\mu \colon Y := \{y_0,\dots,y_n\} \to G$ we write $\mu \models \varphi(y_0,\dots,y_n)$ if $\varphi(y_0,\dots,y_n)$ evaluates to true when setting $y_i$ to $\mu(y_i)$ for all $i \in [0,n]$. 

We claim that $\sol_G(E,L,D) \neq \emptyset$ if and only if $\varphi(y_0,\dots,y_n)$ is satisfiable. For the first direction we assume that $\nu \in \sol_G(E,L,D)$ and let $\pi_{\nu,E} = \pi_1 \cdots \pi_{n+1}$ be the factorized walk induced by $\nu$ on $E$. We define the assignment $\mu \colon Y \to G$ such that $\mu(y_i) := \nu(\alpha_1 \cdots \alpha_i)$ for all $i \in [1,n]$ and $\mu(y_0) := 1$. Then $\mu(y_{i-1}) g_i^{\nu(x_i)} = \mu(y_i)$ for all $i \in P_E$ and $\mu(y_{i-1}) g_i = \mu(y_i)$ for all $i \in Q_E \setminus \{n+1\}$. Moreover, since $\nu(E) = 1$, it holds that $\mu(y_n) g_{n+1} = \mu(y_0)$. Since $\nu$ fulfills the loop constraints in $L$, we have that $\mu(y_i) = \mu(y_j)$ for all $(i,j) \in L$. For all $(i,j) \in D$ we have that $\pi_i$ and $\pi_j$ are disjoint and therefore $(\mu(y_{i-1}),\mu(y_i)) \bot_{g_i,g_j} (\mu(y_{j-1}),\mu(y_j))$ is fulfilled. Thus, $\mu \models \varphi(y_0,\dots,y_n)$.

For the other direction we assume that $\mu \colon Y \to G$ such that $\mu \models \varphi(y_0,\dots,y_n)$. Then we define the valuation $\nu \in \N^X$ such that $\mu(y_{i-1}) g_i^{\nu(x_i)} = \mu(y_i)$ and $\nu(x_i)$ is minimal with this property for all $i \in P_E$. This can be computed by trying all values for $\nu(x_i)$ iteratively since $\mu(y_{i-1}) \tostar[g_i] \mu(y_i)$ evaluates to true. As
\[\bigwedge_{i \in P_E} \mu(y_{i-1}) \tostar[g_i] \mu(y_i) \wedge \bigwedge_{i \in Q_E \setminus \{n+1\}} \mu(y_{i-1}) \xrightarrow{g_i} \mu(y_i) \wedge \mu(y_n) \xrightarrow{g_{n+1}} \mu(y_0)\]
is fulfilled, we have that $\mu(y_0) \nu(\alpha_1) \cdots \nu(\alpha_n) \nu(\alpha_{n+1}) = \mu(y_0)$ and therefore $\nu(E) = 1$. Let $\pi_{\nu,E} = \pi_1 \cdots \pi_{n+1}$ be the factorized walk induced by $\nu$ on $E$. Since $\mu(y_i) = \mu(y_j)$ for all $(i,j) \in L$, it follows that $\mu(y_0) \nu(\alpha_1) \cdots \nu(\alpha_i) = \mu(y_0) \nu(\alpha_1) \cdots \nu(\alpha_j)$, which means that $\pi_{i+1} \cdots \pi_j$ is a loop for all $(i,j) \in L$. Moreover, since $(\mu(y_{i-1}),\mu(y_i)) \bot_{g_i,g_j} (\mu(y_{j-1}),\mu(y_j))$ is fulfilled for all $(i,j) \in D$, the minimality of $\nu(x_i)$ and $\nu(x_j)$ if $i,j \in P_E$ implies that the walks $(\mu(y_0) \nu(\alpha_1 \cdots \alpha_{i-1}) g_i^k)_{0 \leq k \leq \nu(x_i)}$ and $(\mu(y_0) \nu(\alpha_1 \cdots \alpha_{j-1}) g_j^\ell)_{0 \leq \ell \leq \nu(x_j)}$ are disjoint. These walks are also disjoint if $i \in Q_E$ or $j \in Q_E$ by setting $\nu(x_i) := 1$ or $\nu(x_j) := 1$. Therefore, $\pi_i$ and $\pi_j$ are disjoint for all $(i,j) \in D$. Thus, $\nu \in \sol_G(E,L,D)$.
\end{proof}

The reduction from $\SAT^\pm(G)$ to $\KP^\pm(G)$ is established by the next lemma.
\begin{lemma}\label{lem:SAT-KP}
For any finitely generated group $G$ it holds that if $\KP^\pm(G)$ is decidable, then $\SAT^\pm(G)$ is decidable as well.
\end{lemma}
\begin{proof}
Let $\varphi := \exists y_1,\dots,y_n \psi \in \mathcal{F}^\pm$ be a formula in prenex normal form where $\psi$ is quantifier-free with variables $Y = \{y_1,\dots,y_n\}$. If we replace the atoms of $\psi$ by variables and regard the resulting formula as a formula in propositional logic, we can compute all satisfying assignments $\mu_1,\dots,\mu_m$ by trying all combinations of truth assignments of the variables. Then we can write
\[\varphi \equiv \bigvee_{i=1}^m \exists y_1,\dots,y_n \bigwedge_{j=1}^{c_i} a_{i,j}\]
where $a_{i,1},\dots,a_{i,c_i}$ are the atoms of $\psi$ that are set to true in $\mu_i$ for all $i \in [1,m]$. We consider each disjunct separately and write it as
\[\exists y_1,\dots,y_n \bigwedge_{j=1}^{c} a_j.\]
We replace all atoms of the form $a_j = (g_1,g_2) \bot_{g,h} (h_1,h_2)$ by the conjunction $g_1 \tostar[g] g_2 \wedge h_1 \tostar[h] h_2$ and write the resulting formula as
\[\exists y_1,\dots,y_n \bigwedge_{j=1}^{c^\prime} b_j.\]
Furthermore, we define the set
\[B := \{((g_1,g_2),(h_1,h_2)) \mid \exists j \in [1,c] \colon a_j = (g_1,g_2) \bot_{g,h} (h_1,h_2)\}.\]

Let $b_j = s_j \xrightarrow{t_j} e_j$ or $b_j = s_j \tostar[t_j] e_j$ with $s_j,e_j \in Y$ and $t_j \in G$ for all $j \in [1,c^\prime]$. Without loss of generality we assume that for all $j,k \in [1,c^\prime]$ it holds that $s_j \neq s_k$ or $e_j \neq e_k$. We define the graph $\Gamma := (Y,\mathcal{E}^1,\mathcal{E}^\ast,t)$ with vertices $Y$, two sorts of edges
\[\mathcal{E}^1 := \{(s_j,e_j) \mid j \in [1,c^\prime] \wedge b_j =  s_j \xrightarrow{t_j} e_j\}\]
and
\[\mathcal{E}^\ast := \{(s_j,e_j) \mid j \in [1,c^\prime] \wedge b_j =  s_j \tostar[t_j] e_j\}\]
and edge labeling $t \colon \mathcal{E} := \mathcal{E}^1 \cup \mathcal{E}^\ast \to G$ such that $t(s_j,e_j) := t_j$ for all $j \in [1,c^\prime]$. For any subset of edges $\mathcal{S} \subseteq \mathcal{E}$ we write
\[\mathcal{S}^{-1} := \{(v,u) \mid (u,v) \in \mathcal{S}\}\]
to denote the set of reverse edges and $\mathcal{S}^{\pm 1} := \mathcal{S} \cup \mathcal{S}^{-1}$.

Let $C \subseteq Y$ be an undirected connected component of $\Gamma$ and $u \in C$. We interpret $u$ as initial vertex and represent all other vertices in $C$ by a path starting with $u$. Consider an edge $(v,w) \in \mathcal{E} \cap C^2$ that lies in the connected component $C$. We choose an undirected path from $u$ to $v$ and denote it by a tuple $(p_1,\dots,p_\ell)$ with $p_k \in \mathcal{E}^{\pm1}$ for all $k \in [1,\ell]$. We now define a knapsack expression that follows the path and the edge $(v,w)$ to reach $w$ and then goes back to $u$. For all $k \in [1,\ell]$ we define
\[\alpha_k := \begin{cases}
t(p_k)^{x_k}, & \text{if } p_k \in {\mathcal{E}^\ast}^{\pm 1} \\
t(p_k), & \text{otherwise}
\end{cases}\]
where we extend the edge labeling to reverse edges by setting
\[t(p_k) := \begin{cases}
t(p_k), & \text{if } p_k \in \mathcal{E} \\
t(p_k^{-1})^{-1}, & \text{otherwise.}
\end{cases}\]
To follow the edge $(v,w)$ we let
\[\alpha_{\ell+1} := \begin{cases}
t(v,w)^{x_{\ell+1}}, & \text{if } (v,w) \in \mathcal{E}^\ast \\
t(v,w), & \text{otherwise.}
\end{cases}\]
To walk back to $u$ we define
\[\alpha_{\ell+2} := \begin{cases}
(t(v,w)^{-1})^{x_{\ell+2}}, & \text{if } (v,w) \in \mathcal{E}^\ast \\
t(v,w)^{-1}, & \text{otherwise}
\end{cases}\]
and
\[\alpha_{\ell+2+k} := \begin{cases}
(t(p_{\ell+1-k})^{-1})^{x_{\ell+2+k}}, & \text{if } p_{\ell+1-k} \in {\mathcal{E}^\ast}^{\pm 1} \\
t(p_{\ell+1-k})^{-1}, & \text{otherwise}
\end{cases}\]
for all $k \in [1,\ell]$. We then define the knapsack expression $E_{v,w} := \alpha_1 \cdots \alpha_{2\ell+2}$ and loop constraint $L_{v,w} := \{(0,2\ell+2)\}$. If we do this for every edge lying in $C$ we obtain the knapsack expression
\[E_C := \prod_{(v,w) \in \mathcal{E} \cap C^2} E_{v,w}\]
where we make the indices continuous.

Let $\ell_{v,w}$ be the adjusted index $\ell$ in $E_{v,w}$ for all $(v,w) \in \mathcal{E} \cap C^2$. For every $v \in C$ we define the set of indices
\[I_v := \{\ell_{v,w} \mid (v,w) \in \mathcal{E}\} \cup \{\ell_{w,v}+1 \mid (w,v) \in \mathcal{E}\}.\]
We write $I_v = \{\ell_1,\dots,\ell_r\}$ with $\ell_1 < \dots < \ell_r$ and let $L_v := \{(\ell_k,\ell_{k+1}) \mid 1 \leq k < r\}$. Intuitively, the loop constraints in $L_v$ ensure that all edges incident to $v$ start or end at the same point. We can now define the set of loop constraints
\[L_C := \bigcup_{(v,w) \in \mathcal{E} \cap C^2} L_{v,w} \cup \bigcup_{v \in C} L_v\]
where we adjust the indices in $L_{v,w}$ properly.

If we do this for all undirected connected components $C_1,\dots,C_s$ of $\Gamma$ that have size greater than one, we obtain the $\KP^+(G)$-instance
\[(E := E_{C_1} \cdots E_{C_s}, L := L_{C_1} \cup \dots \cup L_{C_s})\]
where we adjust the indices and $\ell_{v,w}$ properly. We define the corresponding disjointness constraints
\[D := \{(\ell_{g_1,g_2}+1,\ell_{h_1,h_2}+1) \mid ((g_1,g_2),(h_1,h_2)) \in B\}.\]
Let $(E_1,L_1,D_1),\dots,(E_m,L_m,D_m)$ be the resulting $\KP^\pm(G)$-instances for all disjuncts. We claim that $\varphi$ is satisfiable if and only if $\bigcup_{i=1}^m \sol_G(E_i,L_i,D_i) \neq \emptyset$.

For the first direction let $\varphi_i(y_1,\dots,y_n) := \bigwedge_{j=1}^{c_i} a_{i,j}$ and assume that $\mu \colon Y \to G$ is a satisfying assignment of $\varphi_i$ for some $i \in [1,m]$. We write $E_i = \alpha_1 \cdots \alpha_d$ and by definition every power $\alpha_j$ with $j \in P_{E_i}$ has base $t(y_k,y_\ell)$ for some $k,\ell \in [1,n]$ with $(y_k,y_\ell) \in {\mathcal{E}^\ast}^{\pm 1}$. We define the valuation $\nu \in \N^X$ such that for all $j \in P_{E_i}$ where $\alpha_j$ has base $t(y_k,y_\ell)$ for some $k,\ell \in [1,n]$  with $(y_k,y_\ell) \in {\mathcal{E}^\ast}^{\pm 1}$ it holds that $\mu(y_k) t(y_k,y_\ell)^{\nu(x_j)} = \mu(y_\ell)$ and $\nu(x_j)$ is minimal with this property. Note that $\nu$ can be computed by trying all values iteratively since the construction of $E_i$ implies that $\mu(y_k) \tostar[t(y_k,y_\ell)] \mu(y_\ell)$ is fulfilled for every power $\alpha_j$ with base $t(y_k,y_\ell)$ where $(y_k,y_\ell) \in {\mathcal{E}^\ast}^{\pm 1}$ as $\mu \models \varphi_i(y_1,\dots,y_n)$. It follows that $\nu(E_i) = 1$ and $\nu$ fulfills the loop constraints in $L_i$ since variables of powers with equal or inverse base are set to the same value. Let $\pi_{\nu,E_i} = \pi_1 \cdots \pi_d$ be the factorized walk induced by $\nu$ on $E_i$. Since
\[(\mu(g_1),\mu(g_2)) \bot_{t(g_1,g_2),t(h_1,h_2)} (\mu(h_1),\mu(h_2))\]
is fulfilled for all $((g_1,g_2),(h_1,h_2)) \in B$ and $\gamma(\alpha_{\ell_{g_1,g_2}+1}) = t(g_1,g_2)$ and $\gamma(\alpha_{\ell_{h_1,h_2}+1}) = t(h_1,h_2)$, the minimality of $\nu(x_{\ell_{g_1,g_2}+1})$ and $\nu(x_{\ell_{h_1,h_2}+1})$, where we set $\nu(x_j) :=1$ if $j \in Q_{E_i}$, implies that $\pi_{\ell_{g_1,g_2}+1}$ and $\pi_{\ell_{h_1,h_2}+1}$ are disjoint. Thus, $\nu \in \sol_G(E_i,L_i,D_i)$.

For the other direction we assume that $\nu \in \sol_G(E_i,L_i,D_i)$ for some $i \in [1,m]$. Let $E_i = \alpha_1 \cdots \alpha_d$ and $\varphi_i(y_1,\dots,y_n) := \bigwedge_{j=1}^{c_i} a_{i,j}$. To show that $\varphi_i(y_1,\dots,y_n)$ is satisfiable, we define the assignment $\mu \colon Y \to G$ such that for all $v \in Y$ with $v$ incident to an edge of $\Gamma$ it holds that
\[\mu(v) := \begin{cases}
\prod_{k=1}^{\ell_{v,w}} \nu(\alpha_k), & \text{if } (v,w) \in \mathcal{E} \text{ for some } w \in Y \\
\prod_{k=1}^{\ell_{w,v}+1} \nu(\alpha_k), & \text{if } (w,v) \in \mathcal{E} \text{ for some } w \in Y
\end{cases}\]
where $\ell_{v,w}$ is the adjusted index in $E_i$. For every $v \in Y$ that is not incident to any edge of $\Gamma$ we set $\mu(v)$ to an arbitrary value of $G$. The loop constraints in $L_i$ ensure that the assignment $\mu$ is well-defined. By definition of $E_i$ it follows that $\mu \models b_j$ for all $j \in [1,c^\prime]$. Since we add for all $j \in [1,c_i]$ with $a_{i,j} = (g_1,g_2) \bot_{g,h} (h_1,h_2)$ the atoms $b_k = g_1 \tostar[g] g_2$ and $b_\ell = h_1 \tostar[h] h_2$ for some $k,\ell \in [1,c^\prime]$, the disjointness constraints in $D_i$ imply that $\mu \models a_{i,j}$ for all $j \in [1,c_i]$. Thus, $\mu \models \varphi_i(y_1,\dots,y_n)$.
\end{proof}

We show next that $\SAT^+(G)$ is effectively equivalent to $\KP^+(G)$ for any finitely generated group $G$. The first direction is shown in the following lemma:
\begin{lemma}\label{lem:pos-KP-SAT}
For any finitely generated group $G$ it holds that if $\SAT^+(G)$ is decidable, then $\KP^+(G)$ is decidable as well.
\end{lemma}
\begin{proof}
We can copy the proof of \Cref{lem:KP-SAT} by setting $D := \emptyset$.
\end{proof}

The reduction from $\SAT^+(G)$ to $\KP^+(G)$ is established by the next lemma.
\begin{lemma}\label{lem:pos-SAT-KP}
For any finitely generated group $G$ it holds that if $\KP^+(G)$ is decidable, then $\SAT^+(G)$ is decidable as well.
\end{lemma}
\begin{proof}
We can copy the proof of \Cref{lem:SAT-KP} by assuming that $\varphi \in \mathcal{F}^+$.
\end{proof}

\subsection{Proof of Corollary \ref{application-virtually-nilpotent}}\label{appendix-application-virtually-nilpotent}

Although we will not refer to this definition in the proof, we include
a definition of nilpotent groups for completeness. For a group $G$, we
define its \emph{lower central series} as the subgroups
$G_1,G_2,\ldots$ with $G_1=G$ and $G_{i+1}=[G_i,G]$ for $i\ge 1$.
Then, $G$ is \emph{nilpotent} if there is a number $n\ge 1$ with
$G_n=\{1\}$.

\begin{proof}[Proof of \cref{application-virtually-nilpotent}]
  Suppose $\KP(G\wr H)$ is decidable. Since $H$ is infinite, we know
  that $\ExpEq(G)$ must be decidable~\cite[Proposition 3.1,
  Proposition 5.1]{GanardiKLZ18}. Towards a contradiction, assume that
  $H$ is not virtually abelian.  As a finitely generated virtually
  nilpotent group, $H$ contains a finite-index nilpotent subgroup $K$
  that is also torsion-free~\cite[Theorem 17.2.2]{KaMe1979}. Since $H$
  is not virtually abelian, $K$ cannot be abelian. Since every
  non-abelian torsion-free nilpotent group has $H_3(\Z)$ as a subgroup
  (see, for example, the proof of \cite[Theorem 12]{holt2005groups}),
  we know that $H_3(\Z)$ is a subgroup of $H$.  Hence, $\KP(G\wr H)$
  is undecidable by \cref{application-heisenberg}, which is a
  contradiction.

  Conversely, suppose $H$ is virtually abelian and $\ExpEq(G)$ is
  decidable.  Since $H$ is virtually abelian, it is
  knapsack-semilinear~\cite[Theorem~7.1]{FL19}.  Therefore, since
  $\ExpEq(G)$ is decidable, decidability of $\KP(G\wr H)$ is shown in
  \cite[Theorem 5.3]{GanardiKLZ18}.
\end{proof}

\subsection{Proof of Corollary \ref{application-magnus}}\label{appendix-application-magnus}

\begin{proof}[Proof of \cref{application-magnus}]
  By the Magnus embedding theorem~\cite[Lemma]{Mag39}, the group
  $F/[N,N]$ embeds in $\Z^r\wr (F/N)$, where $r$ is the rank of
  $F$. By \cref{main-result}, decidability of $\IKPp(F/N)$
  implies decidability of $\KP(\Z^r\wr (F/N))$. Finally, for any
  $G$, $\IKPp(G)$ is a special case of $\ExpEq(G)$.
\end{proof}

\section{Proofs from Section~\ref{sec:wr-kp}}\label{appendix-wr-kp}
\subsection{The modified intersection knapsack problem}\label{app:MKP}
To simplify the constructions in the proofs from \Cref{sec:wr-kp}, we use slight variations of the problems $\KP^{\pm}(H)$ and $\KP^+(H)$.

Let $E = \alpha_1 \cdots \alpha_n$ be a knapsack expression over $G$. For every $i \in [1,n]$ and $\nu \in \mathbb{N}^X$ we define
\[S_E^\nu(i) := \{\nu(\alpha_1 \cdots \alpha_{i-1}) \gamma(\alpha_i)^k \mid 0 \leq k \leq \nu(x_i)-1\}\]
if $i \in P_E$ and
\[S_E^\nu(i) := \{\nu(\alpha_1 \cdots \alpha_{i-1})\}\]
if $i \in Q_E$. Intuitively, $S_E^\nu(i)$ is the set of points visited by the ray associated to $\alpha_i$ under the valuation $\nu$ where we leave out the last point.
\begin{defn} The {\em modified intersection knapsack problem} $\MKP^\pm(G)$ over $G$ is defined as follows:
\begin{description}
\item[Given] a knapsack expression $E$ over $G$, a set $L \subseteq [0,n]^2$ of loop constraints,
and a set $D \subseteq [1,n]^2$ of disjointness constraints.
\item[Question] Is there a valuation $\nu \in \N^X$ with factorized walk $\pi_{\nu,E} = \pi_1 \dots \pi_n$ induced by $\nu$ on $E$ such that the following conditions are fulfilled:
\begin{itemize}
\item $\nu(E) = 1$
\item $\pi_{i+1} \dots \pi_j$ is a loop for all $(i,j) \in L$
\item $S_E^\nu(i) \cap S_E^\nu(j) = \emptyset$ for all $(i,j) \in D$?
\end{itemize}
\end{description}
The {\em positive modified intersection knapsack problem} $\MKP^+(G)$ over $G$
is the restriction of $\MKP^\pm(G)$ to instances where $D = \emptyset$. As before, let $\sol_G(E,L,D)$ (resp. $\sol_G(E,L)$) be the set of solutions of the $\MKP^\pm(G)$-instance $(E,L,D)$ (resp. $\MKP^+(G)$-instance $(E,L)$) over $G$.
\end{defn}
Note that the restricted problems $\KP^+(G)$ and $\MKP^+(G)$ are identical and the only difference between $\KP^\pm(G)$ and $\MKP^\pm(G)$ is that the disjointness constraints of $\MKP^\pm(G)$-instances ignore the last point of walks. The equivalence of $\KP^\pm(G)$ and $\MKP^\pm(G)$ is established by the following lemma:

\begin{lemma}\label{lem:MKP}
For any finitely generated group $G$ we have that $\KP^\pm(G)$ and $\MKP^\pm(G)$ are effectively equivalent.
\end{lemma}
\begin{proof}
We first reduce $\KP^\pm(G)$ to $\MKP^\pm(G)$. Let $(E = \alpha_1 \cdots \alpha_n,L,D)$ be a $\KP^\pm(G)$-instance. We define the knapsack expression
\[E^\prime := \beta_1 \cdots \beta_{2n} := \alpha_1 \cdot 1 \cdots \alpha_n \cdot 1\]
with loop constraints $L^\prime := \{(2i-1,2j-1) \mid (i,j) \in L\}$ and disjointness constraints
\[D^\prime := \bigcup_{(i,j) \in D} \{(2i-1,2j-1), (2i-1,2j), (2i,2j-1), (2i,2j)\}.\]
We regard $(E^\prime,L^\prime,D^\prime)$ as $\MKP^\pm(G)$-instance. Note that with the added 1's we can ensure that $D^\prime$ considers also the last points of the disjointness constraints defined in $D$.

We show that $\sol_G(E,L,D) = \sol_G(E^\prime,L^\prime,D^\prime)$. Let $\nu \in \N^X$ be a valuation and let $\pi_{\nu,E} = \pi_1 \cdots \pi_n$ be the factorized walk induced by $\nu$ on $E$. Clearly, it holds that $\nu(E) = 1$ if and only if $\nu(E^\prime) = 1$ and $\nu$ fulfills the loop constraints in $L$ if and only if it fulfills the loop constraints in $L^\prime$. We now consider the disjointness constraints. Let $g_i := \gamma(\alpha_i)$ for all $i \in [1,n]$ and $\nu(x_i) := 1$ if $i \in Q_E$. For all $(i,j) \in D$ we have that $\pi_i$ and $\pi_j$ are disjoint if and only if
\[\{\nu(\alpha_1 \cdots \alpha_{i-1}) g_i^k \mid 0 \leq k \leq \nu(x_i)\} \cap \{\nu(\alpha_1 \cdots \alpha_{j-1}) g_j^k \mid 0 \leq k \leq \nu(x_j)\} = \emptyset\]
which holds if and only if
\[\{\nu(\alpha_1 \cdots \alpha_{i-1}) g_i^k \mid 0 \leq k \leq \nu(x_i)-1\} \text{ and } \{\nu(\alpha_1 \cdots \alpha_{i-1}) g_i^{\nu(x_i)}\}\]
are disjoint to
\[\{\nu(\alpha_1 \cdots \alpha_{j-1}) g_j^k \mid 0 \leq k \leq \nu(x_j)\} \text{ and } \{\nu(\alpha_1 \cdots \alpha_{j-1}) g_j^{\nu(x_j)}\}\]
which in turn holds if and only if $S_{E^\prime}^\nu(2i-1)$ and $S_{E^\prime}^\nu(2i)$ are disjoint to $S_{E^\prime}^\nu(2j-1)$ and $S_{E^\prime}^\nu(2j)$. Thus, $\nu \in \sol_G(E,L,D)$ if and only if $\nu \in \sol_G(E^\prime,L^\prime,D^\prime)$.

We now reduce $\MKP^\pm(G)$ to $\KP^\pm(G)$. Let $(E = \alpha_1 \cdots \alpha_n,L,D)$ be an $\MKP^\pm(G)$-instance and $g_i := \gamma(\alpha_i)$ for all $i \in [1,n]$. Let $P \subseteq P_E$ be a set of powers whose variables will be set to 0. For all $j \in [1,n]$ we replace
\[\alpha_j \text{ by } \begin{cases}
g_j^{y_{i_{j,1}}} g_j =: \beta_{i_{j,1}} \beta_{i_{j,2}}, & \text{if } j \in P_E \setminus P \\
1 =: \beta_{i_{j,2}}, & \text{if } j \in P \\
1 \cdot g_j =: \beta_{i_{j,1}} \beta_{i_{j,2}}, & \text{if } j \in Q_E
\end{cases}\]
to get the knapsack expression $E_P$ and we write $E_P = \beta_1 \cdots \beta_r$ with variables in $Y := \{y_1,\dots,y_r\}$ by making indices continuous where we adjust $i_{j,1}$ and $i_{j,2}$ accordingly. We define the loop constraints $L_P := \{(i_{j,2},i_{k,2}) \mid (j,k) \in L\}$ and the disjointness constraints
\[D_P := \{(i_{j,1},i_{k,1}) \mid (j,k) \in D \wedge j,k \notin P\}.\]
We interpret $(E_P,L_P,D_P)$ as $\KP^\pm(G)$-instance. The idea is to split progressions at the last point such that the $\KP^\pm(G)$-instance ignores this point. The splitting is not possible if the variable is set to 0. Thus, we need to guess the the set of powers $P$ whose variables are set to 0 beforehand. It remains to show that $\sol_G(E,L,D) \neq \emptyset$ if and only if $\bigcup_{P \in P_E} \sol_G(E_P,L_P,D_P) \neq \emptyset$.

For the first direction let $\nu \in \sol_G(E,L,D)$. We define $P := \{i \in P_E \mid \nu(x_i)=0\}$ and the valuation $\nu_P \in \N^Y$ such that $\nu_P(y_{i_{j,1}}) := \nu(x_j)-1$ for all $j \in P_E \setminus P$. Let $\pi_{\nu_P,E_P} = \pi_1 \cdots \pi_r$ be the factorized walk induced by $\nu_P$ on $E_P$. By definition of $E_P$ it clearly holds that $\nu_P(E_P) = 1$ and $\nu_P$ fulfills all loop constraints in $L_P$. We now consider the disjointness constraints. Let $\nu_P(y_{i_{j,1}}) := 1$ for all $j \in Q_E$ and $h_i := \gamma(\beta_i)$ for all $i \in [1,r]$. For every $(j,k) \in D$ with $j,k \notin P$ we have that $S_E^\nu(j) \cap S_E^\nu(k) = \emptyset$. Therefore, it holds that
\[\{\nu_P(\beta_1 \cdots \beta_{i_{j,1}-1}) h_{i_{j,1}}^\ell \mid 0 \leq \ell \leq \nu_P(y_{i_{j,1}})\} \cap \{\nu_P(\beta_1 \cdots \beta_{i_{k,1}-1}) h_{i_{k,1}}^\ell \mid 0 \leq \ell \leq \nu_P(y_{i_{k,1}})\} = \emptyset\]
which implies that $\pi_{i_{j,1}}$ and $\pi_{i_{k,1}}$ are disjoint. Thus, $\nu_P \in \sol_G(E_P,L_P,D_P)$.

For the other direction let $\nu_P \in \sol_G(E_P,L_P,D_P)$ for some $P \subseteq P_E$. We define the valuation $\nu \in \N^X$ such that
\[\nu(x_j) := \begin{cases}
\nu_P(y_{i_{j,1}})+1, & \text{if } j \in P_E \setminus P \\
0, & \text{if } j \in P
\end{cases}\]
for all $j \in P_E$. Clearly, it holds that $\nu(E) = 1$ and $\nu$ fulfills all loop constraints in $L$. Let $\nu(x_i) := 1$ for all $i \in Q_E$ and $\pi_{\nu_P,E_P} = \pi_1 \cdots \pi_r$ be the factorized walk induced by $\nu_P$ on $E_P$. For every $(j,k) \in D$ with $j,k \notin P$ we have that $\pi_{i_{j,1}}$ and $\pi_{i_{k,1}}$ are disjoint. Therefore, it holds that
\[\{\nu(\alpha_1 \cdots \alpha_{j-1}) g_j^\ell \mid 0 \leq \ell \leq \nu(x_j)-1\} \cap \{\nu(\alpha_1 \cdots \alpha_{k-1}) g_k^\ell \mid 0 \leq \ell \leq \nu(x_k)-1\} = \emptyset\]
which implies that $S_E^\nu(j) \cap S_E^\nu(k) = \emptyset$. For $(j,k) \in D$ with $j \in P$ or $k \in P$ it holds that $\nu(x_j) = 0$ or $\nu(x_k) = 0$ and therefore $S_E^\nu(j) = \emptyset$ or $S_E^\nu(k) = \emptyset$ which implies that $S_E^\nu(j) \cap S_E^\nu(k) = \emptyset$. Thus, $\nu \in \sol_G(E,L,D)$.
\end{proof}

\subsection{Proof of Theorem~\ref{thm:normal}}\label{app:normal}
Let $P$ and $P^\prime$ be two potentially equal decision problems defined so far.
Let $S = \{I_1,\dots,I_s\}$ be a finite set of instances of $P$ and $S^\prime = \{I_1^\prime,\dots,I_t^\prime\}$ be a finite set of instances of $P^\prime$. We say that $S$ is equivalent to $S^\prime$ if $\bigcup_{i=1}^s \sol_{P}(I_i) \neq \emptyset$ if and only if $\bigcup_{i=1}^t \sol_{P^\prime}(I_i^\prime) \neq \emptyset$. Here, $\sol_P$ and $\sol_{P^\prime}$ denote the set of solutions of an instance of the respective problem. We define the equivalence also directly on instances by assuming singleton sets.

We say that a knapsack expression $E = \alpha_1 \cdots \alpha_n$ is \textit{torsion-free} if for all $i \in P_E$ it holds that $\sigma(\gamma(\alpha_i)) = 1$ or $\sigma(\gamma(\alpha_i))$ has infinite order.
\begin{lemma}\label{lem:torsion}
For any knapsack expression one can effectively construct an equivalent finite set of torsion-free knapsack expressions.
\end{lemma}
\begin{proof}
We use the ideas of the proof of Lemma 7.1 from \cite{GKLZ17}. First note that by conjugation we can eliminate constants in a knapsack expression $E$ and assume that $E = g_1^{x_1} \cdots g_d^{x_n} g$ where $g_1,\dots,g_n,g \in G \wr H$. Let $i \in \{1,\dots,n\}$ such that $\sigma(g_i) \neq 1$ and $\ord(\sigma(g_i)) = q < \infty$. Since $\KP(H)$ is decidable we can compute $q$ as follows. We first check if $\sigma(g_i)^x \sigma(g_i) = 1$ has a solution and if so, we try every value for $x$ starting with 0 until we find a solution which is then $q-1$.

We then construct the expression
\[E_r^{\prime\prime} = g_1^{x_1} \cdots g_{i-1}^{x_{i-1}} (g_i^q)^{x_i} g_i^r g_{i+1}^{x_{i+1}} \cdots g_n^{x_n} g\]
and from that the knapsack expression
\[E_r^\prime = g_1^{x_1} \cdots g_{i-1}^{x_{i-1}} (g_i^q)^{x_i} (g_i^r g_{i+1} g_i^{-r})^{x_{i+1}} \cdots (g_i^r g_n g_i^{-r})^{x_n} g_i^r g\]
for all $r \in [0,q-1]$. The idea is to write exponents as multiple of the order of the base with remainder. We then shift the constant factor for the remainder via conjugation to the end of the expression. Note that $E_r^\prime$ has one non-trivial torsion element less than $E$ since $\sigma(g_i^q) = 1$ and conjugation by $g_i^r$ does not change the orders of the elements $g_{i+1},\dots,g_n$. Clearly, it holds that $\sol_{G \wr H}(E_r^{\prime\prime}) = \sol_{G \wr H}(E_r^\prime)$ for all $r \in [0,q-1]$.

If $\nu \in \N^X$ is a solution of $E$, then for $r := \nu(x_i) \text{ mod } q$ we get a solution $\nu^\prime \in \N^X$ of $E_r^\prime$ by setting 
\[\nu^\prime(x_j) := \begin{cases}
s, & \text{if } j = i \\
\nu(x_j), & \text{otherwise}
\end{cases}\]
for all $j \in [1,n]$ where $\nu(x_i) = s q + r$. Conversely, if $\nu^\prime \in \N^X$ is a solution of $E_r^\prime$ for some $r \in [0,q-1]$, then $\nu \in \N^X$ with
\[\nu(x_j) := \begin{cases}
q \nu^\prime(x_i) + r, & \text{if } j = i \\
\nu^\prime(x_j), & \text{otherwise}
\end{cases}\]
for all $j \in [1,n]$ is a solution of $E$. Thus, it holds that $\sol_{G \wr H}(E) \neq \emptyset$ if and only if $\bigcup_{r=0}^{q-1} \sol_{G \wr H}(E_r^\prime) \neq \emptyset$.

Repeating this process for all $E_r^\prime$ until we get torsion-free knapsack expressions $E_1,\dots,E_t$ yields the lemma.
\end{proof}

A knapsack expression $E = \alpha_1 \cdots \alpha_n \alpha_{n+1}$ is in $GH$\textit{-form} if for all $i \in P_E$ it holds that $\sigma(\gamma(\alpha_i)) = 1$ or $\gamma(\alpha_i) \in GH$ and for all $i \in Q_E \setminus \{n+1\}$ it holds that $\alpha_i \in H$.

To do the transformation into $GH$-form, we need an order on the elements in the support of some atom of $E$. Let $h \in H$ be a torsion-free element. We define the binary relation $\preceq_h$ on $H$ as in \cite{GKLZ17}. For $h^\prime, h^{\prime \prime} \in H$ we write $h^\prime \preceq_h h^{\prime \prime}$ if there is a $k \geq 0$ such that $h^\prime = h^k h^{\prime \prime}$. Clearly, $\preceq_h$ is a partial order since $h$ is torsion-free. Moreover, since $\KP(H)$ is decidable, we can decide with a knapsack instance over $H$ whether $h^\prime \preceq_h h^{\prime \prime}$.

To multiply elements $a_i$ for $i \in I$ in a certain order, we write for a finite linearly ordered set $(I = \{i_1,\dots,i_m\},\leq)$ with $i_1 < \dots < i_m$ the product $\prod_{j=1}^m a_{i_j}$ as $\prod_{i \in I}^\leq a_i$. The following lemma is shown in \cite{GKLZ17}.

\begin{lemma}\label{lem:order}
Let $g \in G \wr H$ such that $\ord(\sigma(g)) = \infty$ and let $h \in H$ and $m \in \N$. Moreover, let $F = \supp(g) \cap \{\sigma(g)^{-i} h \mid i \in [0,m-1]\}$. Then $F$ is linearly ordered by $\preceq_{\sigma(g)}$ and
\[\tau(g^m)(h) = \xsideset{}{^{\preceq_{\sigma(g)}}} \prod_{h^\prime \in F} \tau(g)(h^\prime).\]
\end{lemma}

Thus, $\preceq_{\sigma(g)}$ tells us how to evaluate $\tau(g^m)$ at a certain element of $H$. We use this to establish the $GH$-form for $E$.

\begin{lemma}\label{lem:GH}
For any torsion-free knapsack expression one can effectively construct an equivalent torsion-free $\HKP^+(G \wr H)$-instance in $GH$-form.
\end{lemma}
\begin{proof}
We use the idea of the proof of Lemma 29 from \cite{FGLZ20}. Let $u \in G \wr H$ with $\sigma(u)$ torsion-free and for $ h \in \supp(u)$ let $a_h := \tau(u)(h)$. We want to dissect $u^m$ such that every element in the support of $u$ yields a ray. For $h \in \supp(u)$ such a ray visits the points $\sigma(u)^k h$ for all $k \in [0,m-1]$. Note that if $h_1,h_2 \in \supp(u)$ and $\sigma(u)^{k_1} h_1 = \sigma(u)^{k_2} h_2$ for some $0 \leq k_1 \leq k_2 \leq m-1$, that is, the rays of $h_1$ and $h_2$ intersect and the ray of $h_1$ visits the intersection points first, then $h_1 \preceq_{\sigma(u)} h_2$.

We extend the partial order $\preceq_{\sigma(u)}$ to a linear order $\leq_{\sigma(u)}$ on $\supp(u)$. Then by \Cref{lem:order} for all $x \in \N$ it holds that
\[u^x = \Bigg( \xsideset{}{^{\leq_{\sigma(u)}}} \prod_{h \in \supp(u)} h (a_h h^{-1} \sigma(u) h)^x h^{-1} \sigma(u)^{-x}\Bigg) \sigma(u)^x.\]
Note that the part $h (a_h h^{-1} \sigma(u) h)^x$ writes $a_h$ at the points $\sigma(u)^k h$ for $k \in [0,x]$. We then go back with $h^{-1} \sigma(u)^{-x}$ to the beginning which is the starting point for the next element in $\supp(u)$. Finally, we walk with $\sigma(u)^x$ to the end of the progression since also the last factor of the product walks back to the beginning. As in knapsack expressions we cannot use the variable $x$ multiple times, we need loop constraints to ensure that we walk back and forth by the same distance.

Let $\supp(u) = \{h_1,\dots,h_\ell\}$ such that $h_1 \leq_{\sigma(u)} \dots \leq_{\sigma(u)} h_\ell$ and $a_i := a_{h_i}$. Then we can construct the following $\HKP^+(G \wr H)$-instance:
\begin{equation*}
\begin{split}
& \Bigg(\prod_{i=1}^\ell h_i (a_i h_i^{-1} \sigma(u) h_i)^{y_{4i-2}} h_i^{-1} (\sigma(u)^{-1})^{y_{4i}} \Bigg) \sigma(u)^{y_{4\ell+1}} \\
= & \Bigg(\prod_{i=1}^\ell \beta_{4i-3} \beta_{4i-2} \beta_{4i-1} \beta_{4i} \Bigg) \beta_{4\ell+1} =: E_u
\end{split}
\end{equation*}
where for all $j \in [1,4\ell+1]$ it holds that $\gamma(\beta_j) \in GH$ if $j \in P_{E_u}$ and $\beta_j \in H$ if $j \in Q_{E_u}$. We define the corresponding loop constraints
\begin{equation*}
\begin{split}
L_u := & \{(4i-4,4i) \mid 1 \leq i \leq \ell\} \cup \\
& \{(4i-1,4(i+1)-1) \mid 1 \leq i \leq \ell-1\} \cup \\
& \{(4\ell-1,4\ell+1)\}.
\end{split}
\end{equation*}
This means that for any solution $\nu^\prime \in \sol_{G \wr H}(E_u g_u,L_u)$, for some $g_u \in G \wr H$, it must hold that $\nu^\prime(y_{4i-2}) = \nu^\prime(y_{4i}) = \nu^\prime(y_{4\ell+1})$ for all $i \in [1,\ell]$ since $h_i^{-1} \sigma(u) h_i$ is torsion-free. Thus, for all $g_u \in G \wr H$ we have $\sol_{G \wr H}(u^x g_u) = \pi_{u^x}^{(g_u)}(\sol_{G \wr H}(E_u g_u,L_u))$ where we define the projection $\pi_{u^x}^{(g_u)}$ as
\begin{align*}
\pi_{u^x}^{(g_u)} \colon \sol_{G \wr H}(E_u g_u,L_u) &\to \sol_{G \wr H}(u^x g_u)\\
\nu^\prime &\mapsto 
\Bigg(\begin{aligned}
\nu \colon \{x\} &\to \N\\
x &\mapsto \nu^\prime(y_2)
\end{aligned}\Bigg).
\end{align*}
Moreover, since $\sigma(\gamma(\beta_{4i-2})) = h_i^{-1} \sigma(u) h_i$ and $\sigma(u)$ is torsion-free, it follows that $\sigma(\gamma(\beta_{4i-2}))$, $\sigma(\gamma(\beta_{4i}))$ and $\sigma(\gamma(\beta_{4\ell+1}))$ are torsion-free as well for all $i \in [1,\ell]$. Therefore, we have that $\sigma(\gamma(\beta_j))$ is torsion-free for all $j \in P_{E_u}$. Note that the factors $\beta_{4i-3}$ and $\beta_{4i-1}$ with $4i-3, 4i-1 \in Q_{E_u}$ are not torsion-free in general.

We now consider the whole torsion-free knapsack expression $E = \alpha_1 \cdots \alpha_n \alpha_{n+1}$. By conjugation we can eliminate the constants in $E$ and assume that $E = g_1^{x_1} \cdots g_d^{x_n} g$ with $g_1,\dots, g_n, g \in G \wr H$ which is still torsion-free as conjugation does not change the order of an element. Then we construct the following $\HKP^+(G \wr H)$-instance:
\[(E_{g_1} \cdots E_{g_n} g, L_{g_1} \cup \dots \cup L_{g_n})\]
where we choose continuous indices and variables $Y = Y_1 \dot\cup \dots \dot\cup Y_n$ such that $E_{g_i}$ has variables $Y_i$. Here we set $E_{g_i} := g_i^{x_i}$ and $L_{g_i} := \emptyset$ if $\sigma(g_i) = 1$. Note that since $E$ is torsion-free, we have that all $\sigma(g_i) \neq 1$ are torsion-free and therefore $E_{g_i}$ is well-defined. The lemma follows from the following observation:
\[\sol_{G \wr H}(E) = \pi(\sol_{G \wr H}(E_{g_1} \cdots E_{g_n} g, L_{g_1} \cup \dots \cup L_{g_n}))\]
with the projection $\pi$ defined for $\nu^\prime \in \sol_{G \wr H}(E_{g_1} \cdots E_{g_n} g, L_{g_1} \cup \dots \cup L_{g_n})$ as
\begin{align*}
\pi(\nu^\prime) \colon X &\to \N \\
x_i &\mapsto 
\begin{cases}
\pi_{g_i^{x_i}}^{(g_{g_i})}(\nu^\prime |_{Y_i})(x_i), & \text{if } \sigma(g_i) \neq 1\\
\nu^\prime(x_i), & \text{otherwise}
\end{cases}
\end{align*}
where $g_{g_i} := \nu^\prime(E_{g_1} \cdots E_{g_{i-1}})^{-1} \nu^\prime(E_{g_{i+1}} \cdots E_{g_n} g)^{-1}$ and $\nu^\prime |_{Y_i}$ denotes the restriction of $\nu^\prime$ to $Y_i$.
\end{proof}

In the next normalization step we deal with commensurable elements. For a knapsack expression $E = \alpha_1 \cdots \alpha_n$ let us define an equivalence relation $||$ on the set 
\[R_E = \{r \in P_E \mid \gamma(\alpha_r) \notin H \wedge \sigma(\gamma(\alpha_r)) \neq 1\}.\]
For $r_1, r_2 \in R_E$ we say $r_1 || r_2$ if $\sigma(\gamma(\alpha_{r_1}))$ and $\sigma(\gamma(\alpha_{r_2}))$ are commensurable. In the following we write $g_i := \gamma(\alpha_i)$ for all $i \in [1,n]$. 

\begin{lemma}\label{lem:class}
One can compute the $||$-classes for any knapsack expression $E = \alpha_1 \cdots \alpha_n$.
\end{lemma}
\begin{proof}
First note that $R_E$ can be computed since $\KP(H)$ is decidable. For each pair $(i,j) \in R_E^2$ check with $\KP(H)$-instances if $\sigma(g_i)^{x} \sigma(g_j)^{y} \sigma(g_i) = 1$ or $\sigma(g_i)^{x} (\sigma(g_j)^{-1})^{y} \sigma(g_i) = 1$ has a solution with $x,y \in \N$. If so, then there are $a,b \in \Z \setminus \{0\}$ such that $\sigma(g_i)^{a} = \sigma(g_j)^{b}$ which means that $i$ and $j$ are contained in the same $||$-class. If the instances do not have a solution, then $i$ and $j$ are in different $||$-classes.
\end{proof}

\begin{lemma}\label{lem:exp}
For any $||$-class $C$ of a knapsack expression $E = \alpha_1 \cdots \alpha_d$ one can compute natural numbers $e_c \neq 0$ for $c \in C$ such that $\sigma(g_{c_1})^{e_{c_1}} = \sigma(g_{c_2})^{e_{c_2}}$ or $\sigma(g_{c_1})^{e_{c_1}} = \sigma(g_{c_2})^{-e_{c_2}}$ for all $c_1,c_2 \in C$.
\end{lemma}
\begin{proof}
Let $C = \{i_1,\dots,i_m\}$ be a $||$-class with $i_1 < \dots < i_m$. We first compute $a_j,b_j \in \Z \setminus \{0\}$ with $\sigma(g_{i_j})^{a_j} = \sigma(g_{i_{j+1}})^{b_j}$ for all $j \in [1,m-1]$. To this end, we try all values for $x$ and $y$ in $\Z \setminus \{0\}$ until we find $a_j$ and $b_j$. This process terminates since $\sigma(g_{i_j})$ and $\sigma(g_{i_{j+1}})$ are commensurable.

Now we can define integers
\[e_j := \prod_{k=1}^{j-1} b_k \cdot \prod_{k=j}^m a_k\]
for all $j \in [1,m]$. Then for all $j \in [1,m-1]$ it holds that
\begin{equation*}
\begin{split}
\sigma(g_{i_j})^{e_j} & = \sigma(g_{i_j})^{\prod_{k=1}^{j-1} b_k \cdot \prod_{k=j}^m a_k} \\
& = \sigma(g_{i_j})^{a_j \cdot \prod_{k=1}^{j-1} b_k \cdot \prod_{k=j+1}^m a_k} \\
& = \sigma(g_{i_{j+1}})^{b_j \cdot \prod_{k=1}^{j-1} b_k \cdot \prod_{k=j+1}^m a_k} \\
& = \sigma(g_{i_{j+1}})^{e_{j+1}}.
\end{split}
\end{equation*}
Taking $|e_j|$ for all $j \in [1,m]$ yields the lemma.
\end{proof}

\begin{lemma}\label{lem:c-simplified}
For any torsion-free $\HKP^+(G \wr H)$-instance in $GH$-form one can effectively construct an equivalent finite set of c-simplified, torsion-free $\HKP^+(G \wr H)$-instances in $GH$-form.
\end{lemma}
\begin{proof}
Let $(E = \alpha_1 \cdots \alpha_n \alpha_{n+1},L)$ be a torsion-free $\HKP^+(G \wr H)$-instance in $GH$-form and write $g_i := \gamma(\alpha_i)$ for all $i \in [1,n+1]$. Let $g := g_i$ for some $i \in P_E$ with $g_i \notin H$ and $\sigma(g_i) \neq 1$. This means that $\sigma(g)$ is torsion-free since $E$ is torsion-free. Let $e_g \in \N \setminus \{0\}$ be the exponent from \Cref{lem:exp} corresponding to $g$. Note that $e_g$ can be computed since by \Cref{lem:class} we can effectively identify the $||$-class of $g$.

We first show that we can assume a slightly weaker property than to be c-simplified. We allow that for two elements $g_i,g_j \notin H$ with $i,j \in P_E$ such that $\sigma(g_i)$ and $\sigma(g_j)$ are commensurable it holds that $\sigma(g_i) = \sigma(g_j)$ or $\sigma(g_i) = \sigma(g_j)^{-1}$. To this end, we want to write $g^x$ as $g^{e_g y+r}$ for every remainder $r \in [0,e_g-1]$ but we have to make sure that the resulting $\HKP^+(G \wr H)$-instances are still in $GH$-form.

Let us construct the $\HKP^+(G \wr H)$-instance
\begin{equation*}
\begin{split}
F^{(g)} := & (\tau(g) \sigma(g)^{e_g})^{y_1} \sigma(g)^{-1} \sigma(g) (\sigma(g)^{-e_g})^{y_4} \sigma(g) \cdots \\
& (\tau(g) \sigma(g)^{e_g})^{y_{5(e_g-1)-4}} \sigma(g)^{-1} \sigma(g) (\sigma(g)^{-e_g})^{y_{5(e_g-1)-1}} \sigma(g) (\tau(g) \sigma(g)^{e_g})^{y_{5e_g-4}} \sigma(g)^{-1} \sigma(g) \\
= & \beta_1 \cdots \beta_{5e_g-2}
\end{split}
\end{equation*}
with loop constraints
\begin{equation*}
\begin{split}
J^{(g)} := & \{(5i-5,5i-1) \mid 1 \leq i \leq e_g-1\} \cup \\
& \{(5i-2,5(i+1)-3) \mid 1 \leq i \leq e_g-1\}.
\end{split}
\end{equation*}
Intuitively, this means that for all valuations $\nu$ we force that
\begin{align*}
(\sigma(g)^{e_g})^{\nu(y_{5i-4})} \cdot (\sigma(g)^{-e_g})^{\nu(y_{5i-1})} &= 1 \\
(\sigma(g)^{-e_g})^{\nu(y_{5i-1})} \cdot \sigma(g) \cdot (\sigma(g)^{e_g})^{\nu(y_{5i+1})} \sigma(g)^{-1} &= 1
\end{align*}
for all $i \in [1,e_g-1]$. Since $\sigma(g)$ is torsion-free, this implies that $\nu(y_{5i-4}) = \nu(y_{5i-1}) = \nu(y_{5 e_g-4})$ for all $i \in [1,e_g-1]$. Note that $F^{(g)}$ constitutes the part $g^{e_g y}$. The factor $(\tau(g) \sigma(g)^{e_g})^{y_1}$ visits powers of $\sigma(g)$ where the exponents are multiples of $e_g$ with offset 0. With $(\sigma(g)^{-e_g})^{y_4}$ we walk back to the beginning and set with $\sigma(g)$ the offset to 1. We then visit powers of $\sigma(g)$ where the exponents are multiples of $e_g$ with offset 1. We do this for every offset in $[0,e_g-1]$ to reach all the points of the progression associated to $g^{e_g y}$. The factors $\sigma(g)^{-1} \sigma(g)$ are only needed to define the loop constraints.

For the part of the remainder $g^r$ we construct the $\HKP^+(G \wr H)$-instance
\begin{equation*}
\begin{split}
G_r^{(g)} := & (\tau(g) \sigma(g)^{e_g})^{z_1} \sigma(g)^{-e_g} \sigma(g) \cdots (\tau(g) \sigma(g)^{e_g})^{z_{3r-2}} \sigma(g)^{-e_g} \sigma(g) \\
= & \gamma_1 \cdots \gamma_{3r}
\end{split}
\end{equation*}
with loop constraints
\[K_r^{(g)} :=  \{(3i-3,3i-1) \mid 1 \leq i \leq r\}\]
for all $r \in [0,e_g-1]$. Again since $\sigma(g)$ is torsion-free, this means intuitively that for every valuation $\nu$ we force that $\nu(z_{3i-2}) = 1$ for all $i \in [1,r]$. The idea of the construction of $G_r^{(g)}$ is the same as for $F^{(g)}$ but we set the exponents $z_i$ to 1.

We can now combine the two $\HKP^+(G \wr H)$-instances to obtain
\[(E_r^{(g)} := F^{(g)} \cdot G_r^{(g)}, L_r^{(g)} := J^{(g)} \cup K_r^{(g)})\]
for all $r \in [0,e_g-1]$ where we write $E_r^{(g)} = \delta_1 \cdots \delta_{5e_g-2+3r}$ and adjust the loop constraints in $J^{(g)}$ and $K_r^{(g)}$ accordingly.

Let $R_E = \{i_1,\dots, i_m\}$. If we replace every $g_{i_j}$ in $E$ by $E_r^{(g_{i_j})}$ and add the loop constraints $L_r^{(g_{i_j})}$ for all $j \in [1,m]$, we get the $\HKP^+(G \wr H)$-instance
\[(E_{r_1,\dots,r_m} := E_1 \cdots E_n g_{n+1}, L_{r_1,\dots,r_m} := L \cup L_{r_1}^{(g_{i_1})} \cup \dots \cup L_{r_m}^{(g_{i_m})})\]
for all $r_1 \in [0,e_{g_{i_1}}-1], \dots, r_m \in [0,e_{g_{i_m}}-1]$ where
\[E_k := \begin{cases}
E_{r_j}^{(g_k)}, & \text{if } k = i_j \text{ for some } j \in [1,m] \\
\alpha_k, & \text{otherwise}
\end{cases}\]
for all $k \in [1,n]$. To get a well-defined $\HKP^+(G \wr H)$-instance, we write $E_{r_1,\dots,r_m} = \beta_1 \cdots \beta_s \beta_{s+1}$ with variables in $Y = \{y_1,\dots,y_s\}$ and adjust the loop constraints in $L_{r_1,\dots,r_m}$ accordingly. By construction any solution $\nu$ of $(E,L)$ can be transformed into a solution of $(E_{r_1,\dots,r_m},L_{r_1,\dots,r_m})$ where $r_j := \nu(x_{i_j}) \text{ mod } e_{g_{i_j}}$ for all $j \in [1,m]$. Conversely, any solution of $(E_{r_1,\dots,r_m},L_{r_1,\dots,r_m})$ can be transformed into a solution of $(E,L)$. Moreover, note that $(E_{r_1,\dots,r_m},L_{r_1,\dots,r_m})$ is clearly torsion-free and in $GH$-form. If we write $u_i := \gamma(\beta_i)$ for all $i \in [1,s+1]$, then by the choice of $e_{g}$ for all $i,j \in P_{E_{r_1,\dots,r_m}}$ with $u_i,u_j \notin H$ and $\sigma(u_i)$ and $\sigma(u_j)$ commensurable it holds that $\sigma(u_i) = \sigma(u_j)$ or $\sigma(u_i) = \sigma(u_j)^{-1}$.

Finally, we construct for every $(E_{r_1,\dots,r_m},L_{r_1,\dots,r_m})$ an $\HKP^+(G \wr H)$-instance that is c-simplified. Let $C = \{c_1,\dots,c_k\}$ be a $||$-class of $E_{r_1,\dots,r_m}$ with $c_1 < \dots < c_k$. Then for all $i \in [2,m]$ with $\sigma(u_{c_1}) = \sigma(u_{c_i})^{-1}$ we replace $u_{c_i}^{y_{c_i}}$ in $E_{r_1,\dots,r_m}$ by an expression of the form
\[\gamma_1 \cdots \gamma_7 := \sigma(u_{c_i})^{z_1} (\tau(u_{c_i}) \sigma(u_{c_1}))^{z_2} \sigma(u_{c_1})^{-1} \sigma(u_{c_1}) \sigma(u_{c_i})^{z_5} \sigma(u_{c_i})^{-1} \sigma(u_{c_i})\]
and add the corresponding loop constraints $\{(0,3), (1,6)\}$ to $L_{r_1,\dots,r_m}$ by adjusting indices properly. Intuitively, for all valuations $\nu$ we force that $\nu(z_2) = \nu(z_1)+1$ and $\nu(z_5) = \nu(z_2)+1$. The idea is to walk with $\sigma(u_{c_i})^{z_1}$ to the end of the progression associated to $u_{c_i}^{y_{c_i}}$ and then place with $(\tau(u_{c_i}) \sigma(u_{c_1}))^{z_2}$ the elements in direction of $\sigma(u_{c_1})$ and walk with $\sigma(u_{c_i})^{z_5}$ back to the end of the progression again. With this method only factors in $H$
do not satisfy the commensurability property.

Note that all constructed expressions are torsion-free and in $GH$-form. Doing this for all $||$-classes concludes the proof.
\end{proof}

\subsection{Proof of Lemma~\ref{lem:stacking-free2}}
If $(i,h)$ is an address of a knapsack expression $E = \alpha_1 \cdots \alpha_n$ with $i \in P_E$ and $\sigma(\gamma(\alpha_i)) \neq 1$ and $\nu \in \N^X$ is a valuation, then
\[(\sigma(\nu(\alpha_1 \cdots \alpha_{i-1})) h (h^{-1} \sigma(\gamma(\alpha_i)) h)^j)_{0 \leq j \leq \nu(x_i)-1}\]
is the ray associated to $(i,h)$.

For a non-empty set $S := \{E_1,\dots,E_m\}$ of exponent expressions over $G$ with variables in $X$ we define the set of solutions by $\sol_G(S) := \bigcap_{i = 1}^m \sol_G(E_i)$.
Since by assumption $\ExpEq(G)$ is decidable, we have that $\sol_G(S)$ is decidable as well.

Let $(E = \alpha_1 \cdots \alpha_n \alpha_{n+1},L,D)$ be a normalized $\HKP^\pm(G \wr H)$-instance
where $\alpha_i = g_i$ or $\alpha_i = g_i^{x_i}$ for all $i \in [1,n]$
and $\alpha_{n+1} \in G \wr H$.
Let $I \subseteq [1,n+1]$ be the set of stacking indices.
We say that an address $(i,h)$ is {\em stacking} if $i$ is stacking.
Let $C \subseteq A_E$ be a set which contains at least one stacking address.
We will construct a normalized $\HKP^\pm(G \wr H)$-instance $(E_C,L_C,D_C)$
and a set of exponent expressions $S_C$ over the variable set $\{x_i \mid i \in I \}$.
Intuitively, $C$ represents an intersection point of rays with the progression of at least one stacking index.
In $(E_C,L_C,D_C)$ the intersection point is skipped and $S_C$ expresses that the elements at this point multiply to 1.

We will prove that $(E,L,D)$ has a solution
if and only if there exists a set $C \subseteq A_E$ which intersects $I \times H$ and a valuation $\nu_C$
which satisfies both $(E_C,L_C,D_C)$ and $S_C$.
Furthermore, we show that the number of addresses $(i,h) \in I \times H$
decreases from $E$ to $E_C$.
Hence, by iterating this procedure we end up with stacking-free $\HKP^\pm(G \wr H)$-instances.

\subparagraph{Construction}
Let $S := \emptyset$ and let $C := \{(i_1,h_1),\dots,(i_m,h_m)\} \subseteq A_E$ with $i_1 < \dots < i_m$ be a set of addresses
which intersects $I \times H$.
We first add for all $i \in I$ with $\supp(g_i) = \{s_1, \dots, s_{m_i}\}$ the expression
\[s_1 s_1^{-1} \cdots s_{m_i} s_{m_i}^{-1} = \gamma_1 \cdots \gamma_{2m_i}\]
before $\alpha_i$ needed later to define loop and disjointness constraints. By adjusting indices we can assume that $E$ is in that form.
For all $i \in I$ we define the function $\rho_i \colon \supp(g_i) \to [1,n+1]$
such that $\rho_i(s)$ is the index of the added $s$ before $\alpha_i$ for any $s \in  \supp(g_i)$.
We construct the knapsack expression $E_C$ from $E$ by replacing for all $j \in [1,m]$ the atom
\[
	\alpha_{i_j} \text{ by }
\begin{cases}
g_{i_j}^{x_{i_{j,1}}} \sigma(g_{i_j}) g_{i_j}^{x_{i_{j,3}}} =: \beta_{i_{j,1}} \beta_{i_{j,2}} \beta_{i_{j,3}}, & \text{if } i_j \notin I \\
(f^\prime,1)^{x_{i,1}} =: \beta_{i_{j,1}}, & \text{if } {i_j} \in I \setminus \{n+1\} \text{ and } g_{i_j} = (f,1), \\
(f^\prime,h) =: \beta_{i_{j,1}}, & \text{if ${i_j} = n+1$ is stacking and } g_{i_j} = (f,h), \\
\end{cases}\]
where we define
\[
	f^\prime(h') := 
	\begin{cases}
	1, & \text{if } h' = h_{i_j} \\
	f(h'), & \text{otherwise}
	\end{cases}
\]
for all $h' \in H$.
We remark that we can easily compute a representation of the elements $(f^\prime,1)$
and $(f^\prime,h)$ as words over the generators of $G$ and $H$ from the representation of $g_i$.
By making indices continuous and adjusting $i_{j,1}, i_{j,2}, i_{j,3}$ and $\rho_i$ accordingly, we can write
\[
	E_C = \beta_1 \cdots \beta_r \beta_{r+1}
\]
with variables in $Y := \{y_1,\dots,y_r\}$ and $u_i := \gamma(\beta_i)$ for all $i \in [1,r+1]$. 

For $j \in [1,m]$ with $i_j$ non-stacking in $E$ we set $\rho_{i_{j,1}}(1) := i_{j,1}$. Since $E$ is in $GH$-form, we have that $\supp(g_{i_j}) = \{1\}$ if $i_j$ is non-stacking. 

We define the loop constraints
\[L_C := L \cup \{(\rho_{i_{j,1}}(h_j), \rho_{i_{j+1,1}}(h_{j+1})) \mid j \in [1,m-1]\}\]
where we adjust the indices in $L$ properly. Intuitively, the loop constraints ensure that every solution makes every $(i_j,h_j)$ reach the intersection point given by $C$. But after the replacement of $g_{i_j}$ these addresses do not put an element at the intersection point anymore. 

Let $\operatorname{id} \colon A_E \to [1,n]$ be the map defined by $\operatorname{id}((i,h)) := i$ for all $(i,h) \in A_E$.
Let $\alpha^\prime \colon [1,n+1] \setminus \operatorname{id}(C) \to [1,r+1]$ be the map defined by the adjustment of the indices. Then we define $\alpha \colon [1,n+1] \to [1,r+1]$ such that
\[\alpha(k) := \begin{cases}
i_{j,1}, & \text{if } k = i_j \text{ for some } j \in [1,m] \\
\alpha^\prime(k), & \text{otherwise}
\end{cases}\]
for all $k \in [1,n+1]$. To ensure that every address not in $C$ does not reach the point given by the address $(i_j,h_j) \in C \cap (I \times H)$, we define the disjointness constraints
\begin{equation*}
\begin{split}
D_C := & D^\prime \cup \{(\rho_{\alpha(k)}(h)+1, \rho_{i_{j,1}}(h_j)+1) \mid (k,h) \in A_E \setminus C \text{ and } k \in I \} \cup \\
& \{(\alpha(k), \rho_{i_{j,1}}(h_j))+1 \mid (k,h) \in A_E \setminus C \text{ and } k \notin I \}
\end{split}
\end{equation*}
where we adjust the indices in $D$ as follows:
\begin{equation*}
\begin{split}
D^\prime := & \{(i_{j,x},\alpha(\ell)) \mid (k,\ell) \in D \wedge k = i_j \text{ for some } j \in [1,m] \text{ with } {i_j} \notin I \wedge x \in [1,3]\} \cup \\
& \{(i_{j,x},i_{j^\prime,y}) \mid (i_j,i_{j^\prime}) \in D \text{ for some } j,j^\prime \in [1,m] \text{ with } i_j,i_{j^\prime} \notin I \wedge x,y \in [1,3]\} \cup \\
& \{(\alpha(k),\alpha(\ell)) \mid (k,\ell) \in D\}
\end{split}
\end{equation*}
and assume without loss of generality that if $(k,\ell) \in D$ and either $k = i_j$ or $\ell = i_j$ for some $j \in [1,m]$ with ${i_j} \notin I$, then we always have that $k = i_j$. Intuitively, if $D$ contains a disjointness constraint for a ray that is split into parts, then $D^\prime$ ensures this constraint for every such part. Note that $(E_C,L_C,D_C)$ is still normalized. 

We now extend the set of exponent expressions. Let
\[a_j := \begin{cases}
\tau(g_{i_j})(h_j)^{y_{i_{j,1}}}, & \text{if } \sigma(g_{i_j}) = 1 \text{ and } i_j \neq n+1 \\
\tau(g_{i_j})(h_j), & \text{otherwise}
\end{cases}\]
for all $j \in [1,m]$. Then we define
\[S_C := S \cup \Bigg\{ \prod_{j=1}^m a_j \Bigg\}\]
where we replace variables $x_i$ in $S$ by $y_{\alpha(i)}$. The additional exponent expression ensures that the elements written at the point given by $C$ multiply to 1. Here we only need variables for stacking indices since elements of non-stacking indices can visit the point at most once.

We repeat this process for $(E_C,L_C,D_C)$ and $S_C$ until the resulting $\HKP^\pm(G \wr H)$-instance $(E^\prime,L^\prime,D^\prime)$ has no stacking addresses left. If for the corresponding set of exponent expressions $S^\prime$ it holds that $\sol_G(S^\prime) \neq \emptyset$, we construct a stacking-free $\HKP^\pm(G \wr H)$-instance by removing the exponents of powers with base 1. Let $(E_1,L_1,D_1),\dots,(E_t,L_t,D_t)$ be the constructed stacking-free, normalized $\HKP^\pm(G \wr H)$-instances for all possible choices of sets $C$ during the construction. We claim that $\sol_{G \wr H}(E,L,D) \neq \emptyset$ if and only if $\bigcup_{i=1}^t \sol_{G \wr H}(E_i,L_i,D_i) \neq \emptyset$.

\subparagraph{Termination}
We show that in each step of the construction above the number of stacking addresses gets strictly smaller. This means that after a finite number of steps the resulting $\HKP^\pm(G \wr H)$-instance has no stacking addresses left and the construction terminates. For a knapsack expression $E = \alpha_1 \cdots \alpha_n \alpha_{n+1}$ with $g_i := \gamma(\alpha_i)$ for all $i \in [1,n+1]$ let 
\[s(E) := |\{(i,h) \in A_E \mid i \in I \}|\]
be the number of stacking addresses. Let $C \subseteq A_E$ be a set of addresses that contains a stacking address $(i,h)$. During the construction of $E_C$ we replace $g_i$ in $E$ by $(f^\prime,\sigma(g_i))$ where $f^\prime(h) = 1$. This means that $\supp(f^\prime) = \supp(f)-1$. Thus, it holds that
\[s(E_C) = s(E) - |\{(i,h) \in C \mid i \in I \}|.\]
Since $C$ contains at least one stacking address, it follows that $s(E_C) < s(E)$.

\subparagraph{Correctness}
It remains to show that $(E,L,D)$ has a solution if and only if one of the constructed $(E_1,L_1,D_1),\dots,(E_t,L_t,D_t)$ has a solution. We consider each step of the construction separately. Let $(E = \alpha_1 \cdots \alpha_n \alpha_{n+1},L,D)$ be a normalized $\HKP^\pm(G \wr H)$-instance with $s(E) \geq 1$ and $g_i := \gamma(\alpha_i)$ for all $i \in [1,n+1]$. Let $S$ be a set of exponent expressions over $G$ with variables in $X$. We assume that $(E,L,D)$ and $S$ are generated during the construction. We show that $\sol_{G \wr H}(E,L,D) \cap \sol_G(S) \neq \emptyset$ if and only if there exists $C \subseteq A_E$ containing a stacking address such that $\sol_{G \wr H}(E_C,L_C,D_C) \cap \sol_G(S_C) \neq \emptyset$.

For the first direction assume that $\nu \in \sol_{G \wr H}(E,L,D) \cap \sol_G(S) \neq \emptyset$. As $s(E) \geq 1$, there is an address $(i,h) \in A_E$ of a sacking element $g_i$. Let $h_C := \supp_E^\nu(\rho_i(h)+1)$ be the point visited by $(i,h)$ under $\nu$ and
\begin{align*}
C := & \{(j,1) \in A_E \mid j \notin S \text{ and } h_C \in \supp_E^\nu(j)\} \cup \\
& \{(j,h^\prime) \in A_E \mid j \in S \text{ and } h_C \in \supp_E^\nu(\rho_j(h^\prime)+1)\}
\end{align*}
be the set of all addresses reaching $h_C$ under $\nu$. We write $C = \{(i_1,h_1),\dots, (i_m,h_m)\}$ with $i_1 < \dots < i_m$. Let $(E_C,L_C,D_C)$ be the $\HKP^\pm(G \wr H)$-instance and $S_C$ be the set of exponent expressions constructed above from $(E,L,D)$ and $S$ with respect to $C$. We now define a valuation $\nu_C \in \N^Y$. For all $j \in [1,m]$ with ${i_j} \notin I$ we assign $\nu_C(y_{i_{j,1}}) := e$ and $\nu_C(y_{i_{j,3}}) := \nu(x_{i_j})-e-1$ where $e \in [0,\nu(x_{i_j})-1]$ such that $\sigma(\nu(\alpha_1) \cdots \nu(\alpha_{i_j-1}) g_{i_j}^e) = h_C$. For all $k \in P_E$ with $k \in I$ or $k \notin \operatorname{id}(C)$ we assign $\nu_C(y_{\alpha(k)}) := \nu(x_k)$.

Since $\nu \in \sol_{G \wr H}(E,L,D)$ and the construction only splits up some rays, we have that $\sigma(\nu_C(E_C)) = 1$ and $\nu_C$ fulfills all loop constraints in $L_C$ and all disjointness constraints in $D_C$ by definition of $C$. As $\tau(\nu(E))(h^\prime) = 1$ for all $h^\prime \in H$ and there is no address of $E_C$ visiting the point $h_C$ under $\nu_C$, it holds that $\tau(\nu_C(E_C))(h^\prime) = 1$ for all $h^\prime \in H$. Moreover, from $\prod_{j=1}^m a_j^\nu = 1$ with
\[a_j^\nu := \begin{cases}
\tau(g_{i_j})(h_j)^{\nu(x_{i_j})}, & \text{if } \sigma(g_{i_j}) = 1 \text{ and } i_j \neq n+1 \\
\tau(g_{i_j})(h_j), & \text{otherwise}
\end{cases}\]
and from $\nu \in \sol_G(S)$ it follows that $\nu_C \in \sol_G(S_C)$ since the exponent expressions in $S$ only contain variables with stacking indices. Thus, it holds that $\nu_C \in \sol_{G \wr H}(E_C,L_C,D_C) \cap \sol_G(S_C)$.

For the other direction assume that $\nu_C \in \sol_{G \wr H}(E_C,L_C,D_C) \cap \sol_G(S_C)$ for some set of addresses $C = \{(i_1,h_1),\dots,(i_m,h_m)\} \subseteq A_E$ with $i_1 < \dots < i_m$ containing a stacking address $(i,h)$. We now define a valuation $\nu \in \N^X$. For all $j \in [1,m]$ with ${i_j} \notin I$ we assign $\nu(x_{i_j}) := \nu_C(y_{i_{j,1}}) + \nu_C(y_{i_{j,3}}) + 1$. For all $k \in P_E$ with $k \in I$ or $k \notin \operatorname{id}(C)$ we assign $\nu(x_k) := \nu_C(y_{\alpha(k)})$.

Since by construction $S_C$ only contains variables with stacking indices and $\nu_C \in \sol_G(S_C)$, we have that $\nu \in \sol_G(S)$ as $S_C$ extends $S$ by one expression. The additional exponent expression in $S_C$ ensures that by definition of $\nu$ it holds that $\prod_{j=1}^m a_j^\nu = 1$. Therefore, the disjointness constraints in $D_C$ imply that we have $\tau(\nu(E))(h_C) = 1$ where $h_C := \supp_E^\nu(\rho_i(h)+1)$. By construction of $(E_C,L_C,D_C)$ it follows that $\tau(\nu(E))(h^\prime) = 1$ for all $h^\prime \in H$. Moreover, since $\sigma(\nu_C(E_C)) = 1$, it holds that $\sigma(\nu(E)) = 1$ and the definitions of $L_C$ and $D_C$ imply that $\nu$ fulfills all loop constraints in $L$ and all disjointness constraints in $D$. Thus, we have that $\nu \in \sol_{G \wr H}(E,L,D) \cap \sol_G(S)$.

This implies that for a normalized $\HKP^\pm(G \wr H)$-instance $(E,L,D)$ and $S = \emptyset$ it holds that $\sol_{G \wr H}(E,L,D) = \sol_{G \wr H}(E,L,D) \cap \sol_G(S) \neq \emptyset$ if and only if there exist sets of addresses $C_1,\dots,C_m$ such that for the normalized $\HKP^\pm(G \wr H)$-instance $(E^\prime,L^\prime,D^\prime)$ and the set of exponent expressions $S^\prime$ constructed with respect to $C_1,\dots,C_m$ we have $s(E^\prime) = 0$ and $\sol_{G \wr H}(E^\prime,L^\prime,D^\prime) \cap \sol_G(S^\prime) \neq \emptyset$. Since $s(E^\prime) = 0$ and $S^\prime$ only contains variables with stacking indices, it holds that $\sol_{G \wr H}(E^\prime,L^\prime,D^\prime) \cap \sol_G(S^\prime) \neq \emptyset$ if and only if $\sol_{G \wr H}(E^\prime,L^\prime,D^\prime) \neq \emptyset$ and $\sol_G(S^\prime) \neq \emptyset$. Thus, it suffices to construct $\HKP^\pm(G \wr H)$-instances where the corresponding set of exponent expressions has a solution. Removing the exponents of powers in $E^\prime$ that have base 1 yields a stacking-free, normalized $\HKP^\pm(G \wr H)$-instance that fulfills the claim.

\subsection{Proof of Lemma~\ref{lem:interval-KP2}}
If $E = \alpha_1 \cdots \alpha_n$ is a stacking-free knapsack expression in $GH$-form, for all addresses $(i,h) \in A_E$ it holds that $h = 1$. Thus, we can write $A_E = \{i \in P_E \mid \gamma(\alpha_i) \in GH \setminus H\}$.
In the following we often view an addresses $i \in A_E$ as the associated ray 
\[(\sigma(\nu(\alpha_1 \cdots \alpha_{i-1})) \sigma(\gamma(\alpha_i))^j)_{0 \leq j \leq \nu(x_i)-1}\]
under a valuation $\nu \in \N^X$. We say that two rays are {\em parallel} if their periods are commensurable.

We first construct for every splitting of rays into subrays and equivalence relation on these subrays an $\MKP^\pm(H)$-instance. We then show that the resulting instances fulfill the claim. Note that by \Cref{lem:MKP} the resulting $\MKP^\pm(H)$-instances can be transformed to $\KP^\pm(H)$-instances that prove the lemma.

\subparagraph{Construction}
Let $(E = \alpha_1 \cdots \alpha_n \alpha_{n+1}, L, D)$ be a stacking-free, normalized $\HKP^\pm(G \wr H)$-instance with $g_i := \gamma(\alpha_i)$ for all $i \in [1,n+1]$. Let $A_E = \{a_1,\dots,a_m\}$ be the rays of $E$ with $a_1 < \dots < a_m$. Note that if we split a ray at the intersection points with other rays, then every intersection point results in at most two new subrays. As there are $m-1$ other rays, a ray is split into at most $1+2(m-1) = 2m-1$ subrays. Let $N := [1,2m-1]^m$ and for every $\eta := (n_{a_1},\dots,n_{a_m}) \in N$ we define the knapsack expression $E_\eta^\prime$ by replacing $g_{a_i}^{x_{a_i}}$ in $E$ by 
\[g_{a_i}^{y_1} \sigma(g_{a_i})^{-1} \sigma(g_{a_i}) \cdots g_{a_i}^{y_{3n_{a_i}-2}} \sigma(g_{a_i})^{-1} \sigma(g_{a_i}) = \beta_1 \cdots \beta_{3n_{a_i}}.\]
This means we split $g_{a_i}^{x_{a_i}}$ into $n_{a_i}$ parts where the factors $ \sigma(g_{a_i})^{-1} \sigma(g_{a_i})$ are needed later to define loop constraints.
By making indices continuous, we can write
\[E_\eta^\prime = \beta_1 \cdots \beta_r \beta_{r+1}\]
with variables in $Y := \{y_1,\dots,y_r\}$ and $u_i := \gamma(\beta_i)$ for all $i \in [1,r+1]$. We remark that if $E$ is stacking-free and normalized, then so is $E_\eta^\prime$. For every $i \in A_E$ and $j \in [1,n_i]$ let $a_{i,j}$ be the index of the $j$-th subray of $i$ in $E_\eta^\prime$. Furthermore, let $\alpha \colon [1,n+1] \setminus A_E \to [1,r+1]$ be defined by the adjustment of the indices.

Let $\Theta_\eta$ be the set of all equivalence relations on $A_{E_\eta^\prime}$. Note that $\Theta_\eta$ is finite and can be computed by dividing the rays of $E_\eta^\prime$ into equivalence classes. Then for all $\sim \in \Theta_\eta$ we define the loop constraints
\begin{align}
L_\sim^\prime := & L \cup \label{set1} \\
& \{(i-1,j-1) \mid i,j \in A_{E_\eta^\prime} \wedge i < j \wedge i \sim j\} \cup \label{set2} \\
& \{(i+1,j+1) \mid i,j \in A_{E_\eta^\prime} \wedge i < j \wedge i \sim j\} \label{set3}
\end{align}
where we adjust the indices in $L$ properly. For two rays $i$ and $j$ of $E_\eta^\prime$ with $i \sim j$ the loop constraint $(i,j-1)$ in \labelcref{set2} ensures that the starting points of $i$ and $j$ are equal. The loop constraint $(i+2,j+1)$ in \labelcref{set3} ensures that any solution $\nu$ satisfies
\[\sigma(\nu(\beta_1) \cdots \nu(\beta_i)) \sigma(u_i)^{-1} = \sigma(\nu(\beta_1) \cdots \nu(\beta_j)) \sigma(u_j)^{-1}\] 
which means that the endpoints of $i$ and $j$ are equal. Since $E_\eta^\prime$ is normalized, this implies that the rays $i$ and $j$ must be equal. To ensure that the rays in different $\sim$-classes are disjoint, we define the disjointness constraints
\[D_\sim^\prime := D^\prime \cup \{(i,j) \mid i,j \in A_{E_\eta^\prime} \wedge i < j \wedge i \nsim j\}\]
where we adjust $D$ as follows:
\begin{equation*}
\begin{split}
D^\prime := & \{(a_{i,j},a_{k,\ell}) \mid i,k \in A_E \wedge (i,k) \in D \wedge j \in [1,n_i] \wedge \ell \in [1,n_j]\} \cup \\
& \{(a_{i,j},\alpha(k)) \mid i \in A_E \wedge k \in [1,d+1] \setminus A_E \wedge (i,k) \in D \wedge j \in [1,n_i]\} \cup \\
& \{(\alpha(i),\alpha(k)) \mid i,k \in [1,d+1] \setminus A_E \wedge (i,k) \in D\}
\end{split}
\end{equation*}
and assume without loss of generality that if $(i,k) \in D$ with $i \in A_E$ or $k \in A_E$, then we always have that $i \in A_E$. Intuitively, if $D$ contains a disjointness constraint for a ray that is split into parts, then $D^\prime$ ensures this constraint for every such part.

Now by construction it is enough to evaluate $\tau(E_\eta^\prime)$ only within $\sim$-classes. This means that if there is a $\sim$-class $C = \{c_1,\dots,c_k\}$ with $c_1 < \dots < c_k$ such that $\prod_{i=1}^k \tau(u_{c_i})(1) \neq 1$, then we demand that every solution sets $y_{c_i}$ to 0 for all $i \in [1,k]$. To this end, for all such $\sim$-classes we remove $u_{c_i}^{y_{c_i}}$ from $E_\eta^\prime$ for all $i \in [1,k]$ and adjust $L_\sim^\prime$ and $D_\sim^\prime$ properly. Let $(E_\eta = \gamma_1 \cdots \gamma_s \gamma_{s+1},L_\sim,D_\sim)$ be the resulting $\HKP^\pm(G \wr H)$-instance with $v_i := \gamma(\gamma_i)$ for all $i \in [1,s+1]$. Then $(\sigma(E_\eta),L_\sim,D_\sim)$ is an $\MKP^\pm(H)$-instance where we let $\sigma(g^x) := \sigma(g)^x$ for an atom $g^x$. We claim that 
\[\sol_{G \wr H}(E,L,D) \neq \emptyset \text{ if and only if } \bigcup_{\eta \in N \wedge \sim \in \Theta_\eta} \sol_H(\sigma(E_\eta),L_\sim,D_\sim) \neq \emptyset.\]

\subparagraph{Correctness}
It remains to show that the $\HKP^\pm(G \wr H)$-instance $(E,L,D)$ has a solution if and only if there exist $\eta \in N$ and $\sim \in \Theta_\eta$ such that the $\MKP^\pm(H)$-instance $(\sigma(E_\eta),L_\sim,D_\sim)$ has a solution. For the first direction we assume that $\nu \in \sol_{G \wr H}(E,L,D)$. For all $i \in P_E$ let $\sigma_i \colon [0,\nu(x_i)-1] \to H$ such that
\[\sigma_i(e) := \sigma(\nu(\alpha_1 \cdots \alpha_{i-1}) g_i^e)\]
for all $e \in [0,\nu(x_i)-1]$. Furthermore, we define a function $f \colon H \to \mathcal{P}(A_E)$ such that
\[f(h) := \{i \in A_E \mid \exists e \in [0,\nu(x_i)-1] \colon \sigma_i(e) = h\}\]
for all $h \in H$. This means that $f$ maps a point $h \in H$ to the set of rays that visit $h$ under the valuation $\nu$.

We now split the rays into subrays to get $\eta \in N$ and $\sim \in \Theta_\eta$ such that $(\sigma(E_\eta),L_\sim,D_\sim)$ has a solution. For every $i \in A_E$ with $\nu(x_i) \neq 0$ there is a partition of $[0,\nu(x_i)-1]$ into disjoint intervals $[s_1^{(i)},e_1^{(i)}], \dots, [s_{n_i}^{(i)},e_{n_i}^{(i)}]$ such that for all $j \in [1,n_i]$ and $k \in [s_j^{(i)},e_j^{(i)}]$ it holds that
\[f(\sigma_i(s_j^{(i)})) = f(\sigma_i(k))\]
and for all $j \in [1,n_i-1]$ it holds that
\[f(\sigma_i(s_j^{(i)})) \neq f(\sigma_i(s_{j+1}^{(i)})).\]
For $i \in A_E$ with $\nu(x_i) = 0$ we set $n_i := 1$ and $[s_1^{(i)},e_1^{(i)}] := [1,0] = \emptyset$. Intuitively, we split a ray whenever the intersection with another ray starts or ends.

We need to show that $n_i \leq 2m-1$ for all $i \in A_E$. Let $i \in A_E$ be a ray and $i \neq j \in A_E$ be one of the $m-1$ other rays such that $i$ and $j$ intersect. Let $n_{i,j}$ be the number of disjoint intervals in the partition of $[0,\nu(x_i)-1]$ as defined above with respect to the function defined by
\[f_j(h) := \{k \in A_E \setminus \{j\} \mid \exists e \in [0,\nu(x_k)-1] \colon \sigma_k(e) = h\}\]
for all $h \in H$. That is, we do not split $i$ at the intersection with $j$. If $i$ and $j$ are non-parallel, then they intersect in exactly one point $\sigma_i(z)$ for some $z \in [0,\nu(x_i)-1]$. This implies that $j \in f(\sigma_i(z))$ and $j \notin f(\sigma_i(e))$ for all $e \in [0,\nu(x_i)-1] \setminus \{z\}$. Thus, we have $n_i \leq n_{i,j} +2$. If $i$ and $j$ are parallel, then they have the same period since $E$ is c-simplified. So the intersection of $i$ and $j$ is a subray of $i$ with starting point $\sigma_i(z_1)$ and endpoint $\sigma_i(z_2)$ for some $0 \leq z_1 \leq z_2 \leq \nu(x_i)-1$. This implies that $j \in f(\sigma_i(e))$ for all $e \in [z_1,z_2]$ and $j \notin f(\sigma_i(e))$ for all  $e \in [0,\nu(x_i)-1] \setminus [z_1,z_2]$. Thus, we have $n_i \leq n_{i,j} +2$. By induction it follows that $n_i \leq 1+2(m-1)$. Therefore, $\eta := (n_{a_1},\dots,n_{a_m}) \in N$ where we recall that $A_E = \{a_1,\dots,a_m\}$ with $a_1 < \dots < a_m$.

Let $E_\eta^\prime$ be the knapsack expression corresponding to $\eta$ as constructed above and $a_{i,j}$ for $i \in A_E$ and $j \in [1,n_i]$ be the index of the $j$-th subray of $i$ in $E_\eta^\prime$. Then we define the equivalence relation $\sim \in \Theta_\eta$ such that for all $i,k \in A_E, j \in [1,n_i]$ and $\ell \in [1,n_k]$ it holds that $a_{i,j} \sim a_{k,\ell}$ if and only if
\[|[s_j^{(i)},e_j^{(i)}]| = e_j^{(i)} - s_j^{(i)} +1 = e_\ell^{(k)} - s_\ell^{(k)} +1 = |[s_\ell^{(k)},e_\ell^{(k)}]|\]
and
\[\sigma_i(s_j^{(i)} + z) = \sigma_k(s_\ell^{(k)} + z)\]
for all $z \in [0,e_j^{(i)} - s_j^{(i)}]$. This means that $\sim$ relates all equal subrays.

Now we define the valuation $\nu^\prime \in \N^Y$ such that
\[\nu^\prime(y_k) := \begin{cases}
e_j^{(i)} - s_j^{(i)} + 1, & \text{if } k = a_{i,j} \text{ for some } i \in A_E \text{ and } j \in [1,n_i] \\
\nu(x_k), & \text{otherwise}
\end{cases}\]
for all $k \in P_{E_\eta^\prime}$. Let $\beta \colon [1,s+1] \to [1,r+1]$ map the indices of elements of $E_\eta$ to the corresponding indices of elements of $E_\eta^\prime$ that are not removed. Since $\nu$ is a solution of $(E,L,D)$, for every $i \in A_E$ with $\nu(x_i) \neq 0$ and $j \in [1,n_i]$ it holds that $\prod_{\ell=1}^k \tau(u_{c_\ell})(1) = 1$,  where $C = \{c_1,\dots,c_k\}$ with $c_1 < \dots < c_k$ is the $\sim$-class containing $a_{i,j}$, and therefore $u_{a_{i,j}}^{y_{a_{i,j}}}$ is not removed from $E_\eta^\prime$. By construction of $(E_\eta,L_\sim,D_\sim)$ and since $\nu \in \sol_H(\sigma(E),L,D)$, it follows that for $\nu^{\prime\prime} \in \N^Z$ defined by $\nu^{\prime\prime}(z_i) := \nu^\prime(y_{\beta(i)})$ for all $i \in P_{E_\eta}$ we have $\nu^{\prime\prime} \in \sol_H(\sigma(E_\eta),L_\sim,D_\sim)$.

For the other direction assume that $\nu^\prime \in \sol_H(\sigma(E_\eta),L_\sim,D_\sim)$ for some $\eta \in N$ and $\sim \in \Theta_\eta$. Since after the construction $g_i^{x_i}$ is split into $g_i^{z_{i_1}} h_1 g_i^{z_{i_2}} \cdots h_{m_i-1} g_i^{z_{i_{m_i}}}$ for some $m_i \in [0,n_i]$ and products $h_j$ of elements of $H$ such that $h_j = 1$ for all $j \in [1,m_i-1]$, we define the valuation $\nu \in \N^X$ by
\[\nu(x_i) := \nu^\prime(z_{i_1}) + \dots + \nu^\prime(z_{i_{m_i}})\]
for all $i \in P_E$. Since $\sigma(\nu^\prime(E_\eta)) = 1$, we also have that $\sigma(\nu(E)) = 1$. Moreover, as $\nu^\prime \in \sol_H(\sigma(E_\eta),L_\sim,D_\sim)$, the definitions of $L_\sim$ and $D_\sim$ imply that all loop constraints in $L$ and all disjointness constraints in $D$ are fulfilled under $\nu$.

It remains to show that $\tau(\nu(E))(h) = 1$ for all $h \in H$. Note that we can regard $\sim$ also as equivalence relation on $A_{E_\eta}$. For $i,j \in A_{E_\eta}$ we say that $i \sim j$ if and only if $\beta(i) \sim \beta(j)$. Moreover, for $i \in A_{E_\eta}$ let $\supp_{E_\eta}^{\nu^\prime}(i)$ be the support of the ray $i$ under $\nu^\prime$. Since $\nu^\prime$ fulfills the loop constraints in $L_\sim$, for any two rays $i, j \in A_{E_\eta}$ in the same $\sim$-class $C$ it holds that $\supp_{E_\eta}^{\nu^\prime}(i) = \supp_{E_\eta}^{\nu^\prime}(j)$ and we define $\supp_{E_\eta}^{\nu^\prime}(C) := \supp_{E_\eta}^{\nu^\prime}(i)$. By construction of $E_\eta$ we have that $\prod_{i=1}^k \tau(v_{c_i})(1) = 1$ for any $\sim$-class $C = \{c_1,\dots,c_k\}$ with $c_1 < \dots < c_k$. As the disjointness constraints in $D_\sim$ ensure that rays of different $\sim$-classes are disjoint, for any $\sim$-class $C$ it follows that $\tau(\nu^\prime(E_\eta))(h) = 1$ for all $h \in \supp_{E_\eta}^{\nu^\prime}(C)$. Thus, since any $i \in A_{E_\eta}$ is contained in a $\sim$-class, we have that $\tau(\nu^\prime(E_\eta))(h) = 1$ for all $h \in H$. By definition of $\nu$ this implies that $\tau(\nu(E))(h) = 1$ for all $h \in H$ since the rays of $E$ under $\nu$ are built of the rays of $E_\eta$ under $\nu^\prime$.

\subsection{Reduction in abelian case}
As consequence of \Cref{thm:normal} we can start the reduction with a normalized $\HKP^+(G \wr H)$-instance. Again, in the first step we make the instance stacking-free.

\begin{lemma}\label{lem:stacking-free}
For any normalized $\HKP^+(G \wr H)$-instance one can effectively construct an equivalent finite set of stacking-free, normalized $\HKP^+(G \wr H)$-instances.
\end{lemma}
\begin{proof}
We do the same construction as in the proof of \Cref{lem:stacking-free2} but we leave the disjointness constraints out. This results in stacking-free, normalized $\HKP^+(G \wr H)$-instances $(E_1,L_1),\dots,(E_t,L_t)$. For the correctness we assume that the normalized $\HKP^+(G \wr H)$-instance $(E,L)$ with $s(E) \geq 1$ and the set of exponent expressions $S$ are generated during the construction. We need to show that $\sol_{G \wr H}(E,L) \cap \sol_G(S) \neq \emptyset$ if and only if there exists $C \subseteq A_E$ containing an address of a stacking index such that $\sol_{G \wr H}(E_C,L_C) \cap \sol_G(S_C) \neq \emptyset$. The first direction works exactly the same as in the proof of \Cref{lem:stacking-free2}. For the other direction it remains to argue that  $\tau(\nu(E))(h_C) = 1$. But since $G$ is abelian, this follows from the fact that $\prod_{j=1}^m a_j^\nu = 1$ and $\tau(\nu_C(E_C))(h_C) = 1$.
\end{proof}

From now on we assume that $(E,L)$ is a stacking-free, normalized $\HKP^+(G \wr H)$-instance. The next lemma shows how to reduce $\HKP^+(G \wr H)$ for stacking-free, normalized $\HKP^+(G \wr H)$-instances to $\MKP^+(H)$.

\begin{lemma}\label{lem:interval-KP}
For any stacking-free, normalized $\HKP^+(G \wr H)$-instance one can effectively construct an equivalent finite set of $\KP^+(H)$-instances.
\end{lemma}
\begin{proof}
We can again almost copy the proof of \Cref{lem:interval-KP2} by leaving the disjointness constraints out. The result of the construction are $\MKP^+(H)$-instances $(E_1,L_1),\dots,(E_t,L_t)$. For the correctness we have to show that the $\HKP^+(G \wr H)$-instance $(E,L)$ has a solution if and only if there exist $\eta \in N$ and $\sim \in \Theta_\eta$ such that the $\MKP^+(H)$-instance $(\sigma(E_\eta),L_\sim)$ has a solution. The first direction is again the same as in the proof of \Cref{lem:interval-KP2}. For the other direction we only use disjointness constraints to show that $\tau(\nu^\prime(E_\eta))(h) = 1$ for all $h \in \supp_{E_\eta}^{\nu^\prime}(C)$ and $\sim$-classes $C$. But since $G$ is abelian, this also follows from the fact that $\prod_{i \in C} \tau(v_i)(1) = 1$ for any $\sim$-class $C$.
\end{proof}

\section{Proofs from Section~\ref{sec:to-wreath}}

\subsection{Proof of Lemma~\ref{lem:interval-words}}

\begin{lemma}
	\label{lem:pc-words}
	Given $k \in \N$ one can compute $u \in \langle a \rangle^{(\N)}$ with periodic complexity $\ge k$.
\end{lemma}

\begin{proof}
	A function $u \in \langle a \rangle^{(\N)}$ is {\em $(k,s)$-alternating}
	if there are intervals $L_1 = [\ell_1,r_1], \dots, L_k = [\ell_k,r_k]$
	and elements $c_1, \dots, c_k \in \langle a \rangle$
	such that $|L_j| \ge s$, $\ell_1 \le r_1 < \ell_2 \le r_2 < \dots < \ell_k \le r_k$,
	and $u(n) = c_j$ for all $n \in L_j$, $j \in [1,k]$,
	and $c_j \neq c_{j+1}$ and $j \in [1,k-1]$.
	We claim that every $(4k,2^{2^k})$-alternating function $u$ has periodic complexity at least $k$.
	The statement then follows by choosing the word
	$u = (a)^{2^{2^k}} (1)^{2^{2^k}} \dots (a)^{2^{2^k}} (1)^{2^{2^k}}$ consisting of $4k$ blocks.
	Here, the notation $(c)^\ell$ stands for the word consisting of $\ell$ many $c$.
	
	The proof proceeds by induction on $k$.
	Let $L_1, \dots, L_{4k}$ be intervals of size $\ge 2^{2^k}$
	and $c_1, \dots, c_{4k} \in \langle a \rangle$
	such that $u$ is constant $c_j$ on each interval $L_j$ and $c_j \neq c_{j+1}$ for all $j \in [1,4k-1]$.
	Take any basic periodic function $v \neq 1$
	with support $\supp(v) = \{p+qn \mid 0 \le n \le \ell\}$ for some numbers $p,q,\ell$.
	Let $c \in \langle a \rangle$ such that $v(n) = c$ for all $n \in \supp(v)$.
	It suffices to show that $\pc(uv^{-1}) \ge k-1$.
	If the period $q$ is at least $2^{2^{k-1}} + 1$
	then each set $L_j \setminus \supp(v)$ contains an interval of size $2^{2^{k-1}}$.
	Hence $uv^{-1}$ is $(4k,2^{2^{k-1}})$-alternating
	and by induction $\pc(uv^{-1}) \ge k-1$.
	If $q \le 2^{2^{k-1}}$ consider the restriction of $uv^{-1}$ to $D = \{ n \in \N \mid n \equiv p \pmod q \}$.
	Notice that $\supp(v) \subseteq D$ and $\supp(v)$ is convex in $D$,
	i.e. if $n_1 < n_2 < n_3 \in D$ and $n_1,n_3 \in \supp(v)$ then $n_2 \in \supp(v)$.
	Moreover $|L_j \cap D| \ge 2^{2^{k-1}}$ for all $j \in [1,4k]$
	since $|L_j| \ge 2^{2^k}$ and $q \le 2^{2^{k-1}}$.
	Let $J_+ = \{ j \in [1,4k] \mid L_j \cap \supp(v) = L_j \cap D \}$
	and $J_- = \{ j \in [1,4k] \mid L_j \cap \supp(v) = \emptyset \}$,
	which are disjoint sets because $L_j \cap D$ is always nonempty.
	Define $c_j'$ for all $j \in J_+ \cup J_-$ by
	\[
		c_j' = \begin{cases}
			c_j, & \text{if } j \in J_-, \\
			c_jc^{-1}, & \text{if } j \in J_+.
		\end{cases}
	\]
	Notice that $uv^{-1}$ is constant $c_j'$ on each set $L_j \cap D$ for all $j \in J_+ \cup J_-$.
	Morever, $J_+$ is an interval by convexity of $\supp(v)$ in $D$.
	Furthermore, if $j \notin J_+ \cup J_-$
	then $j$ must be adjacent to the interval $J_+$;
	otherwise there would be indices $j_1 < j_2 < j_3$ such that
	$L_{j_1}$ and $L_{j_3}$ both intersect $\supp(v)$, and $L_{j_2}$
	contains a point in $D \setminus \supp(v)$, which again would contradict the convexity of $\supp(v)$ in $D$.
	Therefore $(c_j')_{j \in J_+ \cup J_-}$ is alternating except in at most two positions.
	We can pick a subset $J \subseteq J_+ \cup J_-$ of size $\ge 4k-4$
	such that the sequence $(c_j')_{j \in J}$ is alternating.
	Hence, the periodic subsequence of $uv^{-1}$ induced by $D$
	is $(4(k-1),2^{2^{k-1}})$-alternating.
	By induction we obtain $\pc(uv^{-1}) \ge k-1$,
	concluding the proof.
\end{proof}

First we prove the case $n = 1$.
Let $v = a_1 \dots a_m$ be any function with $\pc(v) \ge k$ (\cref{lem:pc-words}).
Then let $u = a_1 (1)^{m-1} a_2 (1)^{m-1} \dots a_m (1)^{m-1} a_1 \dots a_m$.
Let $p \neq q \in \Z_\infty$.
If $p = \infty$ (or $q = \infty$) then $\loi{p}{u}{}{} \loi{q}{u}{-1}{}$ is $\loi{q}{u}{-1}{}$
(or $\loi{p}{u}{}{}$, respectively)
and has periodic complexity $\ge k$.
Next we can assume that $p,q \in \Z$ and $p < q$ since
$\pc(\loi{p}{u}{}{} \loi{q}{u}{-1}{}) = \pc(\loi{q}{u}{}{} \loi{p}{u}{-1}{})$
(because $\loi{p}{u}{}{} \loi{q}{u}{-1}{}$ is the point-wise inverse of $\loi{q}{u}{}{} \loi{p}{u}{-1}{}$).
If $q-p < m$ then $v$ is a periodic subsequence of $\loi{p}{u}{}{} \loi{q}{u}{-1}{}$.
If $q-p \ge m$ then $v^{-1}$ is a periodic subsequence of $\loi{p}{u}{}{} \loi{q}{u}{-1}{}$.
In any case $\pc(\loi{p}{u}{}{} \loi{q}{u}{-1}{}) \ge \pc(v) \ge k$.

Now let $n \in \N$ be arbitrary and let $u = a_1 \dots a_m$ be any function such that
$\pc(\loi{p}{u}{}{} \loi{q}{u}{-1}{}) \ge k+4(n-1)$ for all $p \neq q \in \Z_\infty$,
which can be constructed as described above.
Then, we set $u_1 = u$ and for all $i \in [2,n]$ we define
\[
	u_i = a_1 1^{|u_{i-1}|} a_2 1^{|u_{i-1}|} \dots a_m.
\]
We claim that for any $p\ne q$ and $i\in[1,n]$, there is a
progression $D\subseteq\Z$ with period $|u_{i-1}|+1$ (if $i = 1$ the period is 1) such that
$\pi_D(\loi{p}{u}{}{i}\loi{q}{u}{-1}{i})$ has
  periodic complexity $\ge k+4(n-1)$.
In particular, we have
$|\supp(\loi{r}{u}{}{j})\cap D|\le 1$ for every $r\in\Z$ and
$j \neq i$.

The claim is obvious for $i=1$. For $i>1$, we distinguish two
cases. First, suppose $p-q$ is divisible by $|u_{i-1}|+1$.  Then the
support of $\loi{p}{u}{}{i}\loi{q}{u}{-1}{i}$ is included in some
progression $D$ with period $|u_{i-1}|+1$. Moreover, we have
$\pi_D(\loi{p}{u}{}{i}\loi{q}{u}{-1}{i})=\loi{p'}{u}{}{}\loi{q'}{u}{-1}{}$
for some $p',q'$ with $p'\ne q'$, hence
$\pc(\pi_D(\loi{p}{u}{}{i}\loi{q}{u}{-1}{i}))=\pc(\loi{p'}{u}{}{}\loi{q'}{u}{-1}{})\ge
k+4(n-1)$. 

Now suppose $p-q$ is not divisible by $|u_{i-1}|+1$. Then there is a progression
$D$ with period $|u_{i-1}|+1$ such that $\pi_D(\loi{p}{u}{}{i}\loi{q}{u}{-1}{i})=u$,
hence $\pc(\pi_D(\loi{p}{u}{}{i}\loi{q}{u}{-1}{i}))=\pc(u)\ge k+4(n-1)$.

Now take numbers $p_1,q_1, \dots, p_n,q_n \in \Z_\infty$
with $p_j \neq q_j$ for some $j \in [1,n]$ and consider
\[
	w = \prod_{i=1}^n \loi{p_i}{u}{}{i}\loi{q_i}{u}{-1}{i}.
\]
We can rewrite the equation to
\begin{equation}
	\label{eq:u}
	\loi{p_j}{u}{}{j}\loi{q_j}{u}{-1}{j} = \Big(\prod_{i < j}^n \loi{p_i}{u}{}{i}\loi{q_i}{u}{-1}{i} \Big)^{-1} \cdot w \cdot \Big(\prod_{i > j}^n \loi{p_i}{u}{}{i}\loi{q_i}{u}{-1}{i} \Big)^{-1}.
\end{equation}
By \eqref{eq:u} and the observation above there exists a progression $D \subseteq \Z$ such that
$\pc(\pi_D(\loi{p_j}{u}{}{j}\loi{q_j}{u}{-1}{j}))\ge k+4(n-1)$
and the functions $\loi{p_j}{u}{}{j}\loi{q_j}{u}{-1}{j}$
and $w$ differ in at most $4(n-1)$ positions in $D$.
Thus $\pc(\pi_D(\loi{p_j}{u}{}{j}\loi{q_j}{u}{-1}{j}))\le
      \pc(\pi_D(w))+4(n-1)$ and hence $\pc(w)\ge k$.

\subsection{Proof of Lemma~\ref{lem:loop}}

	Notice that the ``if''-direction of statement 1.
	is a special case of 2.
	Let $L = \{(i_1,j_1), \dots, (i_\ell,j_\ell)\}$.
	By \cref{lem:interval-words} we can construct functions $u_1, \dots, u_\ell \in \langle a \rangle^{(\N)}$
	such that $\prod_{k=1}^\ell \loi{p_k}{u}{}{k}\loi{q_k}{u}{-1}{k}$ has periodic complexity at least $2m+1$
	for all $(p_1, \dots, p_\ell) \neq (q_1, \dots, q_\ell) \in \Z_\infty^\ell$.
	For $i \in [1,\ell]$ let $\bar f_i \in \langle a \rangle^{(t^*)}$
	with $\bar f_i(h) = u_i(j)$ if $h = t^j$ for some $j \in \Z$
	and $\bar f_i(h) = 1$ otherwise.
	Then for all $i \in [0,m]$ we define
	\[
		f_i = \prod_{k \in [1,\ell], \, i_k = i} \bar f_k \prod_{k \in [1,\ell], \, j_k = i} \bar f_k^{-1}.
	\]
	Let $h_1, \dots, h_m \in H$ and define $\sigma_i = h_1 \dots h_i$.
	If $h_1 \dots h_m = 1$ and $\sigma_{i_k} = \sigma_{j_k}$ for all $k \in [1,\ell]$ then
	\begin{align*}
		f_0 h_1 f_1 \dots h_m f_m &= \prod_{i=0}^m \lo{\sigma_i}{f_i}
		= \prod_{k=1}^\ell \loi{\sigma_{i_k}}{\bar{f}}{}{k} \loi{\sigma_{j_k}}{\bar{f}}{-1}{k} = 1.
	\end{align*}
	For statement 2.
	let $g_1, \dots, g_m \in \mathsf{P}_{a,t}(G \wr H)$
	and define $\sigma_i = \sigma(g_1 \dots g_i)$ for all $i \in [0,m]$.
	In particular, $\sigma_0=1_H$.
	We have
	\begin{equation}
		\label{eq:ws}
		\tau(f_0 g_1 f_1 \dots g_m f_m) = \loi{\sigma_0}{f}{}{0} \prod_{i=1}^m w_i \loi{\sigma_i}{f}{}{i}
	\end{equation}
	where $w_i = \lo{\sigma_{i-1}}{\tau(g_i)}$ for $i \in [1,m]$.
	Assume that $\sigma_{i_s} \neq \sigma_{j_s}$.
	We will apply on \eqref{eq:ws} the homomorphism
	\[
		\varphi \colon G^{(H)} \to G^{(\Z)}, \quad \varphi(f)(n) = f(\sigma_{i_s} t^n).
	\]
	For each $k \in [1,\ell]$ let $p_k \in \Z$ with
	$\sigma_{i_s} t^{p_k} = \sigma_{i_k}$
	if $\sigma_{i_s}^{-1}\sigma_{i_k}\in \langle t\rangle$
	and $p_k = \infty$ otherwise.
	It satisfies $\varphi(\loi{\sigma_{i_k}}{\bar f}{}{k}) = \loi{p_k}{u}{}{k}$:
	If $\sigma_{i_s}^{-1}\sigma_{i_k}\notin \langle t\rangle$ then $p_k = \infty$ and
	\[
		\varphi(\loi{\sigma_{i_k}}{\bar f}{}{k})(n) = \bar f_k(\sigma_{i_k}^{-1} \sigma_{i_s} t^n) = 1.
	\]
	Otherwise, $p_z \in \Z$ and
	\[
		\varphi(\loi{\sigma_{i_k}}{\bar f}{}{k})(n) = \bar f_k(\sigma_{i_k}^{-1} \sigma_{i_s} t^n)
		= \bar f_k(t^{-p_k} t^n) = u_k(n-p_k) = \loi{p_k}{u}{}{k}(n).
	\]
	Similarly, let $q_k \in \Z$ such that
	$\sigma_{i_1} t^{q_k} = \sigma_{j_k}$
	if $\sigma_{i_s}^{-1}\sigma_{i_k}\in \langle t\rangle$
	and $q_k = \infty$ otherwise;
	it satisfies $\varphi(\loi{\sigma_{j_k}}{\bar f}{-1}{k}) = \loi{q_k}{u}{-1}{k}$.
	Since $\sigma_{i_s} \neq \sigma_{j_s}$ we have $p_s \neq q_s$.
	Therefore $\prod_{k=1}^\ell \loi{p_k}{u}{}{k}\loi{q_k}{u}{-1}{k}$ has periodic complexity
	at least $2m+1$.
	Furthermore, $\varphi(w_i)$ is a basic periodic function for all $i \in [1,\ell]$.
	Let $I$ be the set of indices $i \in [1,m]$ where the value of $\tau(g_i)$
	does not belong to $\langle a \rangle$.
	If $i \in I$ then $\tau(g_i)$ has a period that is not commensurable to $t$
	and hence $|\supp(\varphi(w_i))| \le 1$.
	Let $W = \bigcup_{i \in I} \supp(\varphi(w_i))$, which has size at most $m$.
	If $n \in \Z \setminus W$ then $\varphi(w_i)(n) \in \langle a \rangle$ for all $i \in [1,m]$
	and therefore
	\begin{align*}
		\varphi(\tau(f_0 g_1 f_1 \dots g_m f_m))(n) &=
		\varphi(\loi{\sigma_0}{f}{}{0})(n) \prod_{i=1}^m \varphi(w_i)(n) \varphi(\loi{\sigma_i}{f}{}{i})(n) \\
		&= \prod_{i=0}^m \varphi(\loi{\sigma_i}{f}{}{i})(n) \prod_{i=1}^m \varphi(w_i)(n) \\
		&= \prod_{k=1}^\ell \loi{p_k}{u}{}{k}(n) \loi{q_k}{u}{-1}{k}(n) \prod_{i=1}^m \varphi(w_i)(n).
	\end{align*}
	If $f_0 g_1 f_1 \dots g_m f_m = 1$ then
	$\prod_{k=1}^\ell \loi{p_k}{u}{}{k} \loi{q_k}{u}{-1}{k}$ and $\prod_{i=1}^m \varphi(w_i)^{-1}$
	differ in at most $|W| \le m$ positions.
	Since $\prod_{i=1}^m \varphi(w_i)^{-1}$ has periodic complexity at most $m$,
	the periodic complexity of $\prod_{k=1}^\ell \loi{p_k}{u}{}{k} \loi{q_k}{u}{-1}{k}$ is bounded by $2m$,
	which is a contradiction.

\subsection{Proof of Lemma~\ref{lem:kppm-normal}}
	We proceed in two steps.
	Let $E = e_1 \dots e_n$.
	Suppose that there exists a power $e_k = h_k^{x_k}$
	such that $h_k$ has finite order $q \ge 1$.
	We claim that, if $I$ has a solution $\nu$ then there exists one where $\nu(x_k)$ is bounded by $2q-1$.
	If $\nu$ is any solution of $I$
	we can define a solution $\nu'$ by $\nu'(x) = \nu(x)$ for all $x \neq x_k$
	and $\nu'(x_k) = \nu(x_k) - iq$ where $i \in \N$ is minimal such that $0 \le \nu'(x_k) \le 2q-1$.
	Furthermore the induced factorized walks $\pi_{\nu,E} = \pi_1 \dots \pi_n$
	and $\pi_{\nu',E} = \pi_1 \dots \pi_k' \dots \pi_n$ are identical
	up to the $k$-th subwalks $\pi_k$, $\pi_k'$, which have the same support
	(and the same endpoints).
	Therefore $\nu$ and $\nu'$ satisfy the same interval and disjointness constraints.
	Hence, for $c \in \N$ let us define
	\[
		E_c = e_1 \dots e_{k-1} \underbrace{h_k\cdots h_k}_{\text{$c$ atoms $h_k$}} e_{k+1} \dots e_n.
	\]
	Furthermore, we need to adapt the sets $L$ and $D$.
	Every disjointness constraint in $D$ referring to $e_k$
	must be replaced by $c$ disjointness constraints referring to the $c$ atoms $h_k$.
	Formally we set
	\begin{equation}
	\label{eq:id-dd}
	\begin{aligned}
		L_c &= \{ (\iota_c(i),\iota_c(j)) \mid (i,j) \in L \} \\
		D_c &= \{ (i,j) \mid (\delta_c(i),\delta_c(j)) \in D \}
	\end{aligned}
	\end{equation}
	where the functions $\iota_c$ and $\delta_c$ are defined by
	\begin{equation}
	\label{eq:iota-delta}
		\iota_c(i) = \begin{cases}
		i, & \text{if } i < k, \\
		i+c-1, & \text{if } k \le i,
		\end{cases}
		\quad
		\delta_c(i) = \begin{cases}
		i, & \text{if } i < k, \\
		k, & \text{if } k \le i < k+c, \\
		i-c+1, & \text{if } k+c \le i.
		\end{cases}
	\end{equation}
	It is easy to see that $I$ has a solution $\nu$ with $\nu(x_k) = c$
	if and only if $I[x_k = c] = (E_c,L_c,D_c)$ has a solution.
	We construct the set $\mathcal{I} = \{ I[x_k = c] \mid 0 \le c \le 2q-1 \}$.
	This step reduces the number of powers $h_i^{x_i}$ 
	where $h_i$ has finite order, so we can repeat this construction
	until the instances are torsion-free.
	
	Next, to establish orthogonality, suppose that $E$ contains powers $h_\ell^{x_\ell}$
	and $h_r^{x_r}$ such that $(\ell,r) \in D$.
	If $\ell = r$ then the instance is unsatisfiable and we can return $\mathcal{I} = \emptyset$.
	Now assume that $\ell < r$.
	Since $\langle h_\ell \rangle \cap \langle h_r \rangle \neq \{1\}$
	there exist integers $s > 0$ and $t \neq 0$ such that $h_\ell^s = h_r^t$.
	The idea is that, if the $\ell$-th and the $r$-th subwalk intersect
	then they already intersect in the start or the end area
	of one of the rays of constant length.
	
	Assume that $t > 0$ (the case $t < 0$ is similar).
	For the case that $\nu(x_\ell)$ is bounded by $s$ or $\nu(x_r)$ is bounded by $t$,
	we can construct a finite number of instances $I[x_\ell = c]$, $I[x_r = c]$ as above.
	It remains to consider the case that $\nu(x_\ell) \ge s$ and $\nu(x_r) \ge t$.
	We define the following knapsack expression:
	\[
		E' = 
		e_1 \dots e_{\ell-1} h_\ell^s h_\ell^{y_\ell} e_{\ell+1} \dots
		e_{r-1} h_r^t h_r^{y_r} e_{r+1} \dots e_n.
	\]
	Similar to \eqref{eq:id-dd} and \eqref{eq:iota-delta}, we can define sets $L'$, $D'$
	such that $I' = (E',L',D')$ has a solution if and only if
	$I$ has a solution $\nu$ with $\nu(x_\ell) \ge s$ and $\nu(x_\ell) \ge t$.
	In particular, the set $D'$ relates all $s+1$ atoms in $h_\ell^s h_\ell^{y_\ell}$
	to all $t+1$ atoms in $h_r^t h_r^{y_r}$.
	We claim that we can now omit the disjointness constraint between
	$h_\ell^{y_\ell}$ and $h_r^{y_r}$ in $D'$,
	i.e. $I'$ is equivalent to $I'' = (E',L',D'')$
	where
	\[
		D'' = D' \setminus \{ (\ell+s,r+s+t) \}.
	\]
	Clearly, every solution for $I'$ is a solution for $I''$.
	Conversely, assume that $\nu$ is a solution for $I''$.
	and that the disjointness constraint $(\ell+s,r+s+t)$ is violated,
	i.e. there exist $u \in [0,\nu(y_\ell)-1], v \in [0,\nu(y_r)-1]$
	such that
	\[
		\nu(e_1 \dots e_{\ell-1}) h_\ell^s h_\ell^u = \nu(e_1 \dots e_{r-1}) h_r^t h_r^v.
	\]
	We can choose the pair $(u,v)$ to be minimal with respect
	to the partial order $\preceq$ on $\Z^2$ defined by $(u,v) \preceq (u',v')$
	if there exists $d \ge 0$ such that $(u,v) + d \cdot (s,t) = (u',v')$.
	Since we have
	\[
		\nu(e_1 \dots e_{\ell-1}) h_\ell^s h_\ell^{u-s} = \nu(e_1 \dots e_{r-1}) h_r^t h_r^{v-t}
	\]
	we must have $u < s$ or $v < t$ by minimality of $(u,v)$.
	This contradicts the fact that $D''$ is satisfied by $\nu$.
	
	In conclusion, we construct the set
	\[
		\mathcal{I} = \{ I'' \} \cup \{ I[x_\ell = c] \mid 0 \le c < s \} \cup \{ I[x_r = c] \mid 0 \le c < t \}.
	\]
	Notice that the number of disjointness pairs that violate the orthgonality property
	has decreased in each of these instances.
	Furthermore, the transformation preserves torsion-freeness
	so that we can repeat this process until all instances are orthogonal.

 \subsection{Proof of Lemma~\ref{lem:main-lem}}
        \label{to-wreath:loop-constraints}
We begin with a definition of the set $J\subseteq[0,m]^2$ of loop constraints.
We define $J$ on $\hat E$ so as to express the following conditions
\begin{enumerate}
\item all conditions from $L$, which refer to positions in the prefix $E$ of $\hat E$,
\item $E = 1$, \label{cond:2}
\item for every subexpression $E_{i,c,s} = \hat e_{k+1} \dots \hat e_{k+n+2}$
occurring at position $k$ in $E$:
\begin{enumerate}
	\item $E_{i,c,s} = 1$
	\item $e_1 \dots e_{i-1} = \hat e_{k+1} \dots \hat e_{k+i-1}$
	\item $e_1 \dots e_i = \hat e_{k+1} \dots \hat e_{k+i+2}$.
\end{enumerate}
\end{enumerate}

Before we go on to the proof of \cref{lem:main-lem}, we need a lemma.
\begin{lemma}
	\label{lem:t-com}
	For all valuations $\mu$ and $i \in [1,m]$ we have $\mu(\hat e_i) \in \mathsf{P}_{a,t}$.
\end{lemma}

\begin{proof}
	We can verify $\gamma(\hat e_i) \in GH$ easily from \eqref{eq:eics}.
	If $i \in Q_{\hat E}$ then $\tau(\mu(\hat e_i))$ has a support of size $\le 1$
	and hence it has 1 as a period.
	If $i \in P_{\hat E}$ and $\gamma(\hat e_i) \notin \langle a \rangle H$
	then $\sigma(\gamma(\hat e_i)) = t^{-s} h_{j_k} t^s$
	for some $s \in \N$ and $k \in [1,d]$ with $j_k \in P_{\hat E}$.
	By assumption $h_{j_k}$ is not commensurable to $t$
	and therefore $t^{-s} h_{j_k} t^s$ is not commensurable to $t$ either.
\end{proof}

We are now ready to prove \cref{lem:main-lem}.    
Let $\nu$ be an valuation such that $\nu(E) = 1$
and $\nu$ satisfies all loop and disjointness constraints in $L$ and $D$.
We claim that $(\hat E,J)$ has a solution $\mu$.
Let $\pi_{\nu,E} = \pi_1 \dots \pi_n$ be the induced factorized walk.
We extend $\nu$ to a valuation $\mu$ over all variables in $\hat E$
by assigning to the copied variables the same values as the original variables.
Then for all $i \in [1,n]$, $c \in G$, $s \in \N$ we have
$\sigma(\mu(E_{i,c,s})) = 1$ and 
\[
	\tau(\mu(E_{i,c,s}))(h) = \begin{cases}
		c, & \text{if } ht^{-s} \in \supp(\pi_i), \\
		1, & \text{otherwise}.
	\end{cases}
\]
Since all disjointness constraints in $D$ are satisfied we know that
\[
	\mu(E_{i_k,a,s} \cdot E_{j_k,b,s} \cdot E_{i_k,a^{-1},s} \cdot E_{j_k,b^{-1},s}) = 1
\]
for all $1 \le k \le d$ and $s \in S_k$,
and therefore $\mu(\hat E) = 1$.
Furthermore, $\mu$ satisfies all conditions in $J$.

Next we claim $\loi{t^r\!\!}{f}{}{0} \hat e_1 \loi{t^r\!\!}{f}{}{1} \dots \hat e_m \loi{t^r\!\!}{f}{}{m} = 1$
for some $r \le Nm^2$.
Let $g_i = \mu(\hat e_i)$ for $i \in [1,m]$.
Set $\sigma_i = \sigma(g_1 \dots g_i)$ for all $i \in [0,m]$,
which satisfy $\sigma_i = \sigma_j$ for all $(i,j) \in J$.
For all $r \in \N$ we have
\begin{equation}
	\label{eq:ef}
	\loi{t^r\!\!}{f}{}{0} g_1 \loi{t^r\!\!}{f}{}{1} \dots g_m \loi{t^r\!\!}{f}{}{m}
	= \loi{\sigma_0 t^r\!\!}{f}{}{0} \prod_{i=1}^m w_i \loi{\sigma_i t^r\!\!}{f}{}{i}
\end{equation}
where $w_i = \lo{\sigma_{i-1}}{\tau(g_i)}$ for all $i \in [1,m]$.
It suffices to find a number $r \le [0,Nm^2]$ such that each function
$\loi{\sigma_i t^r\!\!}{f}{}{i}$ commutes with each function $w_k$
since
\begin{align*}
	\prod_{i=0}^m \loi{\sigma_i t^r\!\!}{f}{}{i} \prod_{i=1}^m w_i
	=
	\loi{t^r\!\!}{f}{}{0} \sigma(g_1) \loi{t^r\!\!}{f}{}{1} \dots \sigma(g_m) \loi{t^r\!\!}{f}{}{m} g_1 \dots g_m = 1
\end{align*}
where the last equation uses \Cref{lem:loop} and $g_1 \dots g_m = \mu(\hat E) = 1$.
Define the set
\[
	K = \{ k \in [1,m] \mid \gamma(\hat e_k) \notin \langle a \rangle H \}.
\]
If $k \in [1,m] \setminus K$ then $w_k \in \langle a \rangle^{(H)}$
commutes with all functions $\loi{\sigma_i t^r\!\!}{f}{}{i} \in \langle a \rangle^{(H)}$.
  We call a shift $r\in\N$ \emph{good} if
  \[
	\sigma_i t^{r+j} \notin \supp(w_k) \text{ for all $j \in [0,N-1]$ and $i\in[1,m]$}
      \]
   In other words, if we set
      \[ F_i=[0,N-1],~~~A_i=\{s\in\Z \mid \text{$\sigma_it^s\in\supp(w_k)$ for some $k\in K$}\}\]
      for $i\in[1,m]$, then $r$ is good if and only if $(r+F_i)\cap A_i=\emptyset$ for every $i\in[1,m]$. By \cref{lem:t-com}, all $h,h'\in\supp(w_k)$ with $h\ne h'$ satisfy $h^{-1}h'\notin\langle t\rangle$, which means $|A_i|\le |K|\le m$. Thus, \cref{lem:shifting} tells us that there is a good $r\in[0,Nm^2]$.

Assume that $I = (E,L,D)$ has no solution.
Let $\mu$ be any valuation over the variables of $\hat E$.
Let $g_i = \mu(\hat e_i)$ for $i \in [1,m]$
and $\sigma_i = \sigma(g_1 \dots g_i)$ for all $i \in [0,m]$.
Suppose that $\mu$ does not satisfy the conditions in $J$,
say $\sigma_i \neq \sigma_j$ for some $(i,j) \in J$.
Then by \Cref{lem:loop} and \Cref{lem:t-com}
we know that $f_0 g_1 f_1 \dots g_m f_m \neq 1$.
Furthermore,
since $t^{-r}\sigma_i t^r \neq t^{-r} \sigma_j t^r$
we can also apply \Cref{lem:loop} to the product
$(t^{-r} g_1 t^r) \dots (t^{-r} g_m t^r)$
and we obtain
$f_0 (t^{-r} g_1 t^r) f_1 \dots  (t^{-r} g_m t^r) f_m \neq 1$.
Conjugating with $t^r$ yields
$\loi{t^r\!\!}{f}{}{0} g_1 \loi{t^r\!\!}{f}{}{1} \dots g_m \loi{t^r\!\!}{f}{}{m} \neq 1$.

Now let us assume that $\mu$ satisfies all conditions in $J$,
i.e. $\sigma_i = \sigma_j$ for all $(i,j) \in J$.
Similar to \eqref{eq:ef} we have
\begin{equation}
	\label{eq:mu-tau}
	\loi{t^r\!\!}{f}{}{0} g_1 \loi{t^r\!\!}{f}{}{1} \dots g_m \loi{t^r\!\!}{f}{}{m}
	= \loi{\sigma_0 t^r\!\!}{f}{}{0} \prod_{i=1}^m \lo{\sigma_{i-1}}{\tau(g_i)} \loi{\sigma_i t^r\!\!}{f}{}{i}.
\end{equation}
By definition of $J$, for each $j \in [0,m]$ which occurs in $J$ there exists $i \in [0,n-1]$
with $\sigma_i = \sigma_j$.
Therefore
\[
	F = \bigcup_{i = 0}^m \supp(\loi{\sigma_i t^r \!\!}{f}{}{i}) \subseteq \bigcup_{i=0}^{n-1} \{ \sigma_i t^{r+j}  \mid j \in [0,N-1] \}
\]
Consider the following distance function on $H$:
Define $\|g,h\| = |j|$ if $g^{-1}h = t^j$ for some $k \in \Z$
and otherwise $\|g,h\| = \infty$.
We will prove that there exists a set
$U \subseteq H$ such that $\tau(\mu(\hat E))(h) \neq 1$ for all $h \in U$,
$|U| \ge n+1$ and $\|g,h\| \ge N$ for all $g \neq h \in U$.
Since there are at most $n$ elements in $F$ with pairwise distance $\ge N$
there must be an element $h \in U \setminus F$ satisfying
\[
	\tau(\loi{t^r\!\!}{f}{}{0} g_1 \loi{t^r\!\!}{f}{}{1} \dots g_m \loi{t^r\!\!}{f}{}{m})(h)
	\stackrel{\eqref{eq:mu-tau}}{=} \prod_{i=1}^m \lo{\sigma_{i-1}}{\tau(g_i)}(h)
	= \tau(g_1 \dots g_m)(h) = \tau(\mu(\hat E))(h) \neq 1.
\]
Let us now construct such a set $U$.
Let $\pi_{\mu,E} = \pi_1 \dots \pi_n$ be the induced factorized walk
on $E$.
By Condition~\ref{cond:2} we know that $\pi_1 \dots \pi_n$ must be a loop,
and therefore $\mu(E) = 1$.
Furthermore, $\mu$ satisfies all loop constraints in $L$.
Since $I$ has no solution, $\mu$ must violate a disjointness constraint in $D$.
Recall that $I$ is orthogonalized and therefore
$|\supp(\pi_i) \cap \supp(\pi_j)| \le 1$ for all $(i,j) \in D$.
Let $K$ be the set of indices $k \in [1,d]$ where $\supp(\pi_{i_k}) \cap \supp(\pi_{j_k}) \neq \emptyset$
and let $\supp(\pi_{i_k}) \cap \supp(\pi_{j_k}) = \{p_k\}$.
For all $k \in K$ and $s \in S_k$ we have
\[
	\mu(E_{i_k,a,s} \cdot E_{j_k,b,s} \cdot E_{i_k,a^{-1},s} \cdot E_{j_k,b^{-1},s}) =
	(\big[ p_k t^s \mapsto [a,b] \big],1)
\]
in the semidirect product notation where $[h \mapsto g]$ is the function $H \to G$
mapping $h$ to $g$ and all other elements in $H$ to $1$.
For all $k \notin K$ and $s \in S_k$ we have
\[
	\mu(E_{i_k,a,s} \cdot E_{j_k,b,s} \cdot E_{i_k,a^{-1},s} \cdot E_{j_k,b^{-1},s}) = 1.
\]
This implies
\begin{equation}
	\label{eq:tau-mu}
	\tau(\mu(\hat E)) = \prod_{k \in K} \prod_{s \in S_k} \big[ p_k t^s \mapsto [a,b] \big].
\end{equation}
Let $T_k = \{ p_k t^s \mid s \in S_k \}$ for all $k \in K$.
First notice that
\begin{equation}
	\label{eq:distances}
	(n+d)^{2k} N \le \|p_kt^s,p_kt^{s'}\| \le (n+d)^{2k+1} N
\end{equation}
for all $s,s' \in S_k$ with $s\neq s'$, by definition of $S_k$.
We claim that $|T_k \cap T_{k'}| \le 1$ for all $k, k' \in K$ with $k \neq k'$.
Take $k,k' \in K$ with $k < k'$. Any two elements $g,h \in T_k \cap T_{k'}$ with $g\neq h$ satisfy
\[
	(n+d)^{2k'} N \le \|g,h\| \le (n+d)^{2k+1} N,
\]
which contradicts $k < k'$.
Therefore we can take an arbitrary $k \in K$ and
let $U = T_k \setminus \bigcup_{k' \in K \setminus \{k\}} T_{k'}$.
Then $\tau(\mu(\hat E))(h) \neq 1$ for all $h \in U$ by \eqref{eq:tau-mu} and
\[
	|U| \ge |T_k| - |K| + 1 \ge n + d - |K| + 1 \ge n + 1
\]
Furthermore, by \eqref{eq:distances} any two elements in $U$
have distance at least $N$.

\section{Proofs from Section~\ref{sec:applications}}\label{appendix-applications}

\subsection{The discrete Heisenberg group}\label{sec:appendix-h3}
\newcommand{\h}[3]{\begin{psmallmatrix} 1 & #1 & #3 \\ 0 & 1 & #2 \\ 0 & 0 & 1\end{psmallmatrix}}

Let us show (I), (II), and (III). Recall that
\[ A=\h{1}{0}{0},~~B=\h{0}{1}{0},~~C=\h{0}{0}{1} \]
and note that
\[ \h{a}{b}{c}\h{a'}{b'}{c'}=\h{a+a'}{b+b'}{c'+ab'+c}\]
for any $a,b,c\in\Z$.
It is easy to see that $AC=CA$ and $BC=CB$.  Moreover, one readily
checks that the two maps $\alpha,\beta\colon H_3(\Z)\to\Z$ where
$\alpha$ projects to the top-middle and $\beta$ to the right-middle
entry are homomorphisms.  They satisfy $\alpha(A)=1$,
$\alpha(B)=\alpha(C)=0$ and $\beta(B)=1$, $\beta(A)=\beta(C)=0$.  From
this, it follows directly that (I) and (II) hold: Indeed, if
$A^iC^j=A^{i'}C^{j'}$, then applying $\alpha$ yields $i=i'$ and thus
$C^j=C^{j'}$; since $C$ has infinite order, we obtain $j=j'$. A
similar proof establishes (II). Let us now show (III).

\begin{lemma}
  $A^iB^jA^{-i'}B^{-j'}=C^k$ is equivalent to $i=i'$, $j=j'$, and $k=ij$.
\end{lemma}
\begin{proof}
Note that $A^i=\h{i}{0}{0}$, $B^j=\h{0}{j}{0}$, $A^{-i}=\h{-i}{0}{0}$, and $B^{-j}=\h{0}{-j}{0}$. Therefore,
\begin{equation} A^iB^jA^{-i}B^{-j}=\h{i}{j}{ij}A^{-i}B^{-j}=\h{0}{j}{ij}B^{-j}=\h{0}{0}{ij}. \label{commutators-h3}\end{equation}

Now suppose $A^iB^jA^{-i'}B^{-j'}=C^k$. Applying $\alpha$ yields $i=i'$ and applying $\beta$ yields $j=j'$. Hence, we have $A^iB^jA^{-i}B^{-j}=C^k$. By \cref{commutators-h3}, we get
\[ \h{0}{0}{ij}=C^k=\h{0}{0}{k} \]
and thus $ij=k$. Conversely, if $k=ij$, then \cref{commutators-h3}
shows that $A^iB^jA^{-i}B^{-j}=C^k$.
\end{proof}

\subsection{Solvable Baumslag-Solitar groups}\label{sec:appendix-bs}
\newcommand{\bse}[2]{\begin{psmallmatrix} #1 & #2 \\ 0 & 1 \end{psmallmatrix}}
Recall that for $p,q \in \Z \setminus \{0\}$ the \textit{Baumslag-Solitar group} $\BS(p,q)$ is the group presented by
\[\BS(p,q) := \langle a,t \mid t a^p t^{-1} = a^q  \rangle.\]
In the following we consider Baumslag-Solitar groups of the form $\BS(1,q)$ for $q \geq 2$. These groups are solvable and linear. It is well-known (see, for example, \cite{woess2000random}) that $\BS(1,q)$ is isomorphic to the subgroup $T(q)$ of $\GL(2,\Q)$ consisting of the upper triangular matrices
\[\begin{pmatrix}
q^k & u \\
0 & 1
\end{pmatrix}\]
with $k \in \Z$ and $u \in \Z[\tfrac{1}{q}]$. Here $\Z[\tfrac{1}{q}]$ denotes the set of all rational numbers with finite $q$-ary expansion, hence $\Z[\tfrac{1}{q}]=\{m\cdot q^n \mid m,n\in\Z\}$. We identify $\BS(1,q)$ with this subgroup, so that we obtain:
\begin{equation} a=\bse{1}{1},~~~t=\bse{q}{0}. \label{bs-at}\end{equation}
Observe that given two elements $\bse{q^k}{u},\bse{q^\ell}{v}\in T(q)$, their product is
\[\bse{q^k}{u}\bse{q^\ell}{v}=\bse{q^{k+\ell}}{u+q^k\cdot v}. \]

By Lemma 2.1 in \cite{LohreyZ20} the transformation of an element of
$\BS(1,q)$ given as word over the generators $a,t$ into matrix form
and vice versa can be done in polynomial time ($\mathrm{TC}^0$
even). Thus, for algorithmic purposes, we can represent elements of
$\BS(1,q)$ by matrices of $T(q)$ where the entries are given in
$q$-ary encoding.

In this section, we prove \cref{fo-cplus-bs}. To this end, we use an extension of B\"{u}chi arithmetic $(\Z,+,V_q)$~\cite{buechi1960weak}.
Our extension will have the set $\Z[\tfrac{1}{q}]=\{m\cdot q^n \mid m,n\in\Z\}$ as its domain.
$V_q \colon \Z[\tfrac{1}{q}] \to \Z[\tfrac{1}{q}]$ be the function such that $V_q(x)$ is the largest power of $q$ dividing $x$ for any $x \in \Z[\tfrac{1}{q}]$. Here we say that $a \in \Z[\tfrac{1}{q}]$ \emph{divides} $b \in \Z[\tfrac{1}{q}]$ if there is a $k \in \Z$ such that $a k = b$. Furthermore, for each $\ell \in \Z$ we define the binary predicate $S_{\ell}$ on $\Z[\tfrac{1}{q}]$ such that $x S_\ell y$ is fulfilled if and only if there exist $r \in \Z$ and $s \in \N$ such that $x = q^r$ and $y = q^{r+\ell s}$. Then for any $N \in \N$ we define the structure
\[\mathcal{B}_N := (\Z[\tfrac{1}{q}],+,\geq,0,1,V_q,(S_\ell)_{-N \leq \ell \leq N}).\]

\begin{lemma}\label{lem:buechi}
For each given $N$, the first-order theory of $\mathcal{B}_N$ is decidable.
\end{lemma}
\begin{proof}
We show that $\mathcal{B}_N$ is an automatic structure which implies that $\Th(\mathcal{B}_N)$ is decidable (see \cite{khoussainov1994automatic}). We can write each element of $\Z[\tfrac{1}{q}]$ as $\pm \sum_{i=-r}^{r-1} a_i q^i$ where $r \geq 1$ and $a_i \in [0,q-1]$. This representation is unique if we choose $r$ minimal. We encode such an element with the word
\[\pm \begin{pmatrix} a_{-1} \\ a_0 \end{pmatrix} \begin{pmatrix} a_{-2} \\ a_1 \end{pmatrix} \cdots \begin{pmatrix} a_{-r} \\ a_{r-1} \end{pmatrix}\]
over the alphabet $\{+,-\} \cup [0,q-1]^2$. Then all the predicates of $\mathcal{B}_N$ are clearly regular for each $N \in \N$.
\end{proof}

We will also need some preparatory observations. Note that in $\mathcal{B}_N$ we can define the set of integers. It holds that $x \in \Z$ if and only if $V_q(x) \geq 1$. This means that in the following we can quantify over $\Z$ and therefore also over $\N$. We will make use of the following extension of Lemma~4.5 in \cite{LohreyZ20}:
\begin{lemma}\label{lem:mult}
Given the $q$-ary representation of a number $r \in \Z[\tfrac{1}{q}]$ we can effectively construct a formula over $(\Z[\tfrac{1}{q}],+)$ which expresses $y = r \cdot x$ for $x,y \in \Z[\tfrac{1}{q}]$.
\end{lemma}
\begin{proof}
Let $r = \sum_{-k \leq t \leq \ell} a_t q^t$ with $k,\ell \geq 0$ and $a_t \in [0,q-1]$. We have that $y = r x$ if and only if $q^k y = r^\prime x$ where $r^\prime := \sum_{t = 0}^{k+\ell} a_{t-k} q^t \in \Z$. Since $q^k$ and $r^\prime$ are constant integers, we can use iterated addition to express $q^k y$ and $r^\prime x$ by formulas over $(\Z[\tfrac{1}{q}],+)$.
\end{proof}

We are now prepared to prove \cref{fo-cplus-bs}.
\begin{proof}[Proof of \cref{fo-cplus-bs}]
It remains to show that for each finite subset $F\subseteq\BS(1,q)$,
the structure
$(\BS(1,q),(\xrightarrow{g})_{g\in F},(\tostar[g])_{g\in F})$ can be
interpreted in $\mathcal{B}_N$ for some $N$.  We represent each element
$\bse{q^k}{u}$ of $\BS(1,q)$ by the pair $(q^k, u)$ over $\Z[\tfrac{1}{q}]$.
Moreover, we set $N$ to be the maximal value of $|k|$ for which there is an element
$\bse{q^k}{u}$ in $F$ for some $u\in\Z[\tfrac{1}{q}]$.

We now use the idea of the proof of Theorem 4.1 in \cite{LohreyZ20} to
interpret the structure
$(\BS(1,q),(\xrightarrow{g})_{g\in F},(\tostar[g])_{g\in F})$ in
$\mathcal{B}_N$.

Let us fix an element $g = \begin{psmallmatrix} q^\ell & v \\ 0 & 1 \end{psmallmatrix} \in T(q)$. For all $\begin{psmallmatrix} q^k & u \\ 0 & 1 \end{psmallmatrix}, \begin{psmallmatrix} q^m & w \\ 0 & 1 \end{psmallmatrix} \in T(q)$ we have that
\[\begin{pmatrix} q^k & u \\ 0 & 1 \end{pmatrix} \xrightarrow{g} \begin{pmatrix} q^m & w \\ 0 & 1 \end{pmatrix}\]
is fulfilled if and only if
\[q^m = q^k q^\ell \wedge w = u + q^k v\]
which can be expressed by formulas over $\mathcal{B}_N$ for all $N \in \N$ by \Cref{lem:mult}. To express $\tostar[g]$, we use the following observation:
\begin{equation*}
\begin{split}
\begin{pmatrix} q^k & u \\ 0 & 1 \end{pmatrix} \begin{pmatrix} q^\ell & v \\ 0 & 1 \end{pmatrix}^s & = \begin{pmatrix} q^k & u \\ 0 & 1 \end{pmatrix} \begin{pmatrix}
q^{\ell s} & v + q^\ell v + \dots + q^{(s-1)\ell} v \\ 0 & 1 \end{pmatrix} \\
& = \begin{pmatrix} q^k & u \\ 0 & 1 \end{pmatrix} \begin{pmatrix}
q^{\ell s} & v \frac{q^{\ell s}-1}{q^\ell-1} \\ 0 & 1 \end{pmatrix}
= \begin{pmatrix} q^{k+\ell s} & u + v \frac{q^{k+\ell s}-q^k}{q^\ell-1} \\ 0 & 1 \end{pmatrix}
\end{split}
\end{equation*}
for $\ell \neq 0$ and $s \in \N$. Then for $\ell \neq 0$ and all $\begin{psmallmatrix} q^k & u \\ 0 & 1 \end{psmallmatrix}, \begin{psmallmatrix} q^m & w \\ 0 & 1 \end{psmallmatrix} \in T(q)$ we have that
\[\begin{pmatrix} q^k & u \\ 0 & 1 \end{pmatrix} \tostar[g] \begin{pmatrix} q^m & w \\ 0 & 1 \end{pmatrix}\]
is fulfilled if and only if
\[\exists x \in \Z[\tfrac{1}{q}]\colon \exists s \in \N \colon q^m = q^{k+\ell s} \wedge w = u + v x \wedge (q^\ell-1) x = q^m-q^k\]
where we can quantify $x$ over $\Z[\tfrac{1}{q}]$ since $\frac{q^{\ell s}-1}{q^\ell-1}$ is an integer and therefore $q^k \frac{q^{\ell s}-1}{q^\ell-1} \in \Z[\tfrac{1}{q}]$. By \Cref{lem:mult} we have that $w = u + v x$ and $(q^\ell-1) x = q^m-q^k$ are expressible by formulas over $\mathcal{B}_N$ for all $N \in \N$. Moreover, we can express $\exists s \in \N \colon q^m = q^{k+\ell s}$ by $q^k S_\ell q^m$ with $|\ell|\le N$ and therefore in $\mathcal{B}_N$.

If $\ell = 0$, it holds that $g^s = \begin{psmallmatrix} 1 & s v \\ 0 & 1 \end{psmallmatrix}$. Thus, we have that $\begin{psmallmatrix} q^k & u \\ 0 & 1 \end{psmallmatrix} \tostar[g] \begin{psmallmatrix} q^m & w \\ 0 & 1 \end{psmallmatrix}$ is equivalent to
\[\exists s \in \N \colon w = u + q^k s v \wedge q^m = q^k\]
which holds if and only if
\[\exists t \in \N \colon V_q(t) \geq q^k \wedge w = u + v t \wedge q^m = q^k\]
since we can set $t = q^k s$. Again by \Cref{lem:mult} we can express $w = u + v t$ by a formula over $\mathcal{B}_N$ for all $N \in \N$.
\end{proof}

\section{Exponent equations in Baumslag-Solitar groups}\label{appendix-expeq-bs}
The following unpublished proof is due to Moses Ganardi and Markus
Lohrey~\cite{GanardiLohrey2020}.  With their kind permission, we
include the proof for the convenience of the reader.

\begin{theorem} \label{thm-main-KP}
$\ExpEq(\BS(1,2))$ is undecidable.
\end{theorem}
\begin{proof}
  Consider the function $P\colon (x,y) \mapsto x \cdot 2^y$ on the
  natural numbers. B\"uchi and Senger~\cite[Corollary
  5]{BuchiSenger1988} have shown that the existential fragment of the
  first-order theory of $(\mathbb{N}, +, P)$ is undecidable.  We
  reduce this fragment to $\ExpEq(\BS(1,2))$. For this, it suffices to
  consider an existentially quantified conjunction of formulas of the
  following form: $x \cdot 2^y = z$, $x+y = z$, and $x < y$ (the
  latter allow us to express inequalities and thus negations).  We
  replace each of these formulas by an equivalent exponent equation
  over $\BS(1,2)$.  For this we use the two generators $a$ and $t$ as
  in \cref{bs-at}.
  The formula $x+y=z$ is clearly equivalent to $a^x a^y = a^z$, i.e.,
  $a^x a^y a^{-z} = 1$. The formula $x < y$ is equivalent to
  $a^x a^z a a^{-y}=1$ for some fresh variable $z$. Finally,
  $x \cdot 2^y = z$ is equivalent to $t^y a^x t^{-y} a^{-z}=1$.
\end{proof}

\end{document}